\documentclass[a4paper,reqno,10pt,oneside]{amsart}
\usepackage{amsmath,geometry}
\usepackage{graphicx}
\usepackage{mathrsfs}
\usepackage{amssymb,amsmath, amsfonts, amsthm}
\usepackage{latexsym}
\usepackage{color} 
\usepackage[dvips]{epsfig}
\usepackage{graphicx}

\allowdisplaybreaks[1]

\numberwithin{equation}{section}

\renewcommand{\O}{\operatorname{O}}

\renewcommand{\(}{\left(}
\renewcommand{\)}{\right)}
\renewcommand{\[}{\left[}
\renewcommand{\]}{\right]}

\newtheorem{theorem}{Theorem}[section]
\newtheorem{proposition}[theorem]{Proposition}
\newtheorem{lemma}[theorem]{Lemma}
\newtheorem{remark}[theorem]{Remark}
\newtheorem{definition}[theorem]{Definition}

\newcommand{\zo}{H^1_0(\Omega)}
\newcommand{\zn}{H^1(\mathbb R^N)}

\renewcommand{\d }{\delta }

\renewcommand{\l }{\lambda}

\renewcommand{\O}{\Omega}

\newcommand{\A}{{\mathcal A}}

\newcommand{\U}{\mathcal{U}}
\newcommand{\p}{\mathcal{P}}

\newcommand{\beq}{\begin{equation}}
\newcommand{\eeq}{\end{equation}}
\newcommand{\beqs}{\begin{equation*}}
\newcommand{\eeqs}{\end{equation*}}
\newcommand{\beqn}{\begin{eqnarray}}
\newcommand{\eeqn}{\end{eqnarray}}
\newcommand{\beqns}{\begin{eqnarray*}}
\newcommand{\eeqns}{\end{eqnarray*}}
\newcommand{\bdoc}{\begin{document}}
\newcommand{\edoc}{\end{document}}
\newcommand{\be}{\begin{enumerate}}
\newcommand{\ee}{\end{enumerate}}
\newcommand{\bdescr}{\begin{description}}
\newcommand{\edescr}{\end{description}}
\newcommand{\ba}{\begin{array}}
\newcommand{\ea}{\end{array}}
\newcommand{\intR}{\int_{\mathbb R^N}}
\newcommand{\R}{\mathbb R^N}
\newcommand{\e}{\epsilon}
 \renewcommand{\(}{\left(}
\renewcommand{\)}{\right)}
\renewcommand{\[}{\left[}
\renewcommand{\]}{\right]}
\newenvironment{Proof}{\noindent{\bf Proof}}{\hfill$\Box$\\[2mm]}

\begin{document}

\title[Sign-changing blowing-up sol. for the Brezis--Nirenberg prob. in dim. $N=4,5$]{Sign-changing blowing-up solutions for the Brezis--Nirenberg problem in dimensions four and five}

\author{A. Iacopetti, G. Vaira}

\thanks{Research partially supported by MIUR-PRIN project-201274FYK7\underline\ 005.}
\begin{abstract}
We consider the Brezis-Nirenberg problem:
$$-\Delta u =\lambda u + |u|^{p-1}u\qquad \mbox{in}\,\, \Omega,\quad u=0\,\, \mbox{on}\,\,\ \partial\Omega,$$
where $\Omega$ is a smooth bounded domain in $\R$, $N\geq 3$, $p=\frac{N+2}{N-2}$ and $\lambda>0$.

 In this paper we prove that, if $\Omega$ is symmetric and $N=4,5$, there exists a sign-changing solution whose positive part concentrates and blows-up at the center of symmetry of the domain, while the negative part vanishes, as $\lambda\rightarrow \lambda_1$,  where $\lambda_1=\lambda_1(\Omega)$ denotes the first eigenvalue of $-\Delta$ on $\Omega$, with zero Dirichlet boundary condition.
\end{abstract}

\maketitle

\section{Introduction and statement of the main results}\label{intro}

In this paper we deal with the following problem
\begin{equation}\label{BN}
\left\{
\begin{array}{lr}
-\Delta u = \lambda u + |u|^{p-1}u \qquad \mbox{in}\,\, \Omega\\
u=0,\qquad\qquad\qquad\qquad\mbox{on}\,\, \partial\Omega
\end{array}
\right.
\end{equation}
where $\Omega$ is a bounded smooth domain of $\mathbb R^N$, $N=4,5$, $\lambda>0$, 
$p+1=\frac{2N}{N-2}$ is the critical Sobolev exponent for the embedding of $H^1_0(\Omega)$ into $L^{p+1}(\Omega)$.\\

Problem \eqref{BN} is known as the Brezis--Nirenberg problem since the first existence results for positive solutions of \eqref{BN} were given in their celebrated paper \cite{Brezis}. In particular they showed that the dimension $N$ plays a crucial role in the study of Problem \eqref{BN}. In fact they proved that if $N\geq4$ there exist positive solutions of \eqref{BN} for every $\lambda \in (0,\lambda_1)$, where $\lambda_1=\lambda_1(\Omega)$ is the first eigenvalue of $-\Delta$ on $\Omega$ with zero Dirichlet boundary condition, while if $N=3$ there exists $\lambda_*=\lambda_*(\Omega)>0$ positive solutions exist if $\lambda \in (\lambda_*,\lambda_1)$. When $\Omega=B$ is a ball they also proved that $\lambda_*(B)=\frac{\lambda_1(B)}{4}$ and a positive solution of \eqref{BN} exists if and only if $\lambda \in (\frac{\lambda_1(B)}{4}, \lambda_1(B))$. Moreover, as consequence of the classical Pohozaev's identity positive solutions do not exist if $\lambda\leq 0$ and $\Omega$ is star-shaped.

 Since then several results have been obtained for Problem \eqref{BN}, in particular on the asymptotic analysis of positive solutions, mainly for $N\geq 5$, because also the case $N=4$ presents more difficulties compared to the higher dimensional ones. 

 Concerning the case of sign-changing solutions of \eqref{BN}, several existence results have been obtained if $N\geq 4$. In this case one can get sign-changing solutions for every $\lambda \in (0,\lambda_1(\Omega))$, or even $\lambda > \lambda_1(\Omega)$ (see \cite{Capozzi, Cerami, Castro, Devillanova, Devillanova1, ABP, Clapp, Cerami2, Schechter}).  In particular, Capozzi, Fortunato and Palmieri in \cite{Capozzi} showed that for $N=4$, $\lambda>0$ and $\lambda\not\in \sigma(-\Delta)$ (the spectrum of $-\Delta$ in $H_0^1(\Omega)$) Problem \eqref{BN} has a nontrivial solution. The same holds if $N\geq 5$ for all $\lambda>0$.

The case $N=3$ presents the same difficulties enlightened before for positive solutions and even more. In fact, it is not yet known, when $\Omega=B$ is a ball in $\mathbb{R}^3$, if there are nonradial sign-changing solutions of \eqref{BN} when $\lambda$ is smaller than $\lambda_*(B)=\lambda_1(B)/4$. A partial answer to this question posed by H. Brezis has been given in \cite{Ben}.\\

However, even in the case $N=4,5,6$, some apparently strange phenomenon appears for what concerns radial sign-changing solutions in the ball. Indeed it was first proved by Atkinson, Brezis and Peletier in \cite{ABP1} that for $N=4,5,6$ there exists $\lambda^*=\lambda^*(N)$ such that there are no sign-changing radial solutions of \eqref{BN} for $\lambda\in (0, \lambda^*)$. Later this result was proved in \cite{AY} in a different way.

As it will be clear in the sequel, the nonexistence result of Atkinson, Brezis and Peletier is connected to the asymptotic analysis of low-energy sign-changing solutions of \eqref{BN}. Ben Ayed, El Mehdi and Pacella have investigated the latter question in \cite{Ben},\cite{Ben1}. More precisely, denoting by $\|\cdot\|$ the $H^1_0(\Omega)$-norm and by $S$ the best Sobolev constant for the embedding $H_0^{1}(\Omega)$ into $L^{2*}(\Omega)$, they studied the asymptotic behavior of sign-changing solutions $u_\lambda$ of \eqref{BN} such that $\|u_\lambda\|^2\rightarrow 2S^{N/2}$, as $\lambda\rightarrow 0$ if $N\geq 4$, or $\lambda\rightarrow \bar\lambda$, if $N=3$, where $\bar\lambda$ is the infimum of the values of $\lambda$ for which nodal low-energy solutions exist (see \cite{Ben}). They proved that these solutions blow up at two different points $\bar a_1$, $\bar a_2$, which are the limit of the concentration points $a_{\lambda,1}$, $a_{\lambda,2}$ of the positive and negative part of $u_\lambda$. We point out that they need to assume the extra hypothesis that the concentration speeds of the two concentration points are comparable for $N\geq4$ (see \cite{Ben1}), while in dimension three this was derived without any extra assumption (see \cite{Ben}).

In view of the results of Ben Ayed, El Mehdi and Pacella we get that, for $N\geq 4$, the question of proving the existence of sign-changing low-energy solutions (i.e. such that $\|u_\lambda\|_\Omega^2$ converges to $2S^{N/2}$ as $\lambda \rightarrow 0$) whose positive and negative part concentrate and blow up at the same point, was left open.

In \cite{Iac}, by studying the asymptotic behavior, as $\lambda\rightarrow 0$, of low-energy radial sign-changing solutions of \eqref{BN} in the unit ball of $\R$, for $N\geq7$ (for these dimensions they do exist, as proved by Cerami, Solimini and Struwe in \cite{Cerami2}), it has been proved that the positive and the negative part of such solutions concentrate and blow up at the center of the ball, and their concentration speeds are not comparable.
 Moreover, in the recent paper \cite{IacVair}, it has been proved that for $N\geq7$ these so called ``bubble-tower" solutions for \eqref{BN}, exist, as $\lambda\to 0$, in general bounded domains with some symmetry.

We point out that, in the previous result, the assumption $N\geq7$ on the dimension is not only technically crucial but it also is necessary. In fact, in the recent paper \cite{IacPac}, the authors proved that for the low dimensions $N=4,5,6$, and in general bounded domains, there cannot exist sign-changing ``bubble-tower" solutions for \eqref{BN}, as $\lambda \to 0$.
This result is hence the counterpart, in general bounded domains, of the nonexistence theorem of Atkinson, Brezis and Peletier if we think of sign-changing ``bubble-tower'' solutions as the functions which play, in general bounded domains, the same role as the radial solutions in the case of the ball.\\

In view of all these results it is natural to ask what kind of asymptotic profile we can expect for sign-changing solutions in the low dimensions $N=4,5,6$, as $\lambda$ goes to some strictly positive ``limit'' value. The case of radial sign-changing solutions in the ball, having two nodal regions, has been investigated in \cite{IacPac2}. By studying the associated differential equation, and taking into account the results of \cite{ABP}, \cite{AP2}, the authors prove that if $(u_\lambda)$ is a family of radial sign-changing solutions of \eqref{BN} in the unit ball $B_1$ of $\R$, having two nodal regions, such that $u_\lambda(0)>0$, and denoting by $\bar\lambda=\bar\lambda(N)$ the limit value of the parameter $\lambda$, which arises from the study of the related ordinary differential equation, then:

\begin{enumerate}
\item[(i)] if $N=4,5$, then $\bar\lambda=\lambda_1(B_1)$, where $\lambda_1(B_1)$ is the first eigenvalue of $-\Delta$ in $H_0^1(B_1)$, and  $u_\lambda^+$  concentrates and blows-up at the center of the ball having the limit profile of a ``standard bubble'' in $\R$ (i.e. a solution of the critical problem in $\R$, see \eqref{Udelta}), while $u_\lambda^-$ converges to zero uniformly, as $\lambda \rightarrow \bar\lambda$.
 \item[(ii)] if $N=6$, then $\bar\lambda\in (0,\lambda_1(B_1))$ and $u_\lambda^+$ behaves as in (i) while $u_\lambda^-$ converges to the unique positive radial solution of \eqref{BN} in $B_1$, as $\lambda\rightarrow \bar\lambda$.
\end{enumerate}

The aim of this paper is to show that, in general (symmetric) bounded domains of $\R$, when $N=4,5$, there exist sign-changing solutions of Problem \eqref{BN} having an asymptotic profile, as $\lambda \to \lambda_1(\O)$, which is similar to that of radial ones in the ball. \\ The case $N=6$ is more delicate and at the moment we can only make some conjecture (see Remark \ref{congettura}).

In order to state our results, we denote by $e_1$ the first (positive, $L^2$-normalized) eigenfunction of the Laplace operator with Dirichlet boundary condition, namely $e_1$ solves the problem
\beq\label{auto}
\left\{
\begin{array}{lr}
-\Delta e_1 =\lambda_1 e_1\qquad \mbox{in}\,\,\ \O\\\\
e_1=0\qquad\qquad \mbox{on}\,\,\partial\O,
\end{array}
\right.
\eeq
 and $e_1>0$ in $\O$, $|e_1|_2^2=\int_\O |e_1|^2 \ dx=1$.
We construct solutions $u_\lambda$ of Problem \eqref{BN} which, up to a remainder term, are given by a superposition of a ``standard bubble'' (suitably projected) and the first eigenfunction of the Laplace operator, multiplied by a factor tending to zero, as $\lambda \to \lambda_1$.

More precisely, denoting by $\p$ the projection onto $H_0^1(\O)$ (see \eqref{proiezione}), we get:\\

\begin{theorem}\label{principale4}
Let $N=4$. Assume that $0\in\O$ and that $\O$ is symmetric with respect to $x_1,\ldots,x_4$.\\ Then, for all $\lambda>\lambda_1$ sufficiently close to $\lambda_1$, there exists a sign-changing solution $u_\lambda$ of problem \eqref{BN} of the form
\beq\label{sol4}
u_\lambda(x)=\p\left(\alpha_4 \frac{(\lambda-\lambda_1)e^{-\frac{1}{\lambda-\lambda_1}}s_{1\lambda}}{(\lambda-\lambda_1)^2e^{-\frac{2}{\lambda-\lambda_1}} s_{1\lambda}^2+|x|^2}\right)-e^{-\frac{1}{\lambda-\lambda_1}}\left[(s_{2\lambda}-1)^2+1\right]e_1+\Phi_\lambda
\eeq
where $\alpha_4=2\sqrt 2$, $s_{j\lambda}\to \bar{s}_j>0$, $j=1, 2$ as $\lambda\to\lambda_1^+$ and $\Phi_\lambda\to 0$ in $H_0^1(\Omega)$ as $\lambda\to\lambda_1^+$. Moreover $u_\lambda$ is even with respect to the variables $x_1, \ldots, x_4$.
\end{theorem}

\begin{theorem}\label{principale5}
Let $N=5$. Assume that $0\in\O$ and that $\O$ is symmetric with respect to $x_1,\ldots,x_5$.\\ Then, for all $\lambda<\lambda_1$ sufficiently close to $\lambda_1$, there exists a sign-changing solution $u_\lambda$ of problem \eqref{BN} of the form
\beq\label{sol5}
u_\lambda(x)=\p\left[\alpha_5 \left(\frac{(\lambda_1-\lambda)^{\frac 32}d_{2\lambda}}{(\lambda_1-\lambda)^2 d_{2\lambda}^2+|x|^2}\right)^{\frac 32}\right]-(\lambda_1-\lambda)^{\frac 34}d_{1\lambda}e_1+\Phi_\lambda
\eeq
where $\alpha_5=15\sqrt{15}$, $d_{j\lambda}\to \bar{d}_j>0$, $j=1, 2$ as $\lambda\to\lambda_1^-$ and $\Phi_\lambda\to 0$ in $H_0^1(\Omega)$ as $\lambda\to\lambda_1^-$. Moreover $u_\lambda$ is even with respect to the variables $x_1, \ldots, x_5$.
\end{theorem}
\begin{remark}
We observe that the solutions obtained in Theorem \ref{principale4} and Theorem \ref{principale5} are sign-changing because, in the case $N=4$ they solve Problem \ref{BN} for $\lambda>\lambda_1$ and it is well known that for these values of the parameter $\lambda$ there cannot exist positive (or negative) solutions of Problem \eqref{BN} (see Remark 1.1 in \cite{Brezis}). In the case $N=5$, the sign-changingness of the solution is a consequence of the estimates of the $L^\infty$-norm of the remainder term (see the proof of Theorem \ref{principale5} and Proposition \ref{estimatelinfrem}).
\end{remark}

We point out that since $\lambda_1(\O)$ is reached from above, if $N=4$, while,  it is reached from below, if $N=5$, our results agree with those of \cite{AG} and \cite{GazGrun} for radial sign-changing solutions in the ball.

Moreover, we observe that, thanks to the estimates of the $L^\infty$-norm of the remainder term in compact subsets of $\overline{\Omega}\setminus \{0\}$ (see the proof of Theorem \ref{principale5}, Proposition \ref{estimatelinfrem} and Remark \ref{congABwN4}), the main contribution to the negative part of the solutions obtained in  Theorem \ref{principale4}, Theorem \ref{principale5} is given by the first (normalized, positive) eigenfunction of $-\Delta$ in $H_0^1(\Omega)$, multiplied by a factor tending to zero, as $\lambda \to \lambda_1$.
Hence, this family of solutions verifies, in the more general setting of bounded (symmetric) domains, a conjecture made by Atkinson, Brezis and Peletier in \cite{ABP} for nodal radial solutions in the ball, for $N=4,5$, which states that the negative part of these nodal radial solutions, converges to zero, in compact subsets of $\overline{B_1}\setminus \{0\}$, as the first eigenfunction of $-\Delta$ in the unit ball multiplied by a vanishing factor, as $\lambda\to\lambda_1$.\\

We also observe that the energy (see \eqref{funzionale}) of the solutions obtained in Theorem \ref{principale4}, Theorem \ref{principale5} converges, as $\lambda\to\lambda_1(\O)$, to the ``critical''  energy level $\frac{1}{N} S^{N/2}$ for the Palais-Smale condition (as a consequence of \eqref{energialimite4}, \eqref{energialimite5} and since the $H_0^1$-norm of the remainder term goes to zero).\\

The proof of our results is based on the Lyapunov-Schmidt reduction method which allows us to reduce the problem of finding blowing-up solutions to \eqref{BN} to the problem of finding critical points of a functional (the reduced energy) which depends only on the concentration parameters.\\ We point out that, since we deal with the critical exponent, there are serious difficulties with the standard procedure (see Section 1 of \cite{IacVair}), when trying to look for critical points for the energy functional associated to \eqref{BN}, namely
\begin{equation}\label{funzionale}
J_\lambda(u)=\frac 12 \int_\O |\nabla u|^2\, dx-\frac{1}{p+1}\int_\O |u|^{p+1}\, dx-\frac{\lambda}{2}\int_\O u^2\, dx,\qquad u\in H^1_0(\O).
\end{equation}
In order to overcome these difficulties, for the case $N=5$ we use an idea, introduced in our paper \cite{IacVair}, which is based on the splitting of the remainder term in two parts. Usually the remainder term $\Phi_\lambda$ is found by solving an infinite dimensional problem, called ``the auxiliary equation'', here, we look for a remainder term which is the sum of two remainder terms, of different orders. These two functions are found by solving a system of two equations, which is obtained by splitting the auxiliary equation in an appropriate way. We also point out that, in order to make the finite dimensional reduction method work, we use some techniques which usually belongs to the variational framework. In fact, the standard procedure allows us to get only estimates of the $H_0^1$-norm of the remainder term, but in our case it is necessary to improve these estimates (at least for one of its components) in the $L^\infty$-norm.\\

For the case $N=4$ we use the standard procedure, but it requires finer and different estimates, since they are more delicate in this dimension, and it also requires suitable choices of the parameters $\delta$ and $\tau$.\\

We also observe that the symmetry assumption on the domain $\Omega$ is only made in order to simplify the computations which however, even in the symmetric context, are long and tough. But there is no reason, a priori, for the previous result not to hold in general domains.\\ 

\section{ Notations and some preliminary results}

We introduce the functions
\begin{equation}\label{Udelta}
\U_{\delta}(x)=\alpha_N \frac{\delta^{\frac{N-2}{2}}}{\left(\delta^2+|x|^2\right)^{\frac{N-2}{2}}},\qquad \delta>0,\,\, x\in\mathbb R^N
\end{equation}

with $\alpha_N:=[N(N-2)]^{\frac{N-2}{4}}$. Is is well known (see \cite{Aubin}, \cite{Caffarelli}, \cite{Talenti}) that \eqref{Udelta} are the only radial solutions of the equation
\begin{equation}\label{pb0}
-\Delta u= u^p\qquad \mbox{in}\,\, \R.
\end{equation}
We define $\varphi_\delta$ to be the unique solution to the problem
\begin{equation}\label{pbvarphi}
\left\{
\begin{array}{lr}
\Delta\varphi_\delta=0\qquad\quad \mbox{in}\,\, \Omega\\
\varphi_\delta=\U_{\delta}\qquad\quad\mbox{on}\,\, \partial\Omega,
\end{array}
\right.
\end{equation}
and let
\begin{equation}\label{proiezione}
\p \U_\delta:=\U_\delta-\varphi_\delta
\end{equation}
be the projection of $\U_\delta$ onto $H^1_0(\Omega)$, i.e.
\begin{equation}\label{pbproiezione}
\left\{
\begin{array}{lr}
-\Delta\p\U_\delta=\U_\delta^p\qquad \mbox{in}\,\,\ \Omega\\
\p\U_\delta=0\qquad\qquad\mbox{on}\,\,\,\partial\Omega.
\end{array}
\right.
\end{equation}
 Finally, we introduce the Robin's function of a domain $\Omega$ which is defined as $\tau(x)=H(x, x)$.\\ Here $H(x, y)$, $x, y\in \Omega$, is given as follows: for all $y\in \Omega$, $H(x, y)$ satisfies $$-\Delta H(x, y)=0,\qquad \mbox{in}\,\, \Omega, \qquad H(x, y)=\frac{1}{|x-y|^{N-2}},\qquad x\in\partial\Omega.$$ The function $H$ is nothing but the regular part of the Green function. Indeed, if $G(x, y)$ denotes the Green function of the Laplace operator at the boundary $\partial\Omega$, we have:
$$G(x, y)=\gamma_N\left(\frac{1}{|x-y|^{N-2}}-H(x, y)\right)$$ with $\gamma_N:=\frac{1}{(N-2)\omega_N}$, where $\omega_N$ denotes the surface area of the unit sphere in $\mathbb R^N$.\\\\

It is well-known that the following expansions holds (see \cite{Rey})
\begin{equation}\label{expvarphi}
\varphi_\delta(x)=\alpha_N \delta^{\frac{N-2}{2}}H(0, x)+O(\delta^{\frac{N+2}{2}})\qquad \mbox{as}\,\, \delta\rightarrow 0.
\end{equation}
Moreover, from elliptic estimates it follows that
\begin{equation}\label{stimaproiezione}
0<\varphi_\delta(x)<c\delta^{\frac{N-2}{2}}, \qquad \mbox{in}\,\, \Omega
\end{equation}
and
\begin{equation}\label{stimavarphi}
|\varphi_\delta|_{q, \O}\leq C \delta^{\frac{N-2}{2}},
\qquad q\in \left(\frac{p+1}{2}, p+1\right]
\end{equation}
see for instance \cite{Vaira} and references therein.
\\\\

 In what follows we let $$(u, v):=\int_\O \nabla u \cdot \nabla v\, dx,\qquad \|u\|:=\(\int_\O |\nabla u|^2\, dx\)^{\frac 12}$$ as the inner product in $H^1_0(\O)$ and its corresponding norm while we denote by
$(\cdot, \cdot)_{\zn}$ and by $\|\cdot\|_{\zn}$ the scalar product and the standard norm in $\zn$.
Moreover we denote by $$|u|_{r}:=\(\int_\O|u|^r\, dx\)^{\frac{1}{r}}$$ the $L^r(\O)$-standard norm for any $r\in [1, +\infty]$. When $A\neq\O$ is any Lebesgue measurable subset of $\R$, or, when $A=\O$ and we need to specify the domain of integration, we use the alternative notations $\|u\|_A$, $|u|_{r, A}$.\\\\

From now on we assume that $\O$ is a bounded open set with smooth boundary of $\R$, symmetric with respect to $x_1, \ldots, x_N$ and which contains the origin.\\\\ We define then $$H_{sim}:=\left\{u\in H^1_0(\O)\,\,:\,\, u \,\, \mbox{is symmetric with respect to}\,\, x_1, \ldots, x_N\right\},$$ and for $q\in [1, +\infty]$ $$L^q_{sim}:=\left\{u\in L^q(\O)\,\,:\,\,\ u \,\, \mbox{is symmetric with respect to}\,\, x_1, \ldots, x_N\right\}.$$

\section{Setting of the problem}
Let $i^*:L^{\frac{2N}{N+2}}_{sim}\rightarrow H_{sim}$ be the adjoint operator of the embedding $i:H_{sim}(\O)\rightarrow L^{\frac{2N}{N-2}}_{sim}$, namely if $v\in L^{\frac{2N}{N+2}}_{sim}$ then $u=i^*(v)$ in $H_{sim}$ is the unique solution of the equation $$-\Delta u=v \qquad \mbox{in}\,\, \O\qquad u=0\qquad \mbox{on}\,\, \partial\O.$$
By the continuity of $i$ it follows that
\beq\label{stimaistar}
\|i^*(v)\|\leq C |v|_{\frac{2N}{N+2}} \qquad \forall v\in L^{\frac{2N}{N+2}}_{sim}
\eeq
for some positive constant $C$ which depends only on $N$.\\
Hence, we can rewrite problem \eqref{BN} in the following way

\beq\label{pbriscritto}
\left\{
\ba{lr}
u=i^*\[f(u)+ \lambda u\]\\
u\in H_{sim}
\ea
\right.
\eeq
where $f(s)=|s|^{p-1}s$, $p=\frac{N+2}{N-2}$.\\\\

Let $Z_\d$ the following function
\begin{equation}\label{Zdelta}
Z_\d(x):=\partial_{\delta} \U_{\delta}(x)=\alpha_N\frac{N-2}{2}\delta^{\frac{N-4}{2}}\frac{|x|^2-\delta^2}{\left(\delta^2+|x|^2\right)^{\frac N 2}}.
\end{equation}
We remark that the function $Z_\d$  solves the problem (see \cite{Bianchi})
\begin{equation}\label{linearizzatopb0}
-\Delta z= p|\U_\delta|^{p-1}z,\qquad \mbox{in}\,\, \R.
\end{equation}
Let $\p Z_\d$ the projection of $Z_\d$ onto $H^1_0(\O)$. Elliptic estimates give

\begin{equation}\label{stimaproiezionederivata}
\p Z_\d(x)=Z_\d(x)-\underbrace{\alpha_N\frac{N-2}{2}\delta^{\frac{N-4}{2}}H(0, x)+O(\delta^{\frac{N}{2}})}_{:=\psi_\delta(x)}
\end{equation}
uniformly in $\Omega$.\\
 We next describe the shape of the solution we are looking for.\\\\ Let $\delta$, $\tau$ be positive parameters defined in the following way: for $N=4$ we let
 \beq\label{deltatau4}
\delta= \e e^{-\frac{1}{\e}}s_1 \qquad \qquad \tau=  e^{-\frac{1}{\e}}g(s_2)\qquad \mbox{with} \,\,\, \lambda-\lambda_1=\e, \quad g(s_2)=(s_2-1)^2+1, \quad  s_j>0.
\eeq
Instead, for $N=5$ we let
\begin{equation}\label{deltatau5}
\tau=\e^{\frac 34}d_1\qquad\qquad  \delta=\e^{\frac 32}d_2,\qquad \mbox{with}\,\,\ \lambda_1-\lambda=\e,\quad d_j>0.
\end{equation}

Fix a small $\eta>0$ and assume that
\begin{equation}\label{limdjfj}
\eta<d_j, s_j<\frac{1}{\eta}\qquad \mbox{for}\,\,\, j=1, 2.
\end{equation}

We look for an approximate solution to problem \eqref{pbriscritto} which is of the form
\begin{equation}\label{sol}
u_\lambda(x)=\p\U_{\delta}-\tau e_1+\Phi_\lambda(x)
\end{equation}
where $e_1>0$ is the first eigenfunction of $-\Delta$ corresponding to the first eigenvalue  $\lambda_1$, and the remainder term $\Phi_\lambda$ is a small function which is even with respect to the variables $x_1, \ldots, x_N$.
\\\\
Finally let us recall some useful inequality that we will use in the sequel.
Since these are known results, we omit the proof.

\begin{lemma}\label{lem00n2}
Let $\alpha$ be a positive real number.
If $\alpha \leq 1$ there holds
$$(x+y)^\alpha \leq x^\alpha + y^\alpha,$$
 for all $x, y >0$. If $\alpha \geq 1$ we have
$$(x+y)^\alpha \leq 2^{\alpha-1} (x^\alpha + y^\alpha),$$
 for all $x, y >0$.
 \end{lemma}

\begin{lemma}\label{lem0n2}
Let $q$ be a positive real number. There exists a positive constant $c$, depending only on $q$, such that for any $a, b\in \mathbb R$
\begin{equation}\label{lem1n2el}
||a+b|^q-|a|^q| \leq
\begin{cases}
c \min\{|b|^q, |a|^{q-1}|b|\} &\ \hbox{if}\ 0<q<1,\\
c (|a|^{q-1}|b|+|b|^q) & \ \hbox{if}\ 1\leq q\leq 2.
\end{cases}
\end{equation}
Moreover if $q>2$ then
\begin{equation}\label{lem1n2e2}
\left||a+b|^q-|a|^q-q |a|^{q-2}ab\right|\leq c\left(|a|^{q-2}|b|^2+|b|^q\right).
\end{equation}
\end{lemma}
\noindent Recalling that $f(s)=|s|^{p-1}s$, where $p=\frac{N+2}{N-2}$, we have:
\begin{lemma}\label{lem1n2}
Let $N<6$. There exists a positive constant $c$, depending only on $p$, such that for any $a,b \in \mathbb R$
\begin{equation}\label{teslem1n2}
|f(a+b)-f(a)-f^\prime(a)b| \leq c(|a|^{p-2}|b|^2+ |b|^{p}),
\end{equation}
and
\begin{equation}\label{teslem12n2}
|f(a+b)-f(a)| \leq c(|a|^{p-1}|b|+ |b|^{p}+|a|^{p-2}|b|^2).
\end{equation}
\end{lemma}

\begin{lemma}\label{lem0n3}
Let $N<6$. There exists a positive constant $c$ depending only on $p$ such that for any $a, b_1, b_2\in \mathbb R$
we get
\begin{equation}
\left|f(a+b_1)-f(a+b_2)-f'(a)(b_1-b_2)\right|\leq c\left(|a|^{p-2}|b_2-b_1|+|b_1|^{p-1}+|b_2|^{p-1}\right)|b_1-b_2|.
\end{equation}\\
\end{lemma}

\subsection{Scheme of the reduction}
Let us consider $$\mathcal{K}_1:={\rm span} \left\{e_1\right\}\subset H_{sim}; \qquad \mathcal{K}:={\rm span} \left\{\p Z_\d, e_1\right\}\subset H_{sim}$$ and $$\mathcal{K}_1^\bot:=\left\{\phi\in H_{sim}\,:\,\, ( \phi, e_1)_{\zo}=0\right\}; \quad \mathcal{K}^\bot:=\left\{\phi\in H_{sim}\,:\,\, ( \phi, e_1)_{\zo}=0, \ ( \phi, \p Z_\d)_{\zo}=0  \right\}.$$
Let $\Pi_1:H_{sim}\rightarrow \mathcal{K}_1$, $\Pi:H_{sim}\rightarrow \mathcal{K}$ and $\Pi_1^\bot: H_{sim}\rightarrow \mathcal{K}_1^\bot$, $\Pi^\bot: H_{sim}\rightarrow \mathcal{K}^\bot$ be the projections onto $\mathcal {K}_1$,  $\mathcal K$ and $\mathcal {K}_1^\bot$, $\mathcal K^\bot$, respectively.\\
We set
\begin{equation}\label{V}
V_{\lambda}(x):=\p\U_{\delta}(x)-\tau e_1(x).
\end{equation}
We remark that $V_\lambda(x)=V_\lambda({\bar s},x)$ for $N=4$ and $V_\lambda(x)=V_\lambda(\bar d,x)$ for $N=5$ where
$\bar s:=(s_1, s_2)\in \mathbb R^2_+$ and $\bar{d}:=(d_1, d_2)\in \mathbb R^2_+$.
\\\\
In order to solve Problem \eqref{BN} we will solve the couple of equations
\begin{equation}\label{aux}
\Pi^\bot\left\{V_\lambda+\Phi_\lambda-i^*\left[f(V_\lambda+\Phi_\lambda)+\lambda (V_\lambda+\Phi_\lambda)\right]\right\}=0,
\end{equation}
\begin{equation}\label{bif}
\Pi\left\{V_\lambda+\Phi_\lambda-i^*\left[f(V_\lambda+\Phi_\lambda)+\lambda(V_\lambda+\Phi_\lambda)\right]\right\}=0.
\end{equation}

Given $\bar s$ and $\bar{d}$ satisfying condition \eqref{limdjfj}, one has to solve first the equation \eqref{aux} in $\Phi_\lambda$ which is the lower order term in the description of the ansatz and then solve equation \eqref{bif}.\\\\
We recall finally the definition of stable critical point that we will use in the sequel.\\
\begin{definition}\label{def}
Let $h:\mathcal D \to \mathbb R$ be a $C^1-$ function where $\mathcal D\subset \mathbb R^m$ is an open set. We say that $x_0$ is a stable critical point if $$\nabla h(x_0)=0$$ and there exists a neighborhood $U$ of $x_0$ such that $$\nabla h(x) \neq 0 \qquad \forall\,\, x\in\partial U $$ $$\nabla h(x)=0,\,\, x\in U \Longrightarrow h(x)=h(x_0)$$ and $$deg (\nabla h, U, 0)\neq 0$$ where $deg$ denotes Brouwer degree.
\end{definition}
We remark that any non-degenerate critical point of $h$ is a stable critical point in the sense of Definition \ref{def}. \\ Moreover it is easy to see that if $x_0$ is a minimum or a maximum point of $h$ (not necessarily non-degenerate) then $x_0$ is a stable critical point in according to Definition \ref{def}.\\\\\\
{\bf The case $N=4$.}
For $N=4$ we write \eqref{aux} as
\beq\label{eqaux4}
\mathcal R_\lambda +\mathcal L(\Phi_\lambda)+\mathcal N(\Phi_\lambda)=0
 \eeq
 where
 \beq\label{rl4}
 \mathcal R_\lambda:=\Pi^\bot\left\{V_\lambda-i^*\left[f(V_\lambda)+\lambda V_\lambda\right]\right\}
 \eeq
 \beq\label{ll4}\mathcal L(\Phi_\lambda):=\Pi^\bot\left\{\Phi_\lambda-i^*\left[f'(\U_\d)\Phi_\lambda+\lambda\Phi_\lambda\right]\right\}
 \eeq
 and
 \beq\label{nl4}
 \mathcal N(\Phi_\lambda):=\Pi^\bot\left\{-i^*\left[f(V_\lambda+\Phi_\lambda)-f(V_\lambda)-f'(\U_\d)\Phi_\lambda\right]\right\}.
 \eeq
{\bf The case $N=5$.}
As anticipated in the introduction, in the case $N=5$ we split the remainder term as $$\Phi_\lambda=\phi_1+\phi_2$$ with $$\|\phi_2\|=o(\|\phi_1\|).$$
We start with solving the auxiliary equation \eqref{aux}.

 We write \eqref{aux} as
\begin{equation}\label{star}
\mathcal{R}_1+\mathcal{R}_2+\mathcal{L}_1(\phi_1)+\mathcal{L}_2(\phi_2)+\mathcal{N}_1(\phi_1)+\mathcal{N}_2(\phi_1,\phi_2)=0,
\end{equation}
 where


\begin{equation}\label{defelementi2}
\mathcal R_1:=\Pi^\bot_1\left\{-\tau e_1-i^*\left[-\lambda\tau e_1\right]\right\},
\end{equation}
\beq\label{defelementi3}
\mathcal R_2:=\Pi^\bot\left\{\p\U_{\delta}-i^*\left[\lambda\p\U_\delta+f(\p\U_\delta-\tau e_1)\right]\right\}.
\eeq
\begin{equation}\label{defelementi1phi1}
\mathcal L_1(\phi_1):=\Pi^\bot_1\left\{\phi_1-i^*\left[\lambda_1\phi_1\right]\right\}
\end{equation}
\begin{equation}\label{defelementi1phi2}
\mathcal L_2(\phi_2):=\Pi^\bot\left\{\phi_2-i^*\left[f'(\U_\delta)\phi_2+\lambda \phi_2\right]\right\}
\end{equation}
\beq\label{defelementi4}
\mathcal{N}_1(\phi_1):=\Pi^{\perp}_1\{  - i^*[  f(-\tau e_1+\phi_1) - (\lambda_1-\lambda) \phi_1]\},
\eeq and
\beq\label{defelementi5}
\mathcal{N}_2(\phi_1,\phi_2):= \displaystyle   \Pi^{\perp}\{   -i^*[f(V_\lambda+\phi_1+\phi_2)- f^\prime(\U_\delta) \phi_2 -f(-\tau e_1+\phi_1)-f(\p\U_{\delta}-\tau e_1) ]\}.
\eeq

In order to solve equation \eqref{aux} we solve the following system of equations



%
\begin{equation}\label{sistemaaux}
\left\{
\begin{array}{lr}
\mathcal{R}_1+\mathcal{L}_1(\phi_1)+\mathcal{N}_1(\phi_1)=0,\\\\
\mathcal{R}_2+\mathcal{L}_2(\phi_2)+\mathcal{N}_2(\phi_1, \phi_2)=0,
\end{array}
\right.
\end{equation}

We remark that it is not restrictive to consider $\mathcal R_1, \mathcal L_1(\phi_1), \mathcal N_1(\phi_1)\in \mathcal K_1^\bot$ since only $\delta_1$ appears.\\\\
It is clear that a solution of \eqref{sistemaaux} gives a solution of \eqref{star}.\\
Hence we solve the first equation in \eqref{sistemaaux} finding a solution $\bar{\phi}_1$ and after that we solve the second equation in \eqref{sistemaaux} finding also $\bar{\phi}_2$ with $\phi_1=\bar\phi_1$. \\\\

\section{A linear problem}
Let us consider the linear operator $ L: \mathcal K^\bot \to \mathcal K^\bot$ such that
\begin{equation}\label{linearegen}
L(\phi)=\Pi^\bot\left\{\phi-i^*\left[f'(U_\delta)\phi+\lambda\phi\right]\right\}.
 \end{equation}
 Next results states the invertibility of $L$ and provides a uniform estimate on $L^{-1}$.\\
\begin{proposition}\label{p:inv4}
Let $N=4, 5$ and $\d$ as in \eqref{deltatau4} or \eqref{deltatau5}.\\ Then, for any small $\eta>0$, there exists $C\equiv C(\eta)>0$ such that for all $\lambda$ sufficiently close to $\lambda_1$, for any real number $s_1\in (\eta,\frac 1\eta)$ (or $d_1\in (\eta, \frac 1\eta)$) and for any $\phi\in\mathcal K^\bot$ it holds $$\|L(\phi)\|\geq C\|\phi\|.$$
\end{proposition}
\begin{Proof}
We argue by contradiction.  Assume that there exists a small $\eta>0$, a sequence $(\lambda_n)_n$ converging to $\lambda_1$, a sequence of real numbers $(s_n)_n\subset (\eta, \frac 1\eta)$ (or $(d_n)_n\subset (\eta, \frac 1 \eta))$ and a sequence of functions $(\phi_n)_n\subset H^1_0(\O)$ such that for all $n \in \mathbb{N}$
\beq\label{phin4}
\phi_n\in \mathcal{K}^\bot \qquad \mbox{and}\qquad \|\phi_n\|=1
\eeq
and
\beq\label{ln4}
L(\phi_n)=h_n\qquad \mbox{with}\quad\,\, \|h_n\|\rightarrow 0, \ \mbox{as}\  n\to+\infty.
\eeq

Since $h_n\in \mathcal{K}^\bot$ we get that there exist some real numbers $c^n_j$, $j=0, 1$ such that
\beq\label{eq14}
\phi_n-i^*\[f'(\U_{\d_n})\phi_n+\lambda_n\phi_n\]=h_n+w_n\qquad \mbox{in}\,\, \O
\eeq
where $w_n=c_0^n \p Z_{\d_n}+c_{1}^n e_1$.\\\\
{\bf Step 1:} The first aim is to prove that:
\beq\label{normawn}
\lim_{n\rightarrow+\infty}\|w_n\|=0.
\eeq
To this end we multiply \eqref{eq14} by $\p Z_{\d_n}$ and by $e_1$ and we integrate by parts in $\Omega$ deducing that
\begin{eqnarray*}
&&\hskip-1.5cm (\phi_n, \p Z_{\d_n})_{\zo}-\int_\O f'(\U_{\d_n})\phi_n \p Z_{\d_n}\, dx -\lambda_n\int_\O \phi_n \p Z_{\d_n}\, dx\\
&&=(h_n,  \p Z_{\d_n})_{\zo} + c_0^n(\p Z_{\d_n}, \p Z_{\d_n})_{\zo}\, dx +c_1^n( e_1,  \p Z_{\d_n})_{\zo}
\end{eqnarray*}
and
\begin{eqnarray*}
&&\hskip-1.5cm (\phi_n, e_1)_{\zo}-\int_\O f'(\U_{\d_n})\phi_n e_1\, dx -\lambda_n\int_\O \phi_n e_1\, dx\\
&&=( h_n, e_1)_{\zo} + c_0^n( \p Z_{\d_n}, e_1)_{\zo} +c_1^n(e_1, e_1)_{\zo}.
\end{eqnarray*}
We remark that since $\p Z_{\d_n}$ solves \eqref{linearizzatopb0} and $\phi_n\in\mathcal K^\bot$ we have
$$0=(\phi_n, \p Z_{\d_n})_{\zo}=\int_\O f'(\U_{\d_n})\phi_n Z_{\d_n}\, dx$$  and $$(\p Z_{\d_n}, \p Z_{\d_n})_{\zo}=\int_{\O}|\nabla \p Z_{\d_n}|^2\, dx = \int_\O f'(\U_{\d_n})Z_{\d_n}\p Z_{\d_n}\, dx.$$ Moreover since $e_1$ solves \eqref{auto}  $$(\p Z_{\d_n}, e_1)_{\zo}=\int_\O \nabla e_1 \nabla \p Z_{\d_n}\, dx=\lambda_1 \int_\O e_1 \p Z_{\d_n}\, dx.$$  and  (since $e_1\in \mathcal K^\bot$) $$0=(\phi_n, e_1)_{\zo}=\lambda_1\int_\O e_1 \phi_n\, dx.$$ Hence the equations become
\begin{eqnarray*}
c_0^n\underbrace{\int_\O f'(U_{\d_n})Z_{\d_n}\p Z_{\d_n}\, dx}_{(I)} +c_1^n\lambda_1\underbrace{\int_\O e_1 \p Z_{\d_n}\, dx}_{(II)}&=&-\int_\O f'(\U_{\d_n})\phi_n (\p Z_{\d_n}-Z_{\d_n})\, dx \\
&&-\lambda_n\int_\O \phi_n \p Z_{\d_n}\, dx-\underbrace{(h_n, \p Z_{\d_n})_{\zo}}_{\begin{array}{lr}:=0\\ \mbox{since}\,\, h_n\in\mathcal K^\bot\end{array}}
\end{eqnarray*}
and
\begin{eqnarray*}
c_0^n\lambda_1\underbrace{\int_\O e_1 \p Z_{\d_n} \, dx}_{(II)} +c_1^n\lambda_1\underbrace{\int_\O e_1^2\, dx}_{:=D_0>0}&=&-\int_\O f'(\U_{\d_n})\phi_n e_1\, dx -\underbrace{(h_n,  e_1)_{\zo}}_{\begin{array}{lr} :=0 \\ \mbox{since}\,\, h_n\in \mathcal K^\bot\end{array}}
\end{eqnarray*}
By definition of projection we have $\p Z_{\d_n}= Z_{\d_n}-\psi_{\d_n}$, where $\psi_{\d_n}$ is an harmonic function and $\psi_{\d_n}=Z_{\d_n}$ on $\partial\O$. Therefore, by elliptic estimates, it follows that there is a constant $C>0$ depending only on $N$ and $\O$, such that $|\psi_{\d_n}|_{\infty, \O}\leq C \d_n^{\frac{N-4}{2}}$ (see also \eqref{stimaproiezionederivata}).\\ Hence $$\int_\O f'(\U_{\d_n})\p Z_{\d_n} Z_{\d_n}\, dx= \int_\O f'(\U_{\d_n})Z_{\d_n}^2\, dx -\int_\O f'(\U_{\d_n})\psi_{\d_n}Z_{\d_n}\, dx$$ Now  $$\int_\O f'(\U_{\d_n}) Z_{\d_n}^2\, dx = \alpha_N^{p+1} \d_n^{-2}\int_{\mathbb R^N}\frac{(|y|^2-1)^2}{(1+|y|^2)^{N+2}}\, dy + O(\d_n^{N-2})= A\d_n^{-2}+o(1)\quad \mbox{as}\,\, n\to+\infty$$ where $$A:=\alpha_N^{p+1}\int_{\R}\frac{(|y|^2-1)^2}{(1+|y|^2)^{N+2}}\, dy.$$ Moreover
\begin{eqnarray*}
\int_\O f'(\U_{\d_n})\psi_{\d_n}Z_{\d_n}\, dx &=& \d_n^{N-4}\alpha_N^{p} \frac{N-2}{2}H(0,0)\int_{\mathbb R^N}\frac{1-|y|^2}{(1+|y|^2)^{\frac{N+4}{2}}}\, dy + O(\d_n^{N-3})\\
&=& \left\{\begin{array}{lr} A_0+ O(\d_n)\quad\,\, \mbox{if}\,\, N=4\\\\ O(\d_n)\qquad \qquad\mbox{if}\,\, N=5\end{array}\right.
\end{eqnarray*}
where
$$A_0:=\alpha_N^p\frac{N-2}{2}H(0,0)\int_{\mathbb R^N}\frac{1-|y|^2}{(1+|y|^2)^{\frac{N+4}{2}}}.$$
 Therefore $$(I) =\left\{\begin{array}{lr} A\d_n^{-2}-A_0+o(1)\qquad \mbox{if}\,\,N=4\\\\A\d_n^{-2}+o(1)\qquad \qquad\ \mbox{if}\,\, N=5\end{array}\right.$$ as $ n\to+\infty.$ Moreover $$\int_\O e_1 \p Z_{\d_n}\, dx= \int_\O e_1 Z_{\d_n}\, dx -\int_\O e_1 \psi_{\d_n}\, dx$$ and now
\begin{eqnarray*}
\int_\O e_1 Z_{\d_n}\, dx& = & \d_n^{\frac{N-4}{2}}\int_\O e_1 \frac{|x|^2-\d_n^2}{(\d_n^2+|x|^2)^{\frac{N}{2}}}\, dx = \d_n^{\frac{N-4}{2}}\int_\O \frac{e_1}{|x|^{N-2}}\, dx + o(\d_n^{\frac{N-4}{2}})\\
 &=& B \d_n^{\frac{N-4}{2}} + o(\d_n^{\frac{N-4}{2}} )\qquad \mbox{as}\,\, n\to+\infty
\end{eqnarray*}
since $$B=\int_\O e_1\frac{1}{|x|^{N-2}}\, dx \leq |e_1|_{\infty, \O}\int_\O \frac{1}{|x|^{N-2}}\, dx $$ and $\int_\O\frac{1}{|x|^{N-2}}\, dx $ is integrable in $\O \subset \mathbb R^N$ with $N=4, 5$.\\ Moreover
$$\int_\O e_1\psi_{\d_n}\, dx = \d_n^{\frac{N-4}{2}}\underbrace{\alpha_N\frac{N-2}{2}\int_\O e_1H(0, x)\, dx}_{:=B_0}=\left\{\begin{array}{lr} B_0\qquad\,\,\ \mbox{if}\,\,N=4\\ o(1)\qquad \mbox{if}\,\,N=5 \end{array}\right.$$
We then get
$$(II)=\left\{\begin{array}{lr} B-B_0 +o(1)\qquad \mbox{if}\,\,N=4\\ B\d_n^{\frac 12}+o(1)\qquad\,\,\,\ \ \mbox{if}\,\,N=5\end{array}\right.$$
Hence the equations become for $N=4$
\begin{eqnarray*}
c_0^n\left(A-A_0\d_n^2 +o(\d_n^2)\right)+c_1^n\lambda_1\left((B-B_0)\d_n^2+o(\d_n^2)\right)&=&-\underbrace{\d_n^2\int_\O f'(\U_{\d_n})\phi_n (\p Z_{\d_n}-Z_{\d_n})\, dx}_{(III)} \\
&&-\underbrace{\d_n^2\lambda_n\int_\O \phi_n \p Z_{\d_n}\, dx}_{(IV)}
\end{eqnarray*}
and
\begin{equation*}
c_0^n\lambda_1\left(B-B_0+o(1)\right)+c_1^n\lambda_1D_0=-\underbrace{\int_\O f'(\U_{\d_n})\phi_n e_1\, dx }_{(V)}
\end{equation*}
while for $N=5$
\begin{eqnarray*}
c_0^n\left(A +o(\d_n^2)\right)+c_1^n\lambda_1\left(B\d_n^{\frac 52}+o(\d_n^2)\right)&=&-\underbrace{\d_n^2\int_\O f'(\U_{\d_n})\phi_n (\p Z_{\d_n}-Z_{\d_n})\, dx}_{(III)} \\
&&-\underbrace{\d_n^2\lambda_n\int_\O \phi_n \p Z_{\d_n}\, dx}_{(IV)}
\end{eqnarray*}
and
\begin{equation*}
c_0^n\lambda_1\left(B\d_n^{\frac 12}+o(1)\right)+c_1^n\lambda_1D_0=-\underbrace{\int_\O f'(\U_{\d_n})\phi_n e_1\, dx }_{(V)}
\end{equation*}
Now by using \eqref{stimaproiezionederivata} we get that
\begin{eqnarray*}
|(III)|&=& \left|-\d_n^2\int_\O f'(\U_{\d_n})\phi_n ( \p Z_{\d_n}-Z_{\d_n})\, dx\right|\\
&\leq & \d_n^2 |\p Z_{\d_n}-Z_{\d_n}|_{\frac{2N}{N-2},\O}|\phi_n|_{\frac{2N}{N-2}, \O}|f'(\U_{\d_n})|_{\frac{N}{2}, \O}=O(\d_n^2)
\end{eqnarray*}
since  $|\phi_n|_{\frac{2N}{N-2}, \O}\leq C\|\phi_n\|_{\zo}=C$ and $|f'(\U_{\d_n})|_{\frac{N}{2}, \O}$ is bounded.\\

We remark that $$\|\p Z_{\d_n}\|^2=\int_\O f'(\U_{\d_n})\p Z_{\d_n} Z_{\d_n}\, dx\leq \left| \int_\O f'(\U_{\d_n})\p Z_{\d_n} Z_{\d_n}\, dx\right|\leq C\d_n^{-2}.$$ Hence
we get
$$
|(IV)|\leq  C |\phi_n|_{2, \O}\|\p Z_{\d_n}\|\leq C \d_n.$$
Finally
\begin{eqnarray*}
|(V)|&\leq &  |e_1|_{\infty, \O}\int_\O |f'(\U_{\d_n})\phi_n|\, dx\leq  C \left|\U_{\d_n}\right|_{\frac{N+2}{N-2},\O}^{\frac{4}{N+2}}|\phi_n|_{\frac{N+2}{N-2}, \O}\leq C\d_n^{\frac{2(N-2)}{N+2}}=o(1)
\end{eqnarray*}
Then for $N=4$
\begin{equation*}
\left\{
\begin{array}{lr}
c_0^n\left(A-A_0\d_n^2 +o(\d_n^2)\right)+c_1^n\lambda_1\left((B-B_0)\d_n^2+o(\d_n^2)\right)=o(\d_n)\\
\\
c_0^n\lambda_1\left(B-B_0+o(1)\right)+c_1^n\lambda_1D_0=o(1)
\end{array}
\right.
\end{equation*}
while for $N=5$
\begin{equation*}
\left\{
\begin{array}{lr}
c_0^n\left(A +o(\d_n^2)\right)+c_1^n\lambda_1\left(B\d_n^{\frac 52}+o(\d_n^2)\right)=o(\d_n)\\\\
c_0^n\lambda_1\left(B\d_n^{\frac 12}+o(1)\right)+c_1^n\lambda_1D_0=o(1).
\end{array}
\right.
\end{equation*}
In both cases the system is definitely non singular and hence it has a solution $(c_0^n, c_1^n)$ such that $c_j^n\to 0$ as $n\to+\infty$.\\ Moreover $c_0^n=o(\d_n)$. Now we observe that
\begin{eqnarray*}
\|w_n\|_{\zo}^2&=& (\phi_n, w_n)_{\zo}-\int_\O f'(\U_{\d_n})\phi_n w_n\, dx -\lambda_n \int_\O \phi_n w_n\, dx -(h_n, w_n)_{\zo}\\
&=& -c_0^n\int_\O f'(\U_{\d_n})\phi_n \p Z_{\d_n}\, dx -c_1^n \int_\O f'(\U_{\d_n})\phi_n e_1\, dx +(\phi_n, w_n)_{\zo}\\
&&-\lambda_n c_0^n \int_\O \phi_n \p Z_{\d_n}\, dx -\lambda_n c_1^n\int_\O \phi_n e_1\, dx -(h_n, w_n)_{\zo}\\
\end{eqnarray*}
Reasoning as before and using that $c_0^n=o(\d_n), c_1^n =o(1)$ as $n\to+\infty$ we get that $$\|w_n\|^2_{\zo}= o(1)$$ and the thesis easily follows.

{\bf Step 2:} Let  $$\widetilde{\phi}_n(y)=\d_n^{\frac{N-2}{2}}\phi_n(\d_n y).$$ Then $\widetilde{\phi}_n$ solves the problem
\begin{equation}\label{eqrisc}
-\Delta\widetilde{\phi}_n-p\U(y)^{p-1}\widetilde{\phi}_n-\lambda_n \delta_n^2 \widetilde{\phi}_n=\delta_n^{\frac{N+2}{2}}\Delta(h_n(\d_n y)+w_n(\d_n y))\qquad \mbox{in}\quad \frac{\O_n}{\d_n}.
 \end{equation}
 We point out that, since $\|\widetilde{\phi}_n\|_{\frac{\O}{\d_n}}$ is bounded, then $\widetilde{\phi}_n$ converges weakly in $D^{1, 2}(\R)$, up to a subsequence, to some $\phi_0$. This means that $$\int_{\frac{\O}{\d_n}}\nabla \tilde\phi_n\nabla\varphi\, dx \to \int_{\R}\nabla\phi_0\nabla\varphi\, dx \qquad \mbox{as}\,\, n\to+\infty$$ for any $\varphi\in C^\infty_0(\R)$.\\ By multiplying equation \eqref{eqrisc} by $\varphi\in C^\infty_0(\R)$ and integrating we get that $$\int_{\frac{\O}{\d_n}}\nabla\tilde\phi_n\nabla\varphi\, dx -p \int_{\frac{\O}{\d_n}}\U^{p-1}\tilde\phi_n\varphi\, dx -\lambda_n \delta_n^2\int_{\frac{\O}{\d_n}}\tilde\phi_n\varphi\, dx =\delta_n^{\frac{N+2}{2}}\int_{\frac{\O}{\d_n}}\nabla \tilde h_n \nabla \varphi\, dx +\delta_n^{\frac{N+2}{2}}\int_{\frac{\O}{\d_n}}\nabla \tilde w_n\nabla \varphi\, dx $$ where $\tilde h_n(y)= h_n(\delta_n y)$ and $\tilde w_n (y)=w_n(\d_n y)$. So, as $n\to+\infty$, by using also the results of Step 1, we get that $\phi_0$ solves the problem
$$-\Delta\phi_0=p\left|\U(y)\right|^{p-1}\phi_0\qquad \mbox{in} \,\, \R$$ and satisfies the condition $$\int_{\R}\nabla\phi_0\nabla Z\, dx =0$$ and hence $\phi_0\equiv0$.\\\\  Moreover also $\|\phi_n\|_{\zo}$ is bounded and so, up to a subsequence, also $\phi_n$ converges weakly to some $\phi^*$ in $H_0^1(\O)$ and, as before, we get that, as $n\to+\infty$, $\phi^*$  solves
$$-\Delta\phi^*=\lambda_1\phi^*\qquad\mbox{in}\,\,\ \O$$ with the condition $$\int_\O \nabla e_1 \nabla \phi^*\, dx =0.$$ Hence we get that $\phi^*=0$.\\\\

{\bf Step 3:} We claim that $\|\phi_n\|=o(1)$. This immediately gives a contradiction since by assumption $\|\phi_n\|^2=1$. To prove the claim we multiply by $\phi_n$ the equation \eqref{eq14} and we integrate obtaining
\begin{eqnarray*}
\|\phi_n\|^2&=&\int_\O f'(U_{\d_n})\phi_n^2\, dx +\lambda_n\int_\O \phi_n^2\, dx+(h_n, \phi_n)_{H^1_0}\\
&= &  p\int_{\frac{\O}{\d_n}}|\U|^{p-1}\tilde{\phi}_n^2\, dx+\lambda_n\int_\O \phi_n^2\, dx+o(1)\\
\end{eqnarray*}

$\tilde{\phi}_n^2$ is uniformly bounded in $L^{\frac{N}{N-2}}(\R)$ and converges to zero almost everywhere in $\R$. Hence $\tilde{\phi}_n^2$ converges weakly to zero in $L^{\frac{N}{N-2}}(\R)$. Since $\U^{p-1}\in L^{\frac{N}{2}}(\R)$, by definition, we get that $$p\int_{\frac{\O}{\d_n}}|\U|^{p-1}\tilde{\phi}_n^2\, dx=o(1).$$ Moreover $\phi_n$ converges weakly to zero in $H^1_0(\O)$ and then $\phi_n$ converges strongly to zero in $L^2(\O)$ and then $$\lambda_n \int_\O \phi_n^2\, dx = o(1).$$ At the end we get that $\|\phi_n\|=o(1).$
\end{Proof}
\section{The reduction for $N=4$}
Here we estimate the error term $\mathcal R_\lambda$. The following result holds.
\begin{proposition}\label{errore4}
For any $\eta>0$, there exists $\epsilon_0>0$ and $c>0$ such that  for all $\epsilon \in (0,\epsilon_0)$, for all $(s_1,s_2) \in \mathbb R_+^2$ satisfying \eqref{limdjfj}, we have $$\| \mathcal{R}_\lambda\|\leq c \e e^{-\frac{1}{\e}}.$$
\end{proposition}
\begin{proof}
By continuity of $\Pi^\bot$, by definition of $i^*$ and by using \eqref{stimaistar}, we deduce that
\begin{eqnarray}\label{eq1err2N4}
\nonumber
\displaystyle \|\mathcal{R}_\lambda\| &\leq& \displaystyle c |\lambda\p\U_\delta+f(\p\U_\delta-\tau e_1)-f(U_\delta)+\tau (\lambda-\lambda_1) e_1|_{\frac{2N}{N+2}} \\[10pt]
\nonumber
&\leq& \underbrace{\displaystyle c |f(\p\U_\delta-\tau e_1)-f(\p\U_\delta)|_{\frac{2N}{N+2}}}_{(I)}+\underbrace{\displaystyle c |f(\p\U_{\delta})-f(\U_{\delta})|_{\frac{2N}{N+2}}}_{(II)}  \\[10pt]
&& + \underbrace{c\lambda |\p\U_{\delta}|_{\frac{2N}{N+2}}}_{(III)}+ \underbrace{c\tau|\lambda -\lambda_1| |e_1|_{\frac{2N}{N+2}}}_{(IV)}.
\end{eqnarray}
Let us fix $\eta>0$. We begin with estimating $(I)$. 

Thanks to  Lemma \ref{lem1n2}  there exists a positive constant $c$ (depending only on $p$) such that
\begin{equation}\label{eq6err2}
|f(\p\U_\delta-\tau e_1)-f(\p\U_\delta) |\leq c (\p\U_{\delta}^{p-1}\tau e_1 +\tau^pe_1^p +\tau^2\p\U_\d^{p-2}e_1^2).
\end{equation}
Hence, we get
\begin{eqnarray*}
&&\int_{\Omega} |f(\p\U_\delta-\tau e_1)-f(\p\U_\delta) |^{\frac{2N}{N+2}} \ dx\\
& \leq&  c_1 \int_{\Omega}\left[\p\U_{\delta}^{(p-1)(\frac{2N}{N+2})} (\tau e_1)^{\frac{2N}{N+2}} +(\tau e_1)^{p+1}+\tau^{\frac{4N}{N+2}}\p\U_\d^{\frac{2N(6-N)}{(N-2)(N+2)}}e_1^{\frac{4N}{N+2}}\, dx\right]\\
&\leq & c_2\left(\tau^{\frac{2N}{N+2}}|e_1|_\infty^{\frac{2N}{N+2}}\int_{\Omega}\frac{\delta^{\frac{4N}{N+2}}}{(\delta^2+|x|^2)^{\frac{4N}{N+2}}}\, dx+ \tau^{p+1}\int_\Omega e_1^{p+1}  \ dx +\tau^{\frac{4N}{N+2}}|e_1|_\infty^{\frac{4N}{N+2}}\int_\O \frac{\d^{\frac{N(6-N)}{N+2}}}{\left(\d^2+|x|^2\right)^{\frac{N(6-N)}{N+2}}}\, dx\right)\\
&\leq & c_2\left(\tau^{\frac{2N}{N+2}}|e_1|_\infty^{\frac{2N}{N+2}}\int_{\Omega}\frac{\delta^{-\frac{4N}{N+2}}}{\left(1+\left|\frac{x}{\delta}\right|^2\right)^{\frac{4N}{N+2}}}\, dx+ \tau^{p+1}|e_1|_\infty^{p+1} |\Omega| +\tau^{\frac{4N}{N+2}}|e_1|_\infty^{\frac{4N}{N+2}}\d^{\frac{N(6-N)}{N+2}}\int_\O \frac{1}{|x|^{\frac{2N(6-N)}{N+2}}}\, dx\right)\\
&\leq & c_3\left(\tau^{\frac{2N}{N+2}}\delta^{N-\frac{4N}{N+2}} \int_{\Omega/\delta}\frac{1}{\left(1+\left|y\right|^2\right)^{\frac{4N}{N+2}}}\, dy+ \tau^{p+1}+\tau^{\frac{4N}{N+2}}\d^{\frac{N(6-N)}{N+2}}\right)\\
&\leq & c_3\left(\tau^{\frac{2N}{N+2}}\delta^{\frac{N(N-2)}{N+2}} \int_{\R}\frac{1}{\left(1+\left|y\right|^2\right)^{\frac{4N}{N+2}}}\, dy+ O\left(\tau^{\frac{2N}{N+2}}\delta^{\frac{N(N-2)}{N+2}} \int_{1/\delta}^{+\infty}\frac{r^4}{(1+r^2)^{\frac{4N}{N+2}}}\ dr\right)\right.\\
&&\left.+ \tau^{p+1}+\tau^{\frac{4N}{N+2}}\d^{\frac{N(6-N)}{N+2}}\right).\\
\end{eqnarray*}
since $$\int_\O \frac{1}{|x|^{\frac{2N(6-N)}{N+2}}}\, dx<+\infty$$ for $N=4$.
Hence, we have
$$|f(\p\U_\delta-\tau e_1)-f(\p\U_\delta)|_{\frac{2N}{N+2}} \leq c \left(\tau \delta^{\frac{N-2}{2}}+ \tau^p+\tau^2 \d^{\frac{6-N}{2}}\right).$$
Recalling the choice of $\tau$ and $\delta$ (see \eqref{deltatau4}), we deduce that

 \begin{equation}\label{eq5err2N4}
(I) = |f(\p\U_\delta-\tau e_1)-f(\p\U_\delta)|_{\frac{2N}{N+2}} \leq c  \left(\epsilon e^{-\frac{2}{\e}}+\e^3 e^{-\frac{3}{\e}}+\e e^{-\frac 3\e}\right)\leq c \e e^{-\frac{1}{\e}}
\end{equation}
where $c$ depends on $\eta$ (and also on $\Omega$, $N$).

In order to estimate $(II)$, writing $\p\U_\delta=\U_\delta-\varphi_\delta$, using \eqref{stimaproiezione}, \eqref{stimavarphi} and the usual elementary inequalities, we get
\begin{equation}\label{stimautilefuturo1N4}
\begin{array}{lll}
\displaystyle \int_{\Omega}|(\p\U_{\delta})^p -\U_{\delta}^p|^{\frac{2N}{N+2}}\ dx &\leq & \displaystyle c\left( \int_{\Omega}|\U_{\delta}^{p-1}\varphi_{\delta}|^{\frac{2N}{N+2}}\  dx +  \int_{\Omega}|\varphi_{\delta}|^{p+1}\  dx+\int_\O \U_\d^{\frac{2N(6-N)}{(N+2)(N-2)}}\varphi_\d^{\frac{4N}{N+2}}\, dx\right) \\[10pt]
&\leq & \displaystyle c_1\left( \delta^{\frac{N-2}{2}\frac{2N}{N+2}} \int_{\Omega}\left(\frac{\delta^2}{(\delta^2+|x|^2)^2}\right)^{\frac{2N}{N+2}}\  dx + |\varphi_{\delta}|_{p+1, \O}^{p+1}\right.\\
&&\left.\displaystyle+|\varphi_\d|_\infty^{\frac{4N}{N+2}}\d^{\frac{N(6-N)}{N+2}}\int_\O\frac{1}{|x|^{\frac{2N(6-N)}{N+2}}}\, dx\right)\\[10pt]
&\leq & \displaystyle c_2\left( \delta^{\frac{2N(N-2)}{N+2}}\displaystyle\int_{\Omega/\delta}\frac{1}{(1+|y|^2)^{\frac{4N}{N+2}}}\  dy + \delta^{N}+\d^{\frac{2N(N-2)}{N+2}+\frac{N(6-N)}{N+2}}\right).
\end{array}
\end{equation}
We observe that  for $N=4$ we have that $\frac{1}{(1+|y|^2)^{\frac{4N}{N+2}}} \in L^1(\R)$, and hence we get that
$$ \delta^{\frac{2N(N-2)}{N+2}} \int_{\Omega/\delta}\frac{1}{(1+|y|^2)^{\frac{4N}{N+2}}}\  dy =  \delta^{\frac{2N(N-2)}{N+2}} \int_{\R}\frac{1}{(1+|y|^2)^{\frac{4N}{N+2}}}\  dy + o( \delta^{\frac{2N(N-2)}{N+2}}).$$
Hence it follows that
\begin{equation}\label{stimautilefuturo2N4}
\begin{array}{lll}
(II)=\displaystyle \left( \int_{\Omega}|(\p\U_{\delta})^p -\U_{\delta}^p|^{\frac{2N}{N+2}}\ dx \right)^{\frac{N+2}{2N}}
&\leq & \displaystyle c_3 \left( \delta^{N-2} + \delta^{\frac{N+2}{2}}+\d^{N-2+\frac{6-N}{2}}\right)\\[10pt]
&\leq & \displaystyle c_3  \delta^{N-2} \leq c_4 \epsilon e^{-\frac{1}{\e}}.
\end{array}
\end{equation}
for all sufficiently small $\epsilon$.\\
We now estimate $(III)$. Since $\p\U_\delta\leq \U_\delta$ we have:
\begin{equation*}
\begin{array}{lll}
\displaystyle \int_{\Omega} \p\U_{\delta}^{\frac{2N}{N+2}} \ dx &\leq& \displaystyle  \alpha_N^{\frac{2N}{N+2}} \delta^{\frac{N(N-2)}{N+2}} \int_{\Omega}\frac{1}{(\delta^2+|x|^2)^{\frac{N(N-2)}{N+2}}}\  dx\\[12pt]
&\leq& \displaystyle  \alpha_N^{\frac{2N}{N+2}} \delta^{\frac{N(N-2)}{N+2}} \int_{\Omega}\frac{1}{|x|^{\frac{N(N-2)}{N+2}}}\  dx.
\end{array}
\end{equation*}
Since $N=4$ we have that $\frac{1}{|x|^{\frac{N(N-2)}{N+2}}}$ is integrable near the origin and hence we deduce that $$(III) \leq c  \delta^{\frac{N-2}{2}}\leq c \epsilon e^{-\frac{1}{\e}}.$$
Finally $$(IV)\leq c \tau \e \leq c \epsilon e^{-\frac{1}{\e}}.$$

Putting together all these estimates we deduce that there exists a positive constant $c=c(\eta)>0$ and $\epsilon_0=\epsilon_0(\eta)>0$ such that for all $\epsilon \in (0,\epsilon_0)$, for all $(s_1,s_2) \in \mathbb R_+^2$ satisfying \eqref{limdjfj}
$$ \|\mathcal{R}_\lambda\| \leq c \epsilon e^{-\frac{1}{\e}},$$
and the proof is complete.
\end{proof}
\subsection{Solving the equation \eqref{eqaux4}}
We are now in position to find a solution $\Phi_\lambda\in \mathcal K^\bot$ of the equation \eqref{eqaux4}, namely we prove the following result.
\begin{proposition}\label{aux4solve}
Let $N=4$, $\tau$ and $\delta$  as in \eqref{deltatau4}.\\ Then, for any $\eta>0$, there exists $\epsilon_0>0$ and $c>0$ such that for all $\epsilon \in (0,\epsilon_0)$, for all $(s_1,s_2) \in \mathbb R^2_+$ satisfying condition \eqref{limdjfj}, there exists a unique solution $\bar\Phi_\lambda\in\mathcal{K}^\bot$ of the equation \eqref{eqaux4}, such that
\begin{equation}\label{stimaphilambda4}
\|\bar\Phi_\lambda\|\leq c \epsilon e^{-\frac{1}{\e}}.
\end{equation}

Moreover $\bar\Phi_\lambda$ is continuously differentiable with respect to $(s_1, s_2)$.
\end{proposition}
\begin{proof}
Let us fix $\eta>0$ and define the operator $\mathcal{T} :\mathcal{K}^{\perp} \rightarrow \mathcal{K}^{\perp} $ as $$\mathcal{T}(\Phi_\lambda):=- \mathcal{L}^{-1}[\mathcal{N}(\Phi_\lambda)+\mathcal{R}_\lambda].$$ We remark that $\mathcal L$ is invertible by Proposition \ref{p:inv4} ($\mathcal L \equiv L$).\\\\
In order to find a solution of the equation \eqref{eqaux4} we solve the fixed point problem $\mathcal{T}(\Phi_\lambda)=\Phi_\lambda$.  Let us define the proper ball

$$\displaystyle B_{\epsilon}:=\{\Phi_\lambda \in \mathcal{K}^{\perp}: \ \|\Phi_\lambda\| \leq r\ \epsilon e^{-\frac{1}{\e}}\}$$
for $r>0$ sufficiently large.\\
From Proposition \ref{p:inv4}, there exists $\epsilon_0=\epsilon_0(\eta)>0$ and $c=c(\eta)>0$ such that:
\begin{equation}\label{eq1n24}
\|\mathcal{T}(\Phi_\lambda)\| \leq c (\|\mathcal{N}(\Phi_\lambda)\| + \|\mathcal{R}_\lambda\|),
\end{equation}
and
\begin{equation}\label{eq2n24}
\|\mathcal{T}(\Phi_\lambda)-\mathcal{T}(\Psi_\lambda)\| \leq c (\|\mathcal{N}(\Phi_\lambda)-\mathcal{N}(\Psi_\lambda)\|),
\end{equation}
for all $\Phi_\lambda, \Psi_\lambda \in \mathcal{K}^\perp$, for all $(s_1,s_2) \in \mathbb R_+^2$  satisfying \eqref{limdjfj} and for all $\epsilon \in (0,\epsilon_0)$.

We begin with estimating the right hand side of (\ref{eq1n24}).\\ Thanks to Proposition \ref{errore4}  we have that
$$\| \mathcal{R}_\lambda\| \leq c \epsilon e^{-\frac{1}{\e}},$$
for all $\epsilon \in (0,\epsilon_0)$, for all $(s_1,s_2) \in \mathbb R_+^2$ satisfying \eqref{limdjfj}. Thus it remains only to estimate $\|\mathcal{N}_\lambda(\Phi_\lambda)\|$. Thanks to \eqref{stimaistar} and the definition of $\mathcal{N}_\lambda$ we deduce:
\begin{eqnarray}\label{eq3n24}
\nonumber
\|\mathcal{N}(\Phi_\lambda)\| &\leq & c | f(\p\U_\delta - \tau e_1+\Phi_\lambda)-f(\p\U_\d-\tau e_1)-f^\prime(\p\U_\delta-\tau e_1) \Phi_\lambda|_{\frac{2N}{N+2}}\\
\nonumber
&&+|[f'(\p\U_\d-\tau e_1)-f'(\p\U_\d)]\Phi_\lambda|_{\frac{2N}{N+2}}+|[f'(\p\U_\d)-f'(\U_\d)]\Phi_\lambda |_{\frac{2N}{N+2}}\\
\nonumber
&\leq & c|\p\U_\d^{p-2}\Phi_\lambda^2|_{\frac{2N}{N+2}}
+c|\tau^{p-2}e_1^{p-2}\Phi_\lambda^2|_{\frac{2N}{N+2}}+c|\Phi_\lambda^p|_{\frac{2N}{N+2}}
+c|\tau^{p-1}e_1^{p-1}\Phi_\lambda|_{\frac{2N}{N+2}}\\
&&+c|(\tau e_1)^{p-2}\p\U_\d \Phi_\lambda|_{\frac{2N}{N+2}}+c|\varphi_\d^{p-1}\Phi_\lambda|_{\frac{2N}{N+2}}+c|\varphi_\d^{p-2}\U_\d\Phi_\lambda|_{\frac{2N}{N+2}}
\end{eqnarray}

Now
since $p-2=\frac{6-N}{N-2}$, we have $ (\p\U_\delta^{p-2})^{\frac{2N}{N+2}} =\p\U_\delta^{\frac{2N(6-N)}{(N-2)(N+2)}} \leq \U_\delta^{\frac{2N(6-N)}{(N-2)(N+2)}}\leq c  \delta^{-\frac{N(6-N)}{N+2}} $.
Hence we get that
\begin{eqnarray*}
 \left(\int_\Omega \left(\p\U_\delta^{p-2}\Phi_\lambda^2\right)^{\frac{2N}{N+2}} \ dx\right)^{\frac{N+2}{2N}} &\leq& c \left(\delta^{-N\frac{6-N}{N+2}} \int_\Omega \Phi_\lambda^{\frac{4N}{N+2}} \ dx\right)^{\frac{N+2}{2N}}\\
 &\leq & c_1 \delta^{-\frac{6-N}{2}}\left( \int_\Omega \Phi_\lambda^{\frac{4N}{N+2}} \ dx\right)^{\frac{N+2}{2N}}
  \leq c_2 \delta^{-\frac{6-N}{2}}\|\Phi_\lambda\|^{2}.
\end{eqnarray*}
We observe that for $N=4$, and thanks to the choice of $\delta$ we have  $$\delta^{-\frac{6-N}{2}}\|\Phi_\lambda\|^{2}\leq c \epsilon e^{-\frac{1}{\e}} $$
for all sufficiently small $\epsilon$.

The remaining terms of \eqref{eq3n24} are even simpler. In fact, it holds
\begin{eqnarray*}
 \left(\int_\Omega \left(|\tau e_1|^{p-2}\Phi_\lambda^2\right)^{\frac{2N}{N+2}} \ dx\right)^{\frac{N+2}{2N}} &\leq& c \left(\tau^{\frac{2N(6-N)}{(N+2)(N-2)}}|e_1|_\infty^{\frac{2N}{N+2}} \int_\Omega \Phi_\lambda^{\frac{4N}{N+2}} \ dx\right)^{\frac{N+2}{2N}}\\[10pt]
 &\leq& c_1\tau^{\frac{6-N}{N-2}} \left( \int_\Omega \Phi_\lambda^{\frac{4N}{N+2}} \ dx\right)^{\frac{N+2}{2N}}\\[10pt]
  &\leq& c_2\tau^{\frac{6-N}{N-2}}\|\Phi_\lambda\|^{2}.
\end{eqnarray*}
Thanks to choice of $\tau$ and $\delta$ it is clear that $\tau^{\frac{6-N}{N-2}}\|\Phi_\lambda\|^{2}\leq c \epsilon e^{-\frac{1}{\e}}$.

Moreover
\begin{eqnarray*}
 \left(\int_\Omega |\Phi_\lambda|^{p\frac{2N}{N+2}} \ dx\right)^{\frac{N+2}{2N}}
  &\leq& c \|\Phi_\lambda\|^{p},
\end{eqnarray*}
$$|\tau^{p-1} e_1^{p-1}\Phi_\lambda|_{\frac{2N}{N+2}}\leq c\tau^{p-1} \|\Phi_\lambda\|\leq \e e^{-\frac 1\e} \|\Phi_\lambda\|\leq c\epsilon e^{-\frac{1}{\e}},$$
and for the other terms one can reason as before.\\

It remains to prove that $\mathcal{T}:B_{\epsilon} \rightarrow B_{\epsilon} $ is a contraction. It is sufficient to estimate $\|\mathcal{N}(\phi_1)-\mathcal{N}(\phi_2)\|$ for any $\phi_1, \phi_2 \in B_{\epsilon}$. To this end by using Lemma \ref{lem0n3} we get
\begin{eqnarray*}
\|\mathcal{N}(\phi_1)-\mathcal{N}(\phi_2)\| &\leq & c | f(V_\lambda+\phi_1)-f(V_\lambda+\phi_2)- f^\prime(V_\lambda) (\phi_1-\phi_2)|_{\frac{2N}{N+2}}\\
&&+|[f'(V_\lambda)-f'(\p\U_\d)]|\phi_1-\phi_2||_{\frac{2N}{N+2}}+c|[f'(\p \U_\d)-f'(\U_\d)]|\phi_1-\phi_2||_{\frac{2N}{N+2}}\\
&\leq & c|V_\lambda^{p-2}|\phi_1-\phi_2|^2|_{\frac{2N}{N+2}}+c||\phi_1|^{p-1}|\phi_1-\phi_2||_{\frac{2N}{N+2}}\\
&&+c|\phi_2^{p-1}|\phi_1-\phi_2||_{\frac{2N}{N+2}}+c||\p\U_\d|^{p-2}|\tau e_1||\phi_1-\phi_2||_{\frac{2N}{N+2}}   +c||\tau e_1|^{p-1}|\phi_1-\phi_2||_{\frac{2N}{N+2}}\\
&&+c||\U_\d|^{p-2}|\varphi_\d||\phi_1-\phi_2||_{\frac{2N}{N+2}}  +c||\varphi_\d|^{p-1}|\phi_1-\phi_2||_{\frac{2N}{N+2}}\\
&\leq &\epsilon^\alpha \|\phi_1-\phi_2\|
\end{eqnarray*}

for some $\alpha>0$, for all sufficiently small $\epsilon>0$.\\\\ At the end we get that there exists $L\in (0,1)$ such that
$$ \|\mathcal{T}(\phi_2)- \mathcal{T}(\psi_2)\| \leq L \|\phi_1-\phi_2\|.$$ Moreover in a standard way one can prove that the map $\Phi_\lambda$ is differentiable with respect to $(s_1, s_2)$ (see \cite{Ambrosetti}).
The proof is complete.
\end{proof}

\subsection{Estimates for the reduced functional}
We are left now to solve \eqref{bif} for $N=4$. Let $\bar\Phi_\lambda\in \mathcal{K}^\bot$ be the solution found in Proposition \ref{aux4solve}. Hence $V_\lambda+\bar\Phi_\lambda$ is a solution of our original problem \eqref{BN} if we can find $\bar s_\lambda=(\bar s_{1\lambda}, \bar s_{2\lambda})$  which satisfies condition \eqref{limdjfj} and solves equation \eqref{bif}.

 To this end we consider the reduced functional $\tilde J_\lambda:\mathbb R_+^2 \rightarrow \mathbb R$ defined by: $$\tilde J_\lambda(s_1,s_2):=J_\lambda(V_\lambda+\bar\Phi_\lambda),$$ where $J_\lambda$ is the functional defined in \eqref{funzionale}.\\
Our main goal is to show first that solving equation (\ref{bif}) is equivalent to finding critical points $(\bar s_{1,\lambda},\bar s_{2,\lambda})$ of the reduced functional $\tilde J_\lambda(s_1,s_2)$ and then that the reduced functional has a critical point.

\begin{lemma}\label{CP4}
For any small $\eta>0$ there exists $\e_0>0$ such that for all $\lambda \in (\lambda_1, \lambda_1+\e_0)$ if $(\bar{s}_{1,\lambda}, \bar{s}_{2,\lambda})$ is a critical point of $\tilde J_\lambda$ and satisfies \eqref{limdjfj}, then $V_\lambda+\bar\Phi_\lambda$ is a solution of \eqref{BN}.
\end{lemma}
\begin{proof}
Let us fix a small $\eta>0$ and let $\bar\Phi_\lambda$ be the solution of the auxiliary equation \eqref{aux} found in Proposition \ref{aux4solve}. By assumption $(\bar{s}_{1, \lambda}, \bar s_{2,\lambda})$ is a critical point of $\tilde J_\lambda(s_1, s_2)$.  Hence $$\nabla \tilde J_\lambda(\bar s_{1, \lambda}, \bar s_{2, \lambda})=\left(\partial_{s_1} \tilde J_\lambda, \partial_{s_2}\tilde J_\lambda\right)_{|_{(\bar s_{1, \lambda}, \bar s_{2, \lambda})}}=(0, 0).$$
Then for $j=1, 2$
\begin{eqnarray*}
0&=& \partial_{s_j}\tilde J_\lambda = J'_\lambda(V_\lambda +\bar\Phi_\lambda)[\partial_{s_j}V_\lambda+\partial_{s_j}\bar\Phi_\lambda]\\
&=& \left(V_\lambda+\bar\Phi_\lambda-i^*\left[f(V_\lambda+\bar\Phi_\lambda)+\lambda(V_\lambda+\bar\Phi_\lambda)\right], \partial_{s_j}V_\lambda+\partial_{s_j}\bar\Phi_\lambda\right)_{\zo}\\
&\underbrace{=}_{\bar\Phi_\lambda \,\, \mbox{solution of \eqref{aux}}}& \left(c_0^\lambda \p Z_\d+c_1^\lambda e_1, \partial_{s_j}V_\lambda+\partial_{s_j}\bar\Phi_\lambda\right)_{\zo}\\
&=&\left\{\begin{array}{lr}\left(c_0^\lambda \p Z_\d+c_1^\lambda e_1, \partial_{s_1}(\p \U_\d)+\partial_{s_1}\bar\Phi_\lambda\right)_{\zo}\qquad\,\, \mbox{for}\,\, j=1 \\\\
\left(c_0^\lambda \p Z_\d+c_1^\lambda e_1, \partial_{s_2}(-\tau e_1)+\partial_{s_1}\bar\Phi_\lambda\right)_{\zo}\qquad \mbox{for}\,\, j=2\end{array}\right.
\end{eqnarray*}
for some $c_0^\lambda, c_1^\lambda\in\mathbb R$. \\ We shall prove that $c_0^\lambda, c_1^\lambda$ are equal to zero for $\lambda$ near $\lambda_1$.\\ Now we deduce that (recalling the choice of $\tau$) $$\partial_{s_2}(-\tau e_1)=- e^{-\frac 1 \e}g'(s_2) e_1,$$ while (recalling the choice of $\d$) that $$\partial_{s_1}\p\U_\d= \e e^{-\frac 1\e} \p Z_\d.$$
Hence, we get that
\begin{eqnarray}\label{secondaeqc4}
\nonumber
0&=& \left(c_0^\lambda \p Z_\d +c_1^\lambda e_1, \e e^{-\frac 1\e}\p Z_\d+\partial_{s_1}\bar\Phi_\lambda\right)_{\zo}\\
\nonumber
&=&\e e^{-\frac 1\e} c_0^\lambda (\p Z_\d, \p Z_\d)_{\zo}+ c_0^\lambda (\p Z_\d, \partial_{s_1}\bar\Phi_\lambda)_{\zo}+\e e^{-\frac 1\e} c_1^\lambda(e_1, \p Z_\d)_{\zo}\\
&&+c_1^\lambda (e_1, \partial_{s_1}\bar\Phi_\lambda)_{\zo}
\end{eqnarray}
and
\begin{eqnarray}\label{primaeqc4}
\nonumber
0&=& \left(c_0^\lambda \p Z_\d +c_1^\lambda e_1, - e^{-\frac 1\e}g'(s_2) e_1+\partial_{s_2}\bar\Phi_\lambda\right)_{\zo}\\
\nonumber
&=&- e^{-\frac 1\e}g'(s_2) c_0^\lambda (\p Z_\d, e_1)_{\zo}+ c_0^\lambda (\p Z_\d, \partial_{s_2}\bar\Phi_\lambda)_{\zo}- e^{-\frac 1\e}g'(s_2) c_1^\lambda\underbrace{(e_1, e_1)_{\zo}}_{:=D_0>0}\\
&&+c_1^\lambda (e_1, \partial_{s_2}\bar\Phi_\lambda)_{\zo}
\end{eqnarray}

Now:
\begin{eqnarray*}
(\p Z_\d, e_1)_{\zo}&=& \int_\O f'(\U_\d) Z_\d e_1\, dx = p \int_\O \U_\d^{p-1}Z_\d e_1\, dx\\
&=& p\alpha_N^p \int_{\frac{\O}{\d}}\frac{|y|^2-1}{(1+|y|^2)^{\frac{N+4}{2}}}e_1(\d y)\, dy \\
&=& A_0 + o(1)
\end{eqnarray*}
where $$A_0:=p\alpha_N^p e_1(0)\int_{\mathbb R^4}\frac{|y|^2-1}{(1+|y|^2)^{4}}\, dy.$$

Moreover, we point out that, since $\bar\Phi_\lambda\in \mathcal K^\bot$, it holds $$(e_1, \bar\Phi_\lambda)_{\zo}=0;\qquad (\p Z_\d, \bar\Phi_\lambda)_{\zo}=0.$$ These imply $$(e_1, \partial_{s_j}\bar\Phi_\lambda)_{\zo}=-(\partial_{s_j}e_1, \bar\Phi_\lambda)_{\zo}=0$$ and
\begin{equation}\label{numpag20}
(\p Z_\d, \partial_{s_j}\bar\Phi_\lambda)=-(\partial_{s_j}\p Z_\d, \bar \Phi_\lambda)_{\zo}=\left\{\begin{array}{lr} -(\partial_{s_1}\p Z_\d, \bar\Phi_\lambda)_{\zo}\qquad \mbox{for}\,\, j=1\\\\0\qquad\qquad\qquad\qquad\qquad \mbox{for}\,\, j=2\end{array}\right.
\end{equation}
and $$(\partial_{s_1}\p Z_\d, \bar\Phi_\lambda)_{\zo}=O(\|\partial_{s_1}\p Z_\d\|\cdot \|\bar\Phi_\lambda\|)=o(\|\partial_{s_1}\p Z_\d\|).$$

\begin{eqnarray*}
(\p Z_\d, \p Z_\d)_{\zo}&=& \int_\O f'(\U_\d) \p Z_\d Z_\d\, dx \\
&=& \int_\O \frac{\d^2(\d^2-|x|^2)^2}{(\d^2+|x|^2)^6}\, dx -\int_\O \frac{\d^2(\d^2-|x|^2)}{(\d^2+|x|^2)^4}\psi_\d\, dx\\
&=& B\d^{-2}-B_0+O(\d^2)
\end{eqnarray*}
where $$B:=\int_{\mathbb R^4}\frac{(1-|y|^2)^2}{(1+|y|^2)^6}\, dy;\qquad B_0:=\int_{\mathbb R^4}\frac{1-|y|^2}{(1+|y|^2)^4}\, dy$$
Now $$\partial_{s_1} Z_\d=2\e e^{-\frac 1\e}\frac{\delta(-\d^2+3|x|^2)}{(\d^2+|x|^2)^3}$$
and hence
$$\partial_{x_j}(\partial_{s_1}Z_\d)=24\e e^{-\frac 1\e}\frac{\d x_j(\d^2-|x|^2)}{(\d^2+|x|^2)^4}.$$
We have that
$$\|\partial_{s_1}Z_\d\|^2\leq c \d^{-2}.$$

The equations \eqref{primaeqc4}, \eqref{secondaeqc4} become
\begin{equation}\label{sistemacN4}
\left\{
\begin{array}{lr}
c_0^\lambda (A_0+o(1))+ c_1^\lambda D_0=0\\\\
c_0^\lambda(B-B_0\d^2+o(1))+c_1^\lambda (A_0\d^2+o(\d^2) )=0
\end{array}
\right.
\end{equation}
It follows easily that for all $\e$ sufficiently small, meaning for all $\lambda$ sufficiently near to $\lambda_1$, the system \eqref{sistemacN4} has only the trivial solution and the proof is now complete.

\end{proof}

\begin{lemma}
For any $\eta>0$ there exists $\e_0>0$ such that for any $\e\in (0, \e_0)$ it holds $$J_\lambda(V_\lambda +\bar\Phi_\lambda)=J_\lambda (V_\lambda)+ o(\e e^{-\frac 2\e})$$
\end{lemma}
\begin{proof}
Let us fix $\eta>0$. By making some computations it follows that
\begin{eqnarray*}
J_\lambda(V_\lambda+\bar\Phi_\lambda)-J_\lambda(V_\lambda)&=& \frac 12 \int_\O |\nabla \bar\Phi_\lambda|^2\, dx +\int_\O \nabla V_\lambda \nabla\bar\Phi_\lambda\, dx -\frac{\lambda}{2}\int_\O\bar\Phi_\lambda^2\, dx \\\\
&&-\lambda \int_\O \bar\Phi_\lambda V_\lambda\, dx -\int_\O\left[F(V_\lambda+\bar\Phi_\lambda)-F(V_\lambda)\right]\, dx,
\end{eqnarray*}
where we have set $F(s):=\frac{1}{p+1}|s|^{p+1}$.
Now $$\frac 12\int_\O |\nabla \bar\Phi_\lambda|^2\, dx =\frac 12 \|\bar\Phi_\lambda\|^2=o(\e e^{-\frac 2\e})$$ and $$\frac \lambda 2 \int_\O \bar\Phi_\lambda^2\, dx =\frac \lambda 2 |\bar\Phi_\lambda|^2_{2, \O}\leq C\|\bar\Phi_\lambda\|^2=o(\e e^{-\frac 2\e}).$$
Moreover
$$
\int_\O \nabla V_\lambda \nabla \bar\Phi_\lambda\, dx =\int_\O \nabla(\p\U_\d-\tau e_1)\nabla\bar\Phi_\lambda\, dx =
\int_\O \U_\d^p \bar\Phi_\lambda\, dx -\tau\lambda_1\underbrace{\int_\O e_1\bar\Phi_\lambda\, dx}_{:=0\,\,\mbox{since}\,\, \bar\Phi_\lambda\in\mathcal K^\bot}$$ and hence
$$\int_\O \nabla V_\lambda \nabla\bar\Phi_\lambda\, dx -\lambda\int_\O \bar\Phi_\lambda V_\lambda\, dx = \int_\O \U_\d^p\bar\Phi_\lambda\, dx -\lambda\int_\O \bar\Phi_\lambda \p\U_\d\, dx +\lambda\tau \underbrace{\int_\O e_1\bar\Phi_\lambda\, dx}_{:=0\,\,\mbox{since}\,\, \bar\Phi_\lambda\in\mathcal K^\bot}.$$
We estimate the second term and we get that $$\left|\lambda \int_\O \p \U_\d \bar\Phi_\lambda\, dx \right|\leq \left|\lambda\int_\O \bar\Phi_\lambda\U_\d\, dx\right|+\left|\lambda\int_\O \bar\Phi_\lambda\varphi_\d\, dx \right|\leq C\|\bar\Phi_\lambda\||\U_\d|_{2, \O}+ C\|\bar\Phi_\lambda\||\varphi_\d|_{2, \O}=o(\e e^{-\frac 2\e})$$ since $$|\U_\d|_{ 2, \O}\leq c \d\sqrt{\log\frac 1\d}=o(\e^{1-\sigma}e^{-\frac 1\e})$$ for some $\sigma>0$. Finally, for the remaining terms, by using \eqref{lem1n2e2} we deduce that
\begin{eqnarray*}
&&\hskip-1.5cm \left|\int_\O \left[F(V_\lambda+\bar\Phi_\lambda)-F(V_\lambda)-f(V_\lambda)\bar\Phi_\lambda\right]\right|\leq c \left( \int_\O |V_\lambda|^{p-1}\bar\Phi_\lambda^2\, dx +\int_\O |\bar\Phi_\lambda|^{p+1}\, dx\right)\\\\
&&\leq c\left(|\p\U_\d|_{2, \O}|\bar\Phi_\lambda|_{4, \O}+\tau|e_1|_{2, \O}|\bar\Phi_\lambda|_{4, \O}^2+\|\bar\Phi_\lambda\|^4\right)
\leq c_1\left(\d \log\frac 1\delta\|\bar\Phi_\lambda\|^2+\tau \|\bar\Phi_\lambda\|^2+\|\bar\Phi_\lambda\|^4\right)\\\\
&&= o(\e e^{-\frac 2 \e}).
\end{eqnarray*}
At the end by using \eqref{teslem12n2} we get that
\begin{eqnarray*}
&&\left|\left[f(V_\lambda)-f(\U_\d)\right]\bar\Phi_\lambda\, dx \right|\\
&\leq& c_1\left(\int_\O |\U_\d|^{p-2}(\varphi_\d+\tau e_1)^2 |\bar\Phi_\lambda|\, dx+\int_\O |\U_\d|^{p-1}|\varphi_\d+\tau e_1| |\bar\Phi_\lambda|\, dx+\int_\O | \varphi_\d+ \tau e_1|^p|\bar\Phi_\lambda|\, dx\right)\\
&\leq & c_2 \left( |\varphi_\d +\tau e_1|_\infty^2 |\U_\d|_2 |\bar\Phi_\lambda|_4+ |\U_\d|_2|\varphi_\d+\tau e_1|_4|\bar\Phi_\lambda|_4+ |\varphi_\d+\tau e_1|_\infty^3 |\bar\Phi_\lambda|_4 \right)=o(\e e^{-\frac 2\e})
\end{eqnarray*}

\end{proof}
\section{The reduction for $N=5$}
In this section we solve \eqref{sistemaaux}.

\subsection{The auxiliary equation: solution of the system \eqref{sistemaaux}}\label{ausiliaria}
Let us define:
\begin{equation}\label{thetaj}
\theta_1:=\frac{5}{4};\qquad\qquad \theta_2:=3.
\end{equation}

In this section we solve system \eqref{sistemaaux}. More precisely, the aim is to prove the following result.
\begin{proposition}\label{auxsolving}
Let $N=5$, $\tau$ and $\delta$ as in \eqref{deltatau5}. Then, for any $\eta>0$, there exist $\epsilon_0>0$ and $c>0$ such that for all $\epsilon \in (0,\epsilon_0)$, for all $(d_1,d_2) \in \mathbb R_+^2$ satisfying \eqref{limdjfj}, there exists a unique $\bar\phi_1\in\mathcal K_1^\bot$ solution of the first equation of \eqref{sistemaaux} such that $$\|\bar\phi_1\|\leq c \e^{\frac{\theta_1}{2}+\sigma}.$$ Moreover, fixing $\phi_1=\bar\phi_1$, there exists a unique solution $ \bar\phi_2 \in {\mathcal K}^\bot$ of the second equation of \eqref{sistemaaux} such that $$ \|\bar\phi_2\| \leq c \ \epsilon^{\frac{\theta_2}{2}+\sigma},$$
 for some positive real number  $\sigma$ whose choice depends only on $N$.
Furthermore, $\bar\phi_1$ depends only on $d_1$ and it is continuously differentiable with respect to $d_1$, $\bar\phi_2$ is continuously differentiable with respect to $(d_1, d_2)$.
\end{proposition}

In order to prove Proposition \ref{auxsolving} let us first consider the linear operator
$$\mathcal{L}_1: \mathcal{K}_1^\bot\rightarrow \mathcal{K}_1^{\bot}$$ defined as in
\eqref{defelementi1phi1}.\\ 
Next result states the invertibility of the operator $\mathcal{L}_1$ and provides a uniform estimate on the inverse of the operator $\mathcal L_1$.

\begin{proposition}\label{linearprop1}
The linear operator $\mathcal{L}_1: \mathcal{K}_1^\bot\rightarrow \mathcal{K}_1^{\bot}$ is invertible and $\| \mathcal{L}_1^{-1}\|\leq c$ for some constant depending only on $N$ and $\Omega$.
\end{proposition}
\begin{proof}
Let us fix $h\in \mathcal{K}_1^\bot$. We consider the problem
\begin{equation}\label{l1sist}
\left\{
\begin{array}{lr}
-\Delta \phi =\lambda_1 \phi + h \qquad \mbox{in}\,\, \Omega,\\
\phi=0,\qquad\qquad\qquad\mbox{on}\,\, \partial\Omega.
\end{array}
\right.
\end{equation}
Since $h\in \mathcal{K}_1^\bot$ it is well known that \eqref{l1sist} has a solution $\phi \in H_0^1(\Omega)$ (see \cite{AmbrosettiProdi}, Theorem 0.7). Moreover it is elementary to see that the solution is unique in $\mathcal {K}_1^\bot$.

Hence by definition of $\Pi_1^\bot$ and $i^*$ it follows immediately that $\mathcal{L}_1(\phi)=h$ has a unique solution $\bar\phi=\bar\phi(h)\in \mathcal{K}_1^\bot$, and from elliptic estimates we have $\|\bar\phi\| \leq c \|h\|$, which implies the boundedness of $\mathcal{L}_1^{-1}$. The proof is complete.
\end{proof}

For the linear operator $\mathcal{L}_2$ we have the invertibility proved in Proposition \ref{p:inv4}.

At this point the strategy is to solve the first equation of \eqref{sistemaaux} by a fixed point argument, finding a unique $\bar\phi_1$ and then, substituting  $\bar\phi_1$ in the second equation of \eqref{sistemaaux}, we obtain an equation depending only on the variable $\phi_2$. Hence, using again a fixed point argument, we solve the second equation of \eqref{sistemaaux} uniquely.

\subsection{The solution of the first equation of \eqref{sistemaaux}}
The aim is to prove the following proposition.
\begin{proposition}\label{aux1}
Let $N=5$ and $\tau$ as in \eqref{deltatau5}. Then, for any $\eta>0$, there exists $\epsilon_0>0$ and $c>0$ such that for all $\epsilon \in (0,\epsilon_0)$, for all $d_1 \in \mathbb R_+$ satisfying condition \eqref{limdjfj} for $j=1$, there exists a unique solution $\bar\phi_1=\bar\phi_1(d_1)$, $\bar\phi_1\in\mathcal{K}_1^\bot$ of the first equation in \eqref{sistemaaux} which is continuously differentiable with respect to $d_1$ and such that
\begin{equation}\label{stimaphi1}
\|\bar\phi_1\|\leq c \e^{\frac{\theta_1}{2}+\sigma},
\end{equation}
where $\theta_1$ is defined in \eqref{thetaj} and $\sigma$ is some positive real number whose choice depends only on $N$.
\end{proposition}

In order to prove Proposition \ref{aux1} we observe that:

\begin{proposition}\label{errorprop1}
It holds $$\mathcal{R}_1=0.$$ 
\end{proposition}
\begin{proof}
Let us fix $\tau>0$. By linearity we have $\mathcal{R}_1= \tau \Pi^\bot_1\left\{-e_1-i^*\left[-\lambda e_1\right]\right\}$; hence $\mathcal{R}_1=0$ if and only if $-e_1-i^*\left[-\lambda e_1\right]=c e_1$ for some $c \in \mathbb{R}$. This is true, since, by definition of $i^*$ and $e_1$, it holds $-e_1-i^*\left[-\lambda e_1\right]=(-1+\frac{\lambda}{\lambda_1})e_1$. The proof is complete.
\end{proof}

We are ready to prove Proposition \ref{aux1}.\\
\begin{proof}[Proof of Proposition \ref{aux1}.]
Let us fix $\eta>0$ and define $\mathcal{T}_1 :{\mathcal K}_1^\bot \rightarrow {\mathcal K}_1^\bot $ as $$\mathcal{T}_1(\phi_1):=- \mathcal{L}_1^{-1}[\mathcal{N}_1(\phi_1)].$$ Thanks to Proposition \ref{errorprop1} it is clear that solving the first equation of \eqref{sistemaaux} is equivalent to solving the fixed point equation $\mathcal{T}_1(\phi_1)=\phi_1$.\\ Let us define the ball
$$\displaystyle B_{1,\epsilon}:=\{\phi_1 \in \mathcal{K}^{\perp}; \ \|\phi_1\| \leq r \ \epsilon^{\frac{\theta_1}{2}+\sigma} \}\subset \mathcal{K}^\bot$$ with $r>0$ sufficiently large and $\sigma>0$.\\
We want to prove that, for $\e$ small, $\mathcal{T}_1$ is a contraction in the proper ball $\displaystyle B_{1,\epsilon}$, namely we want to prove that, for $\e$ sufficiently small
\begin{enumerate}
\item $\mathcal{T}_1(B_{1,\epsilon}) \subset B_{1,\epsilon}$;
\item $\|\mathcal{T}_1\|<1$.
\end{enumerate}

By Proposition \ref{linearprop1} we get:
\begin{equation}\label{eq1n1}
\|\mathcal{T}_1(\phi_1)\| \leq c \|\mathcal{N}_1(\phi_1)\|
\end{equation}
and
\begin{equation}\label{eq2n1}
\|\mathcal{T}_1(\phi_1)-\mathcal{T}_1(\psi_1)\| \leq c (\|\mathcal{N}_1(\phi_1)-\mathcal{N}_1(\psi_1)\|),
\end{equation}
for all $\phi_1, \psi_1 \in {\mathcal K}_1^\perp$. Thanks to \eqref{stimaistar} and the definition of $\mathcal{N}_1$ we deduce that
\begin{equation}\label{eq3n1}
\|\mathcal{N}_1(\phi_1)\| \leq c | f(-\tau e_1+\phi_1) - (\lambda_1-\lambda) \phi_1|_{\frac{2N}{N+2}},
\end{equation}
and
\begin{equation}\label{eq4n1}
\displaystyle \|\mathcal{N}_1(\phi_1)-\mathcal{N}_1(\psi_1)\| \leq \displaystyle c | f(-\tau e_1+\phi_1) - f(-\tau e_1+\psi_1)-(\lambda_1-\lambda) (\phi_1-\psi_1)|_{\frac{2N}{N+2}} .
\end{equation}
Now we estimate the  right-hand term in (\ref{eq1n1}).

By definition and thanks to Lemma \ref{lem00n2} we have the following:
\begin{equation}\label{eq5n1}
 | f(-\tau e_1+\phi_1)-(\lambda_1-\lambda)\phi_1|\leq c (\tau^p |e_1|+|\phi_1|^p+\epsilon |\phi_1|).
\end{equation}
Since $|e_1| \leq \|e_1\|_\infty$, $p \frac{2N}{N+2} = \frac{2N}{N-2}$ and $|\phi_1^p|_{\frac{2N}{N+2}}=|\phi_1|_{\frac{2N}{N-2}}^{p}$, from (\ref{eq5n1}) and the Sobolev inequality we deduce the following:
\begin{equation}\label{eq6n1}
 |  f(-\tau e_1+\phi_1)-(\lambda_1-\lambda) \phi_1|_{\frac{2N}{N+2}} \leq c_1(\tau^p + \|\phi_1\|^p+\epsilon\|\phi_1\|).
\end{equation}

Observe that since $\tau=O(\epsilon^{3/4})$, $p=\frac{7}{3}$ then $\tau^p=O(\epsilon^{7/4})$, in particular $\tau^p=o(\epsilon^{\theta_1/2})$. Hence, thanks to (\ref{eq1n1}), Proposition \ref{errorprop1}, (\ref{eq3n1}), (\ref{eq6n1}) and since $p>1$, there exist $c>0$ and $\epsilon_0=\epsilon_0(\eta)>0$ such that  $$\|\phi_1\| \leq c \epsilon^{\frac{\theta_1}{2}+ \sigma} \Rightarrow \|\mathcal{T}_1(\phi_1) \|\leq c \epsilon^{\frac{\theta_1}{2}+ \sigma},$$
for all $\epsilon \in (0,\epsilon_0)$, for all $d_1 \in \mathbb R_+$ satisfying \eqref{limdjfj} (with $j=1$), for some positive real number  $\sigma$, whose choice depends only on $N$.
In other words $\mathcal{T}_1$ maps the ball $B_{1,\epsilon}$ into itself and (1) is proved.

We want to show that $\mathcal{T}_1$ is a contraction. To do this for any $\phi_1, \psi_1\in B_{1, \e}$ we write
\begin{eqnarray*}
&&\left|f(-\tau e_1+\phi_1)-f(-\tau e_1+\psi_1)-(\lambda_1-\lambda) (\phi_1-\psi_1)\right| \\[12pt]
=&&\left|f(-\tau e_1+\psi_1+(\phi_1-\psi_1))-f(-\tau e_1+\psi_1)-f^\prime(-\tau e_1+\psi_1)(\phi_1-\psi_1) - (\lambda_1-\lambda) (\phi_1-\psi_1)\right|
\end{eqnarray*}

 By using Lemma \ref{lem0n3}, we get that
 \begin{equation}\label{stimaeln1}
\begin{array}{lll}
&&\displaystyle\left|f(-\tau e_1+\psi_1+(\phi_1-\psi_1))-f(-\tau e_1+\psi_1)-f^\prime(-\tau e_1+\psi_1)(\phi_1-\psi_1)\right|\\[10pt]
 &\leq&\displaystyle c\left(\tau^{p-2}|e_1|^{p-2}|\phi_1-\psi_1|+|\phi_1|^{p-1}+|\psi_1|^{p-1}\right)|\phi_1-\psi_1|\\[10pt]
& \leq &\displaystyle c\left(\tau^{p-2}\|e_1\|_\infty^{p-2}|\phi_1-\psi_1|+|\phi_1|^{p-1}+|\psi_1|^{p-1}\right)|\phi_1-\psi_1|\\[10pt]
&  \leq &\displaystyle c_1\left(\tau^{p-2}|\phi_1-\psi_1|+|\phi_1|^{p-1}+|\psi_1|^{p-1}\right)|\phi_1-\psi_1|
\end{array}
\end{equation}
By direct computation $(p-1) \frac{2N}{N+2}=\frac{8N}{(N-2)(N+2)}$, so, since $\tau^{p-2}|\phi_1-\psi_1|^{ \frac{2N}{N+2}}$, $ |\phi_1|^{(p-1) \frac{2N}{N+2}}, |\psi_1|^{(p-1) \frac{2N}{N+2}} \in L^{\frac{N+2}{4}}$, $|\phi_1-\psi_1|^{\frac{2N}{N+2}} \in L^{p}$ and $1= \frac{4}{N+2}+\frac{N-2}{N+2}$, by H\"older inequality, we get that
\begin{eqnarray}\label{eq8n1}
\nonumber
&&\displaystyle \left| \left((\tau^{p-2}|\phi_1-\psi_1|+ |\phi_1|^{p-1}+|\psi_1|^{p-1}\right) (\phi_1 - \psi_1) \right|_{\frac{2N}{N+2}}\\
\nonumber
& \leq & c \displaystyle\left[ \left((\tau^{p-2}|\phi_1-\psi_1|_{\frac{2N}{N-2}}+| \phi_1|_{\frac{2N}{N-2}}^{\frac{4}{N-2}} + |\psi_1|_{\frac{2N}{N-2}}^{\frac{4}{N-2}}\right)^{\frac{2N}{N+2}} \ \left( | \phi_1 - \psi_1|_{\frac{2N}{N-2}}^{{\frac{2N}{N-2}}} \right)^{\frac{N-2}{N+2}} \right]^{\frac{N+2}{2N}}\\
&=& c \displaystyle \left(\tau^{p-2}|\phi_1-\psi_1|_{\frac{2N}{N-2}}+| \phi_1|_{\frac{2N}{N-2}}^{\frac{4}{N-2}} + |\psi_1 |_{\frac{2N}{N-2}}^{\frac{4}{N-2}}\right) \ | \phi_1 - \psi_1|_{\frac{2N}{N-2}}.
\end{eqnarray}
We point out that since

Hence by (\ref{eq2n1}), (\ref{eq4n1}), \eqref{stimaeln1}, (\ref{eq8n1}) and Sobolev inequality we get that there exists $L\in (0,1)$ such that
$$\|\phi_1\| \leq c \epsilon^{\frac{\theta_1}{2}+ \sigma}, \|\psi_1\| \leq c \epsilon^{\frac{\theta_1}{2}+ \sigma} \Rightarrow \|\mathcal{T}_1(\phi_1)- \mathcal{T}_1(\psi_1)\| \leq L \|\phi_1-\psi_1\|.$$

Hence by the Contraction Mapping Theorem we can
uniquely solve $\mathcal{T}_1(\phi_1)=\phi_1$ in $B_{1,\epsilon}$. We denote by $\bar\phi_1 \in B_{1,\epsilon}$ this solution.  A standard argument shows that $d_1 \rightarrow \bar\phi_1(d_1)$ is a $C^1$-map (see also \cite{Musso1}). The proof is then concluded.
\end{proof}
\begin{remark}\label{stimafixpoint1}
Let us fix a small $\eta>0$. We observe that if $\bar\phi_1$ is the fixed point of  $\mathcal{T}_1$ in $B_{1,\epsilon}$ found in the previous proposition, then, by using the estimate \eqref{eq6n1}, taking into account that $\|\bar\phi_1\| \leq c \e^{5/8}$ and by making elementary computations, we get that $\|\bar \phi_1\|\leq c \tau^{1+\gamma}$, for all $0<\gamma\leq \frac{17}{24}$, for all $\e$ sufficiently small, for all $d_1 \in (\eta, \frac{1}{\eta})$.
\end{remark}
\subsection{The solution of the second equation of \eqref{sistemaaux}}
The aim of this section is to prove the following result:

\begin{proposition}\label{aux2}
Let $N=5$, $\tau$ and $\delta$  as in \eqref{deltatau5}. 
Let $\bar\phi_1\in\mathcal K_1^\bot$ be the solution of the first equation in \eqref{sistemaaux} found in Proposition \ref{aux1}.\\ Then, for any $\eta>0$, there exists $\epsilon_0>0$ and $c>0$ such that for all $\epsilon \in (0,\epsilon_0)$, for all $(d_1,d_2) \in \mathbb R^2_+$ satisfying condition \eqref{limdjfj}, there exists a unique solution $\bar\phi_2\in\mathcal{K}^\bot$ of the second equation in \eqref{sistemaaux}, such that
\begin{equation}\label{stimaphi2}
\|\bar\phi_2\|\leq c \e^{\frac{\theta_2}{2}+\sigma},
\end{equation}
where $\theta_2$ is defined in \eqref{thetaj} and $\sigma$ is some positive real number depending only on $N$.
Moreover $\bar\phi_2$ is continuously differentiable with respect to $(d_1, d_2)$.
\end{proposition}

In order to prove Proposition \ref{aux2} we have to estimate the error term $\mathcal{R}_2$ defined in \eqref{defelementi3}. It holds the following result.

\begin{proposition}\label{errorprop2}
For any $\eta>0$, there exists $\epsilon_0>0$ and $c>0$ such that  for all $\epsilon \in (0,\epsilon_0)$, for all $(d_1,d_2) \in \mathbb R_+^2$ satisfying \eqref{limdjfj}, we have $$\| \mathcal{R}_2\|\leq c \ \epsilon^{\frac{\theta_2}{2}+\sigma},$$ for some positive real number  $\sigma$, whose choice depends only on $N$.
\end{proposition}
\begin{proof}
By continuity of $\Pi^\bot$ and by definition of $i^*$, using \eqref{stimaistar}, we deduce that
\begin{eqnarray}\label{eq1err2}
\nonumber
\displaystyle \|\mathcal{R}_2\| &\leq& \displaystyle c |\lambda\p\U_\delta+f(\p\U_\delta-\tau e_1)-f(U_\delta)|_{\frac{2N}{N+2}} \\[10pt]
\nonumber
&\leq& \underbrace{\displaystyle c |f(\p\U_\delta-\tau e_1)-f(\p\U_\delta)|_{\frac{2N}{N+2}}}_{(I)}+\underbrace{\displaystyle c |f(\p\U_{\delta})-f(\U_{\delta})|_{\frac{2N}{N+2}}}_{(II)}  \\[10pt]
&& + \underbrace{c\lambda |\p\U_{\delta}|_{\frac{2N}{N+2}}}_{(III)}.
\end{eqnarray}
Let us fix $\eta>0$. Reasoning as in the proof of Proposition \ref{errore4} we get that

$$(I)=|f(\p\U_\delta-\tau e_1)-f(\p\U_\delta)|_{\frac{2N}{N+2}} \leq c \left(\tau \delta^{\frac{N-2}{2}}+ \tau^p+\tau^2\d^{\frac{6-N}{2}}\right).$$
Recalling the choice of $\tau$ and $\delta$ (see \eqref{deltatau5}), we deduce that

 \begin{equation}\label{eq5err2}
(I) = |f(\p\U_\delta-\tau e_1)-f(\p\U_\delta)|_{\frac{2N}{N+2}} \leq c  \left(\epsilon^{\frac{3}{4}}\epsilon^{\frac{3}{2}\frac{3}{2}}+ \epsilon^{\frac{3}{4} \frac{7}{3}}+\e^{\frac 94}\right)\leq c \e^{\frac{3}{2}+\sigma}
\end{equation}
where $c$ depends on $\eta$ (and also on $\Omega$, $N$), $\sigma$ is some positive real number depending only on $N$.

In order to estimate $(II)$, writing $\p\U_\delta=\U_\delta-\varphi_\delta$, using \eqref{stimaproiezione}, \eqref{stimavarphi}, and arguing as in \eqref{stimautilefuturo1N4}, we get
\begin{equation}\label{stimautilefuturo1}
\begin{array}{lll}
\displaystyle \int_{\Omega}|(\p\U_{\delta})^p -\U_{\delta}^p|^{\frac{2N}{N+2}}\ dx &\leq &\displaystyle c_2\left( \delta^{\frac{2N(N-2)}{N+2}} \int_{\Omega/\delta}\frac{1}{(1+|y|^2)^{\frac{4N}{N+2}}}\  dy + \delta^{N} + \d^{\frac{2N(N-2)}{N+2}+\frac{N(6-N)}{N+2}}\right).
\end{array}
\end{equation}
We observe that  for $N=5$ we have that $\frac{1}{(1+|y|^2)^{\frac{4N}{N+2}}} \in L^1(\R)$, and hence we get that
$$ \delta^{\frac{2N(N-2)}{N+2}} \int_{\Omega/\delta}\frac{1}{(1+|y|^2)^{\frac{4N}{N+2}}}\  dy =  \delta^{\frac{2N(N-2)}{N+2}} \int_{\R}\frac{1}{(1+|y|^2)^{\frac{4N}{N+2}}}\  dy + o( \delta^{\frac{2N(N-2)}{N+2}}).$$
Hence it follows that
\begin{equation}\label{stimautilefuturo2}
\begin{array}{lll}
(II)=\displaystyle \left( \int_{\Omega}|(\p\U_{\delta})^p -\U_{\delta}^p|^{\frac{2N}{N+2}}\ dx \right)^{\frac{N+2}{2N}}
&\leq & \displaystyle c_3 \left( \delta^{N-2} + \delta^{\frac{N+2}{2}}+\d^{N-2+\frac{6-N}{2}}\right)\\[10pt]
&\leq & \displaystyle c_4  \delta^{N-2} \leq c_5 \epsilon^{\frac{\theta_2}{2}+\sigma}.
\end{array}
\end{equation}
for all sufficiently small $\epsilon$, for some positive $\sigma$ depending only on $N$.\\
It remains to estimate $(III)$. Since $\p\U_\delta\leq \U_\delta$ we have:
\begin{equation*}
\begin{array}{lll}
\displaystyle \int_{\Omega} \p\U_{\delta_{2}}^{\frac{2N}{N+2}} \ dx
&\leq& \displaystyle  \alpha_N^{\frac{2N}{N+2}} \delta^{\frac{N(N-2)}{N+2}} \int_{\Omega}\frac{1}{|x|^{\frac{N(N-2)}{N+2}}}\  dx.
\end{array}
\end{equation*}
Since $N=5$ we have that $\frac{1}{|x|^{\frac{N(N-2)}{N+2}}}$ is integrable near the origin and hence we deduce that $(III) \leq c \ \epsilon \delta^{\frac{N-2}{2}}$.
Thanks to the choice $\delta$, by an elementary computation, we get that:
$$(III)\leq c \ \epsilon^{\frac{9}{4}}\leq c\e^{\frac{\theta_2}{2}+\sigma},$$
for all sufficiently small $\epsilon$, for some positive $\sigma$ depending only on $N$.

Finally, putting together all these estimates we deduce that there exists a positive constant $c=c(\eta)>0$ and $\epsilon_0=\epsilon_0(\eta)>0$ such that for all $\epsilon \in (0,\epsilon_0)$, for all $(d_1,d_2) \in \mathbb R_+^2$ satisfying \eqref{limdjfj}
$$ \|\mathcal{R}_2\| \leq c \epsilon^{\frac{\theta_2}{2} + \sigma},$$
for some positive real number $\sigma$ (whose choice depends only on $N$).
The proof is complete.
\end{proof}
We are now in position to prove Proposition \ref{aux2}.\\

\begin{proof}[proof of Proposition \ref{aux2}]
Let us fix $\eta>0$ and define the operator $\mathcal{T}_2 :\mathcal{K}^{\perp} \rightarrow \mathcal{K}^{\perp} $ as $$\mathcal{T}_2(\phi_2):=- \mathcal{L}_2^{-1}[\mathcal{N}_1(\bar\phi_1,\phi_2)+\mathcal{R}_2],$$
where $\bar\phi_1 \in \mathcal{K}_1^\perp \cap B_{1,\epsilon}$ be the unique solution of the first equation of \eqref{sistemaaux} found in Proposition \ref{aux1}.\\
In order to find a solution of the second equation of \eqref{sistemaaux} we solve the fixed point problem $\mathcal{T}_2(\phi_2)=\phi_2$.  Let us define the proper ball

$$\displaystyle B_{2,\epsilon}:=\{\phi_2 \in \mathcal{K}^{\perp}: \ \|\phi_2\| \leq r\ \epsilon^{\frac{\theta_2}{2}+\sigma} \}$$
for $r>0$ sufficiently large and $\sigma>0$ to be chosen later.\\
From Proposition \ref{p:inv4}, there exist $\epsilon_0=\epsilon_0(\eta)>0$ and $c=c(\eta)>0$ such that:
\begin{equation}\label{eq1n2}
\|\mathcal{T}_2(\phi_2)\| \leq c (\|\mathcal{N}_2(\bar\phi_1,\phi_2)\| + \|\mathcal{R}_2\|),
\end{equation}
and
\begin{equation}\label{eq2n2}
\|\mathcal{T}_2(\phi_2)-\mathcal{T}_2(\psi_2)\| \leq c (\|\mathcal{N}_2(\bar\phi_1, \phi_2)-\mathcal{N}_2(\bar\phi_1, \psi_2)\|),
\end{equation}
for all $\phi_2, \psi_2 \in \mathcal{K}^\perp$, for all $(d_1,d_2) \in \mathbb R_+^2$  satisfying \eqref{limdjfj} and for all $\epsilon \in (0,\epsilon_0)$.

We begin with estimating the right hand side of (\ref{eq1n2}).\\ Thanks to Proposition \ref{errorprop2}  we have that
$$\| \mathcal{R}_2\| \leq c \epsilon^{\frac{\theta_2}{2}+\sigma},$$
for all $\epsilon \in (0,\epsilon_0)$, for all $(d_1,d_2) \in \mathbb R_+^2$ satisfying \eqref{limdjfj}. Thus it remains only to estimate $\|\mathcal{N}_2(\bar\phi_1,\phi_2)\|$. Thanks to \eqref{stimaistar} and the definition of $\mathcal{N}_2$ we have:
\begin{equation}\label{eq3n2}
\|\mathcal{N}_2(\bar\phi_1,\phi_2)\| \leq c | f(\p\U_\delta - \tau e_1+\bar\phi_1+\phi_2)- f^\prime(\U_\delta) \phi_2 -f(-\tau e_1+\bar\phi_1)-f(\p\U_{\delta}-\tau e_1)) |_{\frac{2N}{N+2}}.
\end{equation}
We estimate the right-hand side of (\ref{eq3n2}): 
\begin{equation*}
\begin{array}{lll}
&&\displaystyle | f(\p\U_\delta - \tau e_1+\bar\phi_1+\phi_2)- f^\prime(\U_\delta) \phi_2 -f(-\tau e_1+\bar\phi_1)-f(\p\U_{\delta}-\tau e_1))|_{\frac{2N}{N+2}} \\[10pt]
&\leq&\displaystyle | f(\p\U_\delta - \tau e_1+\bar\phi_1+\phi_2) - f(\p\U_\delta - \tau e_1+\bar\phi_1) - f^\prime(\p\U_\delta - \tau e_1+\bar\phi_1) \phi_2|_{\frac{2N}{N+2}} \\[10pt]
&+&  \displaystyle | f(\p\U_\delta - \tau e_1+\bar\phi_1) - f(\p\U_\delta - \tau e_1) - f^\prime(\p\U_\delta - \tau e_1)  \bar\phi_1 |_{\frac{2N}{N+2}}\\[10pt]
&+&  \displaystyle | [f^\prime(\p\U_\delta - \tau e_1+\bar\phi_1) - f^\prime(\p\U_\delta - \tau e_1)]\phi_2|_{\frac{2N}{N+2}}\\[10pt]
&+&\displaystyle|  [f^\prime(\p\U_\delta - \tau e_1) - f^\prime(\p\U_\delta) ]\phi_2  |_{\frac{2N}{N+2}}+\displaystyle|  [f^\prime(\p\U_\delta ) - f^\prime(\U_\delta) ]\phi_2  |_{\frac{2N}{N+2}}\\[10pt]
&+&\displaystyle|  [-f(- \tau e_1+\bar\phi_1) + f(- \tau e_1) + f^\prime(- \tau e_1)\bar\phi_1|_{\frac{2N}{N+2}}\\[10pt]
&+&|f^\prime(- \tau e_1)\bar\phi_1|_{\frac{2N}{N+2}} + |[f^\prime(\p\U_\delta - \tau e_1)-f^\prime(-\tau e_1)]\bar\phi_1|_{\frac{2N}{N+2}}.
\end{array}
\end{equation*}

In order to estimate these terms, by  Lemma \ref{lem0n2} and Lemma \ref{lem1n2} we deduce that:
\begin{equation*} \label{eq4bisn2}
|f(\p\U_\delta - \tau e_1+\bar\phi_1+\phi_2) - f(\p\U_\delta - \tau e_1+\bar\phi_1) - f^\prime(\p\U_\delta - \tau e_1+\bar\phi_1) \phi_2| \leq  c ( |\p\U_\delta - \tau e_1+\bar\phi_1 |^{p-2}\phi_2^2 + |\phi_2|^p),
\end{equation*}

\begin{equation*} \label{eq4bisn2a}
|f(\p\U_\delta - \tau e_1+\bar\phi_1) - f(\p\U_\delta - \tau e_1) - f^\prime(\p\U_\delta - \tau e_1)  \bar\phi_1| \leq  c ( |\p\U_\delta - \tau e_1 |^{p-2}\bar\phi_1^2 + |\bar\phi_1|^p),
\end{equation*}

\begin{equation*} \label{eq4bisn2b}
| [f^\prime(\p\U_\delta - \tau e_1+\bar\phi_1) - f^\prime(\p\U_\delta - \tau e_1)]\phi_2| \leq  c ( |\p\U_\delta - \tau e_1 |^{p-2}|\bar\phi_1||\phi_2| + |\bar\phi_1|^{p-1}|\phi_2|),
\end{equation*}

\begin{equation*} \label{eq4bisn2c}
| f^\prime(\p\U_\delta - \tau e_1) - f^\prime(\p\U_\delta) ]\phi_2| \leq  c ( |\p\U_\delta|^{p-2}| \tau e_1||\phi_2| + | \tau e_1|^{p-1}|\phi_2|),
\end{equation*}

\begin{equation*} \label{eq4bisn2d}
| f^\prime(\p\U_\delta) - f^\prime(\U_\delta) ]\phi_2| \leq  c ( |\U_\delta|^{p-2}|\varphi_\delta||\phi_2| + |\varphi_\delta|^{p-1}|\phi_2|),
\end{equation*}

\begin{equation*}\label{eq4bisn21}
|-f(- \tau e_1+\bar\phi_1) + f(- \tau e_1) + f^\prime(- \tau e_1)\bar\phi_1|\leq  c\left( | \tau e_1|^{p-2}|\bar\phi_1|^2+ |\bar\phi_1|^p\right),
\end{equation*}

\begin{equation} \label{eq4bisn2e}
| [f^\prime(\p\U_\delta - \tau e_1)-f^\prime(-\tau e_1)]\bar\phi_1| \leq  c ( |\p\U_\delta|^{p-2}| \tau e_1||\bar\phi_1| + | \tau e_1|^{p-1}|\bar\phi_1|).
\end{equation}\\

Now we have to estimate the $L^{\frac{2N}{N+2}}$-norm of every term. We give a detailed proof of the first two terms; for the remaining ones the proof is similar. \\

By the triangular inequality and Lemma \ref{lem00n2}, it holds
\begin{equation}\label{stiman2pt1}
  |\p\U_\delta - \tau e_1+\bar\phi_1 |^{p-2}\phi_2^2 + |\phi_2|^p \leq |\p\U_\delta|^{p-2}\phi_2^2+|\tau e_1|^{p-2}\phi_2^2+|\bar\phi_1|^{p-2}\phi_2^2 + |\phi_2|^p.
\end{equation}
Since $p-2=\frac{6-N}{N-2}$, we have $ (\p\U_\delta^{p-2})^{\frac{2N}{N+2}} =\p\U_\delta^{\frac{2N(6-N)}{(N-2)(N+2)}} \leq \U_\delta^{\frac{2N(6-N)}{(N-2)(N+2)}}\leq c  \delta^{-\frac{N(6-N)}{N+2}} $.
Hence, by the previous estimate and thanks to the Sobolev inequality, we get that
\begin{eqnarray*}
 \left(\int_\Omega \left(\p\U_\delta^{p-2}\phi_2^2\right)^{\frac{2N}{N+2}} \ dx\right)^{\frac{N+2}{2N}} &\leq& c \left(\delta^{-N\frac{6-N}{N+2}} \int_\Omega \phi_2^{\frac{4N}{N+2}} \ dx\right)^{\frac{N+2}{2N}}\\[10pt]
 &\leq& c_1 \delta^{-\frac{6-N}{2}}\left( \int_\Omega \phi_2^{\frac{4N}{N+2}} \ dx\right)^{\frac{N+2}{2N}}\\[10pt]
  &\leq& c_2 \delta^{-\frac{6-N}{2}}\|\phi_2\|^{2}.
\end{eqnarray*}
We observe that for $N=5$, and thanks to the choice of $\delta$ we have  $$\delta^{-\frac{6-N}{2}}\|\phi_2\|^{2}\leq c \epsilon^{-\frac{1}{2}\frac{3}{2}}\epsilon^{3} \leq c \epsilon^{\frac{\theta_2}{2}+\sigma},$$
for all sufficiently small $\epsilon$.

The remaining terms of \eqref{stiman2pt1} are even simpler. In fact, it holds
\begin{eqnarray*}
 \left(\int_\Omega \left(|\tau e_1|^{p-2}\phi_2^2\right)^{\frac{2N}{N+2}} \ dx\right)^{\frac{N+2}{2N}} &\leq& c \left(\tau^{\frac{2N(6-N)}{(N+2)(N-2)}}\|e_1\|_\infty^{\frac{2N}{N+2}} \int_\Omega \phi_2^{\frac{4N}{N+2}} \ dx\right)^{\frac{N+2}{2N}}\\[10pt]
 &\leq& c_1\tau^{\frac{6-N}{N-2}} \left( \int_\Omega \phi_2^{\frac{4N}{N+2}} \ dx\right)^{\frac{N+2}{2N}}\\[10pt]
  &\leq& c_2\tau^{\frac{6-N}{N-2}}\|\phi_2\|^{2}.
\end{eqnarray*}
Thanks to choice of $\tau$ and $\delta$ it is clear that $\tau^{\frac{6-N}{N-2}}\|\phi_2\|^{2}\leq c \epsilon^{\frac{\theta_2}{2}+\sigma}$.

For the third term, applying H\"older inequality and Sobolev inequality we have
\begin{eqnarray*}
 \left(\int_\Omega \left(|\bar\phi_1|^{p-2}\phi_2^2\right)^{\frac{2N}{N+2}} \ dx\right)^{\frac{N+2}{2N}}
  &\leq& c \|\bar\phi_1\|^{p-2}\|\phi_2\|^{2},
\end{eqnarray*}
and it is clear that $\|\bar\phi_1\|^{p-2}\|\phi_2\|^{2} \leq c  \epsilon^{\frac{\theta_2}{2}+\sigma}$.
Finally, since $\frac{2N}{N+2} \cdot p = p+1$ by the Sobolev inequality we get that
\begin{eqnarray*}
 \left(\int_\Omega |\phi_2|^{p\frac{2N}{N+2}} \ dx\right)^{\frac{N+2}{2N}}
  &\leq& c \|\phi_2\|^{p},
\end{eqnarray*}
and the estimate of the first term is done. For the second term of \eqref{eq4bisn2e}, reasoning as in \eqref{stiman2pt1}, we get that
\begin{equation}\label{stiman2pt2}
|f(\p\U_\delta - \tau e_1+\bar\phi_1) - f(\p\U_\delta - \tau e_1) - f^\prime(\p\U_\delta - \tau e_1)  \bar\phi_1| \leq c \left( |\p\U_\delta|^{p-2}\bar\phi_1^2+|\tau e_1|^{p-2}\bar\phi_1^2+ |\bar\phi_1|^p\right).
  \end{equation}
As be before we have
\begin{eqnarray*}
 \left(\int_\Omega \left(\p\U_\delta^{p-2}\bar\phi_1^2\right)^{\frac{2N}{N+2}} \ dx\right)^{\frac{N+2}{2N}} &\leq& c  \delta^{-\frac{6-N}{2}}\|\bar\phi_1\|^{2}.
\end{eqnarray*}
Thanks to the choice of $\delta$ and since $\bar\phi \in B_{1,\epsilon}$, then, for all sufficiently small $\epsilon$,  we have
$$\delta^{-\frac{6-N}{2}}\|\bar\phi_1\|^{2} \leq c \epsilon^{-\frac{1}{2}\frac{3}{2}}\epsilon^{\frac{5}{2}}\leq c \epsilon^{\frac{\theta_2}{2}+ \sigma} .$$
This estimate is trivially true for $||\tau e_1|^{p-2}\bar\phi_1^2|_{\frac{2N}{N+2}}$, let us check the estimate of the $L^{\frac{2N}{N+2}}$-norm of the last term of \eqref{stiman2pt2}. As before we have
\begin{eqnarray*}
 \left(\int_\Omega |\bar\phi_1|^{p\frac{2N}{N+2}} \ dx\right)^{\frac{N+2}{2N}}
  &\leq& c \|\bar\phi_1\|^{p},
\end{eqnarray*}
and for all small $\epsilon$ we have $\|\bar\phi_1\|^{p} \leq c \epsilon^{\frac{5}{4}\frac{7}{3}}\leq c \epsilon^{\frac{\theta_2}{2}+\sigma}$.\\
The same techniques used so far apply to the other remaining estimates, we limit ourselves to checking the lowest order term.\\
Thanks to H\"older and Sobolev inequalities we have $$\big|  \p\U_\delta^{p-2} \tau e_1 \bar\phi_1\big|_{\frac{2N}{N+2}}\leq c |\p\U_\delta|_{\frac{2N}{N-2}}  \tau |e_1|_\infty  \|\bar\phi_1\|$$
and thanks to the choice of $\tau$ and since $\bar\phi_1\in B_{1,\epsilon}$ we have
$$ |\p\U_\delta|_{\frac{2N}{N-2}}  \tau \|e_1\|_\infty  \|\bar\phi_1\| \leq c  \epsilon^{\frac{3}{4}}\epsilon^{\frac{5}{4}} \leq c \epsilon^{\frac{\theta_2}{2}+\sigma}.$$

It remains to prove that $\mathcal{T}_2:B_{2,\epsilon} \rightarrow B_{2,\epsilon} $ is a contraction. Thanks to (\ref{eq2n2}) it suffices to estimate $\|\mathcal{N}_2(\bar\phi_1, \phi_2)-\mathcal{N}_2(\bar\phi_1, \psi_2)\|$ for any $\psi_2, \phi_2 \in B_{2,\epsilon}$. To this end, thanks to \eqref{stimaistar}, the definition of $\mathcal{N}_2$ we have:
\begin{equation*}
\|\mathcal{N}_2(\bar\phi_1,\phi_2)-\mathcal{N}_2(\bar\phi_1,\psi_2)\| \leq c | f(\p\U_\delta - \tau e_1+\bar\phi_1+\phi_2)-f(\p\U_\delta - \tau e_1+\bar\phi_1+\psi_2)- f^\prime(\U_\delta) (\phi_2-\psi_2)|_{\frac{2N}{N+2}},
\end{equation*}
then, reasoning as in the proof of Proposition \ref{aux1}, we get that
\begin{equation}\label{n2contract}
\|\mathcal{N}_2(\bar\phi_1,\phi_2)-\mathcal{N}_2(\bar\phi_1,\psi_2)\| \leq \e^\alpha \|\phi_2-\psi_2\|,
\end{equation}
for some $\alpha>0$, for all sufficiently small $\epsilon>0$.\\\\ At the end, thanks to \eqref{eq2n2} and \eqref{n2contract}, we get that there exists $L\in (0,1)$ such that
$$ \|\mathcal{T}_2(\phi_2)- \mathcal{T}_2(\psi_2)\| \leq L \|\phi_2-\psi_2\|.$$
The proof is complete.
\end{proof}

Putting together Proposition \ref{aux1} and \ref{aux2} we get the proof of Proposition \ref{auxsolving}.

\begin{remark}\label{stimafixpoint2}
Let us fix a small $\eta>0$. We observe that the fixed point $\bar\phi_2$  of  $\mathcal{T}_2$ in $B_{2,\epsilon}$, found in Proposition \ref{aux2}, verifies $\|\bar\phi_2\| \leq c \e^{3/2}$, hence, by definition we get that $\|\bar \phi_2\|\leq c \tau^{1+\gamma}$, for all $0<\gamma\leq \frac{3}{4}$, for all $\e$ sufficiently small, for all $(d_1,d_2) \in (\eta, \frac{1}{\eta})\times (\eta, \frac{1}{\eta})$.
\end{remark}

\subsection{Estimates for the reduced functional}
We are left now to solve \eqref{bif}. Let $(\bar\phi_1, \bar\phi_2)\in \mathcal{K}^\bot_1\times\mathcal{K}^\bot$ be the solution found in Proposition \ref{auxsolving}. Hence $V_\lambda+\bar\phi_1+\bar\phi_2$ is a solution of our original problem \eqref{BN} if we can find $\bar d_\lambda=(\bar d_{1\lambda}, \bar d_{2\lambda})$  which satisfies condition \eqref{limdjfj} and solves equation \eqref{bif}.

 To this end we consider the reduced functional $\tilde J_\lambda:\mathbb R_+^2 \rightarrow \mathbb R$ defined by: $$\tilde J_\lambda(d_1,d_2):=J_\lambda(V_\lambda+\bar\phi_1+\bar\phi_2),$$ where $J_\lambda$ is the functional defined in \eqref{funzionale}.\\
Our main goal is to show first that solving equation (\ref{bif}) is equivalent to finding critical points $(\bar d_{1,\lambda},\bar d_{2,\lambda})$ of the reduced functional $\tilde J_\lambda(d_1,d_2)$ and then that the reduced functional has a critical point.

\begin{lemma}\label{CP}
For any small $\eta>0$ there exists $\e_0>0$ such that for all $\lambda \in (\lambda_1-\e_0, \lambda_1)$ if $(\bar{d}_{1,\lambda}, \bar{d}_{2,\lambda})$ is a critical point of $\tilde J_\lambda$ and satisfies \eqref{limdjfj}, then $V_\lambda+\bar\Phi_\lambda$ is a solution of \eqref{BN}.
\end{lemma}
\begin{proof}
Let us fix a small $\eta>0$ and let $\bar\Phi_\lambda=\bar\phi_1+\bar\phi_2$ be the solution of the auxiliary equation \eqref{aux}, found Proposition \ref{auxsolving}. By assumption $(\bar{d}_{1, \lambda}, \bar d_{2,\lambda})$ is a critical point of $\tilde J_\lambda(d_1, d_2)$.  Hence $$\nabla \tilde J_\lambda(\bar d_{1, \lambda}, \bar d_{2, \lambda})=\left(\partial_{d_1} \tilde J_\lambda, \partial_{d_2}\tilde J_\lambda\right)_{|_{(\bar d_{1, \lambda}, \bar d_{2, \lambda})}}=(0, 0).$$
Then for $j=1, 2$
\begin{eqnarray*}
0&=& \partial_{d_j}\tilde J_\lambda = J'_\lambda(V_\lambda +\bar\Phi_\lambda)[\partial_{d_j}V_\lambda+\partial_{d_j}\bar\Phi_\lambda]\\
&=& \left(V_\lambda+\bar\Phi_\lambda-i^*\left[f(V_\lambda+\bar\Phi_\lambda)+\lambda(V_\lambda+\bar\Phi_\lambda)\right], \partial_{d_j}V_\lambda+\partial_{d_j}\bar\Phi_\lambda\right)_{\zo}\\
&\underbrace{=}_{\bar\Phi_\lambda \,\, \mbox{solution of \eqref{aux}}}& \left(c_0^\lambda \p Z_\d+c_1^\lambda e_1, \partial_{d_j}V_\lambda+\partial_{d_j}\bar\Phi_\lambda\right)_{\zo}\\
&=&\left\{\begin{array}{lr}\left(c_0^\lambda \p Z_\d+c_1^\lambda e_1, \partial_{d_1}(-\tau e_1)+\partial_{d_1}\bar\Phi_\lambda\right)_{\zo}\qquad \mbox{for}\,\, j=1\\\\
\left(c_0^\lambda \p Z_\d+c_1^\lambda e_1, \partial_{d_2}(\p \U_\d)+\partial_{d_2}\bar\Phi_\lambda\right)_{\zo}\qquad \mbox{for}\,\, j=2 \end{array}\right.
\end{eqnarray*}
for some $c_0^\lambda, c_1^\lambda\in\mathbb R$. \\ We shall prove that $c_0^\lambda, c_1^\lambda$ are equal to zero for $\lambda$ near $\lambda_1$. Let $\e=\lambda_1-\lambda$.\\ Now we deduce that (recalling that $\tau=\e^{\frac 34}d_1$) $$\partial_{d_1}(-\tau e_1)=-\e^{\frac 34}e_1,$$ while (recalling that $\d= \e^{\frac 32}d_2$) we get that $$\partial_{d_2}\p\U_\d= \e^{\frac 32}\p Z_\d.$$
Hence
\begin{eqnarray}\label{primaeqc}
\nonumber
0&=& \left(c_0^\lambda \p Z_\d +c_1^\lambda e_1, -\e^{\frac 34}e_1+\partial_{d_1}\bar\Phi_\lambda\right)_{\zo}\\
\nonumber
&=&-\e^{\frac 34}c_0^\lambda (\p Z_\d, e_1)_{\zo}+ c_0^\lambda (\p Z_\d, \partial_{d_1}\bar\Phi_\lambda)_{\zo}-\e^{\frac 34}c_1^\lambda(e_1, e_1)_{\zo}\\
&&+c_1^\lambda (e_1, \partial_{d_1}\bar\Phi_\lambda)_{\zo}
\end{eqnarray}
and
\begin{eqnarray}\label{secondaeqc}
\nonumber
0&=& \left(c_0^\lambda \p Z_\d +c_1^\lambda e_1, \e^{\frac 32}\p Z_\d+\partial_{d_2}\bar\Phi_\lambda\right)_{\zo}\\
\nonumber
&=&\e^{\frac 32}c_0^\lambda (\p Z_\d, \p Z_\d)_{\zo}+ c_0^\lambda (\p Z_\d, \partial_{d_2}\bar\Phi_\lambda)_{\zo}+\e^{\frac 32}c_1^\lambda(e_1, \p Z_\d)_{\zo}\\
&&+c_1^\lambda (e_1, \partial_{d_2}\bar\Phi_\lambda)_{\zo}
\end{eqnarray}
Now:
\begin{eqnarray*}
(\p Z_\d, e_1)_{\zo}&=& \int_\O f'(\U_\d) Z_\d e_1\, dx = p \int_\O \U_\d^{p-1}Z_\d e_1\, dx\\
&=& p\alpha_N^p \d^{\frac{N-4}{2}}\int_{\frac{\O}{\d}}\frac{|y|^2-1}{(1+|y|^2)^{\frac{N+4}{2}}}e_1(\d y)\, dy \\
&=& A_0 \d^{\frac{N-4}{2}}+ o(\delta^{\frac{N-4}{2}})
\end{eqnarray*}
where $$A_0:=p\alpha_N^p e_1(0)\int_{\R}\frac{|y|^2-1}{(1+|y|^2)^{\frac{N+4}{2}}}\, dy.$$

By using \eqref{stimaproiezionederivata} we get
\begin{eqnarray*}
(\p Z_\d, \p Z_\d)_{\zo}&=& \int_\O f'(\U_\d) Z_\d\p Z_\d, dx\\
&=& \int_\O f'(\U_\d)Z_\d^2\, dx -\int_\O f'(\U_\d)Z_\d\psi_\d\, dx \\
&=&p\alpha_N^p \d^{N-2}\int_\O\frac{(|x|^2-\d^2)^2}{(\d^2+|x|^2)^{N+2}}\, dx + O\left(\d^{\frac N 2}|\psi_\d|_{\infty, \O}\int_\O\frac{|x|^2-\d^2}{(\d^2+|x|^2)^{\frac{N+4}{2}}}\, dx\right)\\
&=& B_0 \d^{-2}+ O(\d^{N-4})
\end{eqnarray*}
where $$B_0:=p\alpha_N^p \int_{\R}\frac{(|y|^2-1)^2}{(1+|y|^2)^{N+2}}\, dy.$$ Secondly, we point out that, since $\bar\Phi_\lambda\in \mathcal K^\bot$, it holds $$(e_1, \bar\Phi_\lambda)_{\zo}=0;\qquad (\p Z_\d, \bar\Phi_\lambda)_{\zo}=0.$$ These imply $$(e_1, \partial_{d_j}\bar\Phi_\lambda)_{\zo}=-(\partial_{d_j}e_1, \bar\Phi_\lambda)_{\zo}=0$$ and
\begin{equation*}
(\p Z_\d, \partial_{d_j}\bar\Phi_\lambda)=-(\partial_{d_j}\p Z_\d, \bar \Phi_\lambda)=\left\{\begin{array}{lr} 0\qquad\qquad\qquad\qquad\qquad \  \mbox{for}\,\, j=1\\\\ -(\partial_{d_2}\p Z_\d, \bar\Phi_\lambda)_{\zo}\qquad \mbox{for}\,\, j=2\end{array}\right.
\end{equation*}
and $$(\partial_{d_2}\p Z_\d, \bar\Phi_\lambda)_{\zo}=O(\|\partial_{d_2}\p Z_\d\|\cdot \|\bar\Phi_\lambda\|)=o(\|\partial_{d_2}\p Z_\d\|).$$  Now, since $$\|\partial_{d_2} \p Z_\d\|=O(\|\partial_{d_2}Z_\d\|),$$ by making some computations and using the fact that $N=5$, we get that

$$\partial_{d_2}Z_\d=\d^{\frac 52}\frac{3\d^2-7|x|^2}{(\d^2+|x|^2)^{\frac 72}}+\frac 12 \d^{\frac 12}d_2^{-1}\frac{|x|^2-\d^2}{(\d^2+|x|^2)^{\frac 52}}$$
Hence
\begin{eqnarray*}
\|\partial_{d_2} Z_\d\|^2 &=& \int_\O \frac{\d |x|^2}{(\d^2+|x|^2)^7}\left[35\d^2\frac{|x|^2-\d^2}{\d^2+|x|^2}+\frac 12 d_2^{-1}(7\d^2-3|x|^2)\right]^2\, dx\\
&=& C_0 \d^{-2}+O(\d)
\end{eqnarray*}
where $$C_0:=\int_{\R}\frac{|y|^2}{(1+|y|^2)^7}\left[35\frac{|y|^2-1}{1+|y|^2}+\frac 12 d_2^{-1}(7-3|y|^2)\right]^2\, dy$$ and hence $$(\partial_{d_2}\p Z_\d, \bar\Phi_\lambda)_{\zo}=o(\d^{-1}).$$
The equations \eqref{primaeqc}, \eqref{secondaeqc} become
\begin{equation}\label{sistemac}
\left\{
\begin{array}{lr}
c_0^\lambda (A_0\e^{\frac 32}\bar d_{2,\lambda}+O(\e^2))+\e^{\frac 34} c_1^\lambda D=0\\\\
c_0^\lambda(B_0 \bar d_{2,\lambda}+ o(1))+c_1^\lambda (A_0 \bar{d}_{2,\lambda}\e^{\frac 94}+O(\e^{\frac 15 4})=0
\end{array}
\right.
\end{equation}
Since $(\bar d_{1,\lambda},\bar d_{2,\lambda})$ satisfies \eqref{limdjfj}, it follows easily that for $\e$ sufficiently small, the system \eqref{sistemac} has only the trivial solution and the proof is now complete.
\end{proof}
We prove a technical result on the behavior of the $L^\infty$-norm of $\bar\phi_1$. After that, we prove two lemmas about the $C^0$-expansion of the reduced functional $\tilde J_\epsilon (d_1,d_2):=J_\epsilon(V_\lambda+\bar\phi_1+\bar\phi_2)$, where $\bar\phi_1 \in \mathcal{K}^\perp_1\cap B_{1,\epsilon}$ and $\bar\phi_2 \in  \mathcal{K}^\perp \cap B_{2,\epsilon}$ are the functions given by Proposition \ref{auxsolving}.

\begin{lemma}\label{techlembehvphi1}
Let $\eta$ be a small positive number and $\bar\phi_1\in\mathcal K_1^\bot$ be the solution of the first equation in \eqref{sistemaaux}, found in Proposition \ref{aux1} . Then, up to a subsequence, as $\e\to 0$, we have
$$|\bar\phi_1|_\infty \rightarrow 0, $$
uniformly with respect to $d_1$ such that $\eta<d_1<\frac{1}{\eta}.$
\end{lemma}
\begin{proof}
Let us fix a small $\eta>0$ and remember that $\tau=d_1 \e^{\frac{3}{4}}$, where $d_1\in]\eta,\frac{1}{\eta} [$ . We observe that by definition, since  $\bar\phi_1\in\mathcal K_1^\bot$ solves the first equation of (\ref{sistemaaux}), then, for all $\epsilon$ sufficiently small, for all $d_1 \in]\eta,\frac{1}{\eta} [ $, there exists a constant $c_\epsilon=c_\e(d_1)$ such that $\bar\phi_1$ weakly solves
\begin{equation}\label{eqdebolec}
-\Delta \bar\phi_1 = (\lambda_1-\epsilon) \bar\phi_1 + f(-\tau e_1 + \bar\phi_1) -\lambda_1 c_\epsilon e_1.
\end{equation}

Testing \eqref{eqdebolec} with $e_1$, and taking into account that $\bar\phi_1\in\mathcal K_1^\bot$, we deduce that $c_\epsilon \rightarrow 0$, as $\e\to 0$, uniformly with respect to $d_1 \in]\eta,\frac{1}{\eta} [$.\\
We observe that $\bar\phi_1$ is a classical solution of \eqref{eqdebolec}. This comes from standard elliptic regularity theory, the application of a well-known lemma by Brezis and Kato, taking into account that $\bar\phi_1 \in H_0^1(\Omega)$ weakly solves \eqref{eqdebolec} and the smoothness of $e_1$, f.\\

We consider the quantity $\sup_{d_1 \in ]\eta,\frac{1}{\eta}[}|\bar\phi_1|_\infty$, which is defined for all $\epsilon \in (0,\e_0)$, where $\epsilon_0>0$ is given by Proposition \ref{aux1}.
We want to prove that
\begin{equation}\label{sefunge}
 \displaystyle \lim_{\e\to0^+} \sup_{d_1 \in ]\eta,\frac{1}{\eta}[}|\bar\phi_1|_\infty=0.
\end{equation}

In order to prove \eqref{sefunge} we argue by contradiction. Assume that \eqref{sefunge} is false. Then, there exists a positive number $ m \in \mathbb{R}^+$, a sequence $(\e_k)_k \subset \mathbb{R}^+$, $\e_k\to0$ as $k\rightarrow +\infty$, such that
\begin{equation}\label{sefunge2}
 \sup_{d_1 \in ]\eta,\frac{1}{\eta}[}|\bar\phi_{1,k}|_\infty > m,
\end{equation}
for any $k\in \mathbb{N}$, where we have set $\bar\phi_{1,k}:=\bar\phi_{1}(\epsilon_k,d_1) \in B_{1,\epsilon_k}$. We observe that \eqref{sefunge2} contemplates the possibility that $ \sup_{d_1 \in ]\eta,\frac{1}{\eta}[}|\bar\phi_{1,k}|_\infty=+\infty$. 
From (\ref{sefunge2}), for any $k \in \mathbb{N}$, thanks to the definition of $\sup$, we get that there exists $d_{1,k} \in ]\eta,\frac{1}{\eta}[$ such that
$$|\bar\phi_{1,k}|_\infty (d_{1,k})> \frac{m}{2}.$$
Hence, if we consider the sequence $\left(|\bar\phi_{1,k}|_\infty (d_{1,k})\right)_k$, then, up to a subsequence, as $k\rightarrow +\infty$, there are only two possibilities:
\begin{description}
\item[(a)] $ |\bar\phi_{1,k}|_\infty (d_{1,k})\to +\infty$;
\item[(b)] $|\bar\phi_{1,k}|_\infty (d_{1,k}) \to l $, for some $l \geq  \frac{m}{2}>0$.
\end{description}
We will show that (a) and (b) cannot happen.\\\

Assume (a).
We point out that, since $\eta>0$ is fixed, then, $d_{1,k} \in ]\eta,\frac{1}{\eta}[$ for all $k$, in particular this sequence stays definitely away from $0$ and from $+\infty$. Hence, in order to simplify the notation of this proof, we omit the dependence from $d_{1,k}$ in $\bar\phi_{1,k}(d_{1,k})$, $c_{\e_k}(d_{1,k})$ and thus we simply write $\bar\phi_{1,k}$, $c_{\e_k}$. In particular, we observe that, for any fixed $k$, $\bar\phi_{1,k}$ is a function depending only on the space variable $x \in \Omega$.

Then, for any $k\in\mathbb{N}$, let $a_k\in \Omega$ such that
$|\bar\phi_{1,k}(a_k)|=|\bar\phi_{1,k}|_\infty$ and set $M_k:=|\bar\phi_{1,k}|_\infty$. We consider the rescaled function
$$\widetilde\phi_{1,k}(y):=
\frac{1}{M_k}\bar\phi_{1,k}\left(a_k +
\frac{y}{M_k^{\beta}}\right),\qquad \beta=\frac{2}{N-2},$$
defined for $y \in \widetilde
\Omega_k:=M_k^{\frac{2}{N-2}}(\Omega-a_k)$.
Moreover let us set $\widetilde e_{1,k}(y):=\frac{1}{M_k} e_1 \left(a_k  + \frac{y}{M_k^{\frac{2}{N-2}}}\right)$, $\tau_k:=d_{1,k}\e_k^{\frac{3}{4}}$. It is clear that $\|\widetilde e_{1,k}\|_{\infty,\widetilde\Omega_k} \rightarrow 0$, $\tau_k\to 0$, as $k\to+\infty$. By elementary computations we see that $\widetilde\phi_{1,k}$ solves

\begin{equation}\label{PSBN}
\left\{
\begin{array}{lr}
\displaystyle -\Delta  \widetilde\phi_{1,k} = \displaystyle \frac{\lambda_{1}-\epsilon_k}{M_k^{\frac{4}{N-2}}}  \widetilde\phi_{1,k} + f(-\tau_k \widetilde e_{1,k}+ \widetilde\phi_{1,k}) - \frac{\lambda_{1}c_{\e_k}}{M_k^{\frac{4}{N-2}}} \widetilde e_{1,k}   \qquad \mbox{in}\,\, \widetilde\Omega_k\\
\displaystyle \widetilde\phi_{1,k}=0,\qquad\qquad\qquad\qquad\ \ \qquad\qquad\qquad\qquad\qquad\qquad\qquad\mbox{on}\,\, \partial\widetilde\Omega_k
\end{array}
\right.
\end{equation}

Let us denote by $\Pi$ the limit domain of $\widetilde\Omega_k$. Since we are assuming (a) we have $M_k \rightarrow +\infty$, as $k \to +\infty$, and hence $\Pi$ is the whole $\R$ or an half-space. Moreover, since the family $(\widetilde\phi_{1,k})_k$ is
uniformly bounded and solves \eqref{PSBN}, then, by the same proof of Lemma 2.2 of \cite{Ben1}, we get that $0
\in \Pi$ (in particular $0 \notin \partial \Pi$), and, by standard
elliptic theory, it follows that, up to a subsequence,  as $k\rightarrow + \infty$, we have that
$\widetilde\phi_{1,k}$ converges in $C_{loc}^2(\Pi)$ to a
function $w$ which
satisfies

\begin{equation}\label{scplc}
-\Delta w =f(w) \ \hbox{in} \ \Pi, \ \ \ w(0)=1\ (\hbox{or} \ w(0)=-1),\ \ \ |w|\ \leq 1 \ \hbox{in}\ \Pi, \ \ \ w=0 \
\hbox{on} \ \partial \Pi.
\end{equation}

We observe that, thanks to the definition of the chosen rescaling, by elementary computations (see Lemma 2 of \cite{Iac}), it holds $\|\widetilde
\phi_{1,k}\|_{\widetilde\Omega_\epsilon}^2 =
\|\bar\phi_{1,k}\|_\Omega^2$. Now, since $\|\bar\phi_{1,k}\|\leq c \e_k^{\frac{\theta_1}{2}+\sigma}$, where $c$ depends only on $\eta$ and $\sigma$ is some positive number (see Proposition \ref{aux1}), we have $\|\widetilde \phi_{1,k}\|_{\widetilde\Omega_k}^2 = \|\bar\phi_{1,k}\|_\Omega^2\rightarrow 0$, as
$k\rightarrow +\infty$. Hence, since $\widetilde\phi_{1,k}\rightarrow w$ in $C_{loc}^2(\Pi)$, by Fatou's lemma, it follows that

\begin{equation}\label{stimaenwlim}
\|w\|_\Pi^2 \leq \liminf_{k\rightarrow +\infty} \|\widetilde \phi_{1,k}\|_{\widetilde\Omega_k}^2=0.
\end{equation}

Therefore, since $\|w\|_\Pi^2=0$ and $w$ is smooth, it follows that $w$ is constant, and from $w(0)=1$ (or $w(0)=-1$) we get that $w\equiv1$ (or $w\equiv -1$) in $\Pi$. But, since $w$ is constant and solves $-\Delta w =f(w)$ in $\Pi$, then necessarily $f(w)\equiv0$ in $\Pi$, and hence $w$ must be the null function, but this contradicts $w\equiv1$ (or $w\equiv -1$).

Assume (b). We use the same convention on the notation as in previous case.  Then $(\bar\phi_{1,k})_k$ is uniformly bounded, in particular there exist two positive constants $c_1,c_2$ such that for all $k \in \mathbb{N}$ it holds
\begin{equation}\label{limitatezzapensiero}
 c_1 < |\bar\phi_{1,k}|_\infty < c_2.
\end{equation}
By definition, $\bar\phi_{1,k}$ solves
\begin{equation}\label{eqdebolec2k}
-\Delta \bar\phi_{1,k} = (\lambda_1-\epsilon_k) \bar\phi_{1,k} + f(-\tau_k e_{1} + \bar\phi_{1,k}) -\lambda_1 c_{\epsilon_k} e_1.
\end{equation}
 Hence, by standard elliptic theory, it follows that, up to a subsequence, $\bar\phi_{1,k}$ converges in $C_{loc}^2(\Omega)$ to a function $w$ which satisfies
\begin{equation}\label{eqlimlamb1}
\begin{cases}
- \Delta w = \lambda_1 w + f(w) &\ \hbox{in}\ \Omega,\\
w=0 & \ \hbox{on}\ \partial \Omega.
\end{cases}
\end{equation}
Now, since $\|\bar\phi_{1,k}\|\leq c \e_k^{\frac{\theta_1}{2}+\sigma}$, where $c>0$ depends only on $\eta$ and $\bar\phi_{1,k} \rightarrow w$ in $C_{loc}^2(\Omega)$, then, by Fatou's Lemma and Sobolev inequality we have that $$|w|_{p+1}\leq \liminf_{k\rightarrow +\infty} |\bar\phi_{1,k}|_{p+1}=0,$$
 thus, since $w$ is smooth, it follows that $w\equiv0$. But, if $a_k \in \Omega$ is such that $|\bar\phi_{1,k}|_\infty=\bar\phi_{1,k}(a_k)$, by slightly modifications to the proof of Lemma 2.2 of \cite{Ben1} we have that $d(a_k,\partial\Omega) \nrightarrow 0$ as $k\to +\infty$. Hence, this fact, $\bar\phi_1 \rightarrow w$ in $C_{loc}^2(\Omega)$ and $w\equiv0$ contradict \eqref{limitatezzapensiero}.

 Alternatively, assuming that $\partial \Omega$ is of class $C^{2,\alpha}$, for some $\alpha \in (0,1)$, without using the arguments of Lemma 2.2 in \cite{Ben1}, but using standard elliptic regularity theory and Lemma 6.36 in \cite{GT}, since $\bar\phi_{1,k}$ is uniformly bounded, we get that, up to a subsequence $\bar\phi_{1,k}$ converges to $w$ in $C^2(\overline{\O})$, where $w$ solves \eqref{eqlimlamb1}. As before it holds $w\equiv 0$ and hence we contradicts  \eqref{limitatezzapensiero}.   The proof is then concluded.

\end{proof}

\begin{lemma}\label{lem1exp1}
For any $\eta>0$ there exists $\epsilon_0>0$ such that for any $\epsilon \in (0,\epsilon_0)$ it holds:
$$J_\lambda(V_\lambda+\bar\phi_1)=J_\lambda(V_\lambda) +O (\epsilon^{\theta_1+ \sigma}), $$
with
\begin{equation}\label{esplem1exp1}
O(\epsilon^{\theta_1 + \sigma}) = \epsilon^{\theta_1+\sigma} M_1(d_1)+  o\left( \epsilon^{\theta_2}\right),
\end{equation}
for some function $M_1$ depending only on $d_1$ (and uniformly bounded with respect to $\epsilon$), where $\theta_1, \theta_2$ are defined in \eqref{thetaj}, $\sigma$ is some positive real number (depending only on $N$).
These expansion are $C^0$-uniform with respect to $(d_1,d_2)$ satisfying condition \eqref{limdjfj}.
\end{lemma}
\begin{proof}
Let us fix $\eta>0$. By direct computation we immediately see that
\begin{equation}\label{eq1exp1}
\begin{array}{lll}
J_\lambda(V_\lambda+\bar\phi_1) - J_\lambda(V_\lambda)&=&\frac{1}{2} \int_\Omega |\nabla \bar\phi_1|^2 \ dx + \int_\Omega \nabla V_\lambda \cdot \nabla \bar\phi_1 \ dx -\frac{\lambda}{2} \int_\Omega |\bar\phi_1|^2 \ dx - \lambda \int_\Omega V_\lambda \bar\phi_1 \ dx\\[10pt]
&&- \frac{1}{p+1} \int_\Omega (|V_\lambda+\bar\phi_1|^{p+1} -|V_\lambda|^{p+1} ) \ dx.
\end{array}
\end{equation}

By definition we have $$\int_\Omega \nabla V_\lambda \cdot \nabla \bar\phi_1 \ dx=  \int_\Omega \nabla (\p\U_{\delta}- \tau e_1) \cdot \nabla \bar\phi_1 \ dx =  \int_\Omega (\U_{\delta}^p -\lambda_1 \tau e_1) \bar\phi_1 \ dx =  \int_\Omega [f(\U_{\delta}) -  \lambda_1\tau e_1] \bar\phi_1 \ dx,$$
moreover, since $F(s)=\frac{1}{p+1} |s|^{p+1}$ is a primitive of $f$, we can write (\ref{eq1exp1}) as

\begin{equation}\label{eq2exp1}
\begin{array}{lll}
J_\lambda(V_\lambda+\bar\phi_1) - J_\lambda(V_\lambda)&=& \frac{1}{2} \|\bar\phi_1\|^2   -\frac{\lambda}{2} |\bar\phi_1|_2^2  - \lambda \int_\Omega V_\lambda \bar\phi_1 \ dx +  \int_\Omega [f(\U_{\delta}) - \lambda_1 \tau e_1] \bar\phi_1 \ dx\\[10pt]
&&-  \int_\Omega [F(V_\lambda+\bar\phi_1) -F(V_\lambda)] \ dx\\[10pt]
&=& \frac{1}{2} \|\bar\phi_1\|^2   -\frac{\lambda}{2} |\bar\phi_1|_2^2  - \lambda \int_\Omega \p\U_\delta \bar\phi_1 \ dx   +(\lambda-  \lambda_1)\int_\Omega \tau e_1  \bar\phi_1 \ dx \\[10pt]
&&+ \int_\Omega [f(\U_{\delta}) -f(V_\lambda)] \bar\phi_1 \ dx-  \int_\Omega [F(V_\lambda+\bar\phi_1) -F(V_\lambda)-f(V_\lambda)\bar\phi_1] \ dx\\[10pt]
&&A+B+C+D+E+F.
\end{array}
\end{equation}
\textbf{A,B:} Thanks to Proposition \ref{auxsolving}, for all sufficiently small $\epsilon$, we have $\|\bar\phi_1\|\leq c \epsilon^{\frac{\theta_1}{2}+ \sigma}$, for some $c>0$ and for some $\sigma>0$ depending only on $N$. Hence we deduce that $A=O(\epsilon^{\theta_1+ 2 \sigma})$, $B=O(\epsilon^{\theta_1+ 2 \sigma+1})$. We point out that, since only $\bar\phi_1$ is involved in $A$ and $B$, these terms depend only on $d_1$.

\textbf{C:}  Thanks to H\"older inequality, we have
\begin{equation*}
 \lambda \int_\Omega \p\U_\delta |\bar\phi_1| \ dx \leq \lambda |\p\U_{\delta}|_{\frac{2N}{N+2}}|\bar\phi_1|_{\frac{2N}{N-2}}
\end{equation*}
For $N=5$ we have that $\displaystyle\int_{\mathbb{R}^N} \frac{1}{(1+|y|^2)^{\frac{N(N-2)}{(N+2)}}} \ dy $ is finite, so, it follows that $|\p\U_{\delta}|_{\frac{2N}{N+2}}=O(\delta^{2})$. Hence, from our choice of $\delta$ (see (\ref{deltatau5})) and since $\|\bar\phi_1\|\leq c \epsilon^{\frac{\theta_1}{2}+ \sigma}$,  we deduce that
 \begin{equation}\label{eq3exp1}
|C| \leq c \epsilon ( \epsilon^{3}  \epsilon^{\frac{3}{4}+\sigma})\leq c \e^{\theta_2+\sigma},
\end{equation}
for all sufficiently small $\epsilon$.

\textbf{D:} First we observe that this term depends only on $d_1$, and hence it will suffice that it is of order $\theta_1+\sigma$. Since $\epsilon=\lambda_1-\lambda$, and thanks to the H\"older and Sobolev inequalities, we have
$$|D| \leq \epsilon\int_\Omega \tau e_1  |\bar\phi_1| \ dx \leq c \epsilon \tau \|e_1\|_\infty \|\phi_1\|.$$
Now, thanks to the choice of $\tau$ and since $\bar\phi_1 \in B_{1,\epsilon}$ we get that
$$|D| \leq c \epsilon \epsilon^{\frac{5}{4}+\sigma} \leq c  \epsilon^{\theta_1+\sigma},$$
for all sufficiently small $\epsilon$.

\textbf{E:}
Using the definition of $V_\epsilon$ we write E as
\begin{equation}\label{eq4exp1}
\begin{array}{lll}
E&=&    \underbrace{\displaystyle\int_\Omega [f(\p\U_{\delta})-f(\p\U_\delta-\tau e_1)]\bar\phi_1\, dx}_{I_1} +\underbrace{ \displaystyle\int_\Omega [f(\U_{\delta}) - f(\p\U_{\delta})] \bar\phi_1 \ dx}_{I_2}.
\end{array}
\end{equation}
Applying H\"older and Sobolev inequalities it follows that
$$|I_1| \leq c |f(\p\U_{\delta})-f(\p\U_\delta-\tau e_1)|_{\frac{2N}{N+2}} \|\bar\phi_1\|,$$
$$|I_2| \leq c |f(\p\U_{\delta})-f(\U_\delta)|_{\frac{2N}{N+2}} \|\bar\phi_1\|.$$
By the same computations of \eqref{eq5err2}, \eqref{stimautilefuturo1}, \eqref{stimautilefuturo2} and since $\bar\phi_1 \in B_{1,\epsilon}$ we get that
$$|I_1| \leq c \epsilon^{\frac{7}{4}}\epsilon^{\frac{5}{4}+ \sigma} = c \epsilon^{\theta_2+\sigma},$$
$$|I_2| \leq c \epsilon^{\frac{9}{2}}\epsilon^{\frac{5}{4}+ \sigma} \leq c \epsilon^{\theta_2+\sigma}.$$
And we are done.\\
\textbf{F:} Applying \eqref{lem1n2e2} we get that
\begin{eqnarray*}
|F| &\leq& c \int_\Omega \left( |V_\epsilon|^{p-1} \bar\phi_1^2 + |\bar\phi_1|^{p+1}\right) \ dx \\
&\leq& c \int_\Omega( \p\U_\delta^{p-1} \bar\phi_1^2 + (\tau e_1)^{p-1} \bar\phi_1^2 + |\bar\phi_1|^{p+1}) \ dx\\
&=& F_1+F_2+F_3.
\end{eqnarray*}
We begin with estimating $F_1$.

Now, applying Lemma \ref{techlembehvphi1}, as $\epsilon \rightarrow 0$, we have $\bar\phi_1=o(1)$ in $\O$. Hence, taking into account that $\int_\O \frac{1}{|x|^4} \ dx$ is finite, we get that

\begin{eqnarray*}
\int_\O \p\U_\delta^{p-1} \bar\phi_1^2 \ dx&\leq& \int_\O \U_\delta^{p-1} \bar\phi_1^2 \ dx\leq c \int_\O \frac{\delta^2}{|x|^4}\bar\phi^2 \ dx
= o\left(\delta^2 \int_\O \frac{1}{|x|^4} \ dx\right)=o(\e^{\theta_2}).
\end{eqnarray*}


 At the end we have proved that $F_1=\int_\Omega \p\U_\delta^{p-1} \bar\phi_1^2 \ dx=o(\epsilon^{\theta_2})$, for all sufficiently small $\epsilon$.\\
For $F_2$, thanks to the definition of $\tau$ and since $\bar\phi_1 \in B_{1,\epsilon}$, we have
\begin{eqnarray*}
\int_\Omega (\tau e_1)^{p-1} \bar\phi_1^2  \ dx &\leq& \tau^{p-1} \|e_1\|_\infty^{p-1} \int_\Omega \bar\phi_1^2  \ dx
\leq c \ \tau^{p-1} \int_\Omega |\nabla \bar\phi_1|^2  \ dx \leq c_1 \epsilon^{\frac{3}{4} \frac{4}{3}} \epsilon^{2(\frac{5}{4}+ \sigma)}
\leq c_1 \epsilon^{\theta_2 + \sigma}.
\end{eqnarray*}
Finally, for $F_3$, we have
\begin{eqnarray*}
\int_\Omega  |\bar\phi_1|^{p+1} \ dx &\leq& c \| \bar\phi_1\|^{p+1}\leq  c_1 \epsilon^{\frac{10}{3}(\frac{5}{4}+\sigma)}
 \leq c_1 \epsilon^{\theta_2 + \sigma}.
\end{eqnarray*}\\

From \eqref{eq2exp1} and \textbf{A}-\textbf{F} we conclude that
$$J_\lambda(V_\lambda+\bar\phi_1) - J_\lambda(V_\lambda)= \e^{\theta_1+\sigma}M_1(d_1)+o(\e^{\theta_2}),$$
for all sufficiently small $\epsilon$, for some function $M_1$ depending only on $d_1$ (and uniformly bounded with respect to $\e$). The proof is complete.
\end{proof}

\begin{lemma}\label{lem2exp1}
For any $\eta>0$ there exists $\epsilon_0>0$ such that for any $\epsilon \in (0,\epsilon_0)$ it holds:
$$J_\lambda(V_\lambda+\bar\phi_1+\bar\phi_2)=J_\lambda(V_\lambda+\bar\phi_1) +O (\epsilon^{\theta_2+ \sigma}), $$
$C^0$-uniformly with respect to $(d_1,d_2)$ satisfying condition \eqref{limdjfj}, for some positive real number $\sigma$ depending only on $N$.
\end{lemma}

\begin{proof}
As we have seen in the proof of Lemma \ref{lem1exp1}, by direct computation we get that
\begin{equation}\label{eq1exp2}
\begin{array}{lll}
&&\hskip-1.0cm   J_\lambda(V_\lambda+\bar\phi_1+\bar\phi_2) - J_\lambda(V_\lambda+\bar\phi_1)
=\frac{1}{2} \int_\Omega |\nabla \bar\phi_2|^2 \ dx + \int_\Omega \nabla (V_\lambda+\bar\phi_1) \cdot \nabla \bar\phi_2 \ dx \\[10pt]
&&-\frac{\epsilon}{2} \int_\Omega |\bar\phi_2|^2 \ dx - \lambda \int_\Omega (V_\lambda +\bar\phi_1) \bar\phi_2 \ dx- \frac{1}{p+1} \int_\Omega (|V_\lambda+\bar\phi_1+\bar\phi_2|^{p+1} -|V_\lambda+\bar\phi_1|^{p+1} ) \ dx\\[10pt]
&=&- \frac{1}{2} \|\bar\phi_2\|^2   +\frac{\lambda}{2} |\bar\phi_2|_2^2 + \int_\Omega \nabla (V_\lambda + \bar\phi_1 + \bar\phi_2) \cdot \nabla \bar\phi_2 \ dx  \\[10pt]
&& - \lambda \int_\Omega (V_\lambda+ \bar\phi_1+\bar\phi_2) \bar\phi_2 \ dx - \int_\Omega f(V_\lambda+\bar\phi_1)\bar\phi_2 \ dx \\[10pt]
&&-  \int_\Omega [F(V_\lambda+\bar\phi_1+\bar\phi_2) -F(V_\lambda+\phi_1)-f(V_\lambda+\bar\phi_1)\bar\phi_2] \ dx\\[10pt]
\end{array}
\end{equation}
Since $\bar\phi_1+\bar\phi_2$ is a solution of (\ref{aux}) we have
$$\Pi^{\perp}\{V_\lambda+\bar\phi_1+\bar\phi_2 - i^*[\lambda (V_\lambda+\bar\phi_1+\bar\phi_2) + f(V_\lambda+\bar\phi_1+\bar\phi_2)]\}=0, $$
hence, for some $\psi \in \mathcal K$, we get that $V_\lambda+\bar\phi_1+\bar\phi_2$ weakly solves
\begin{equation}\label{eq2exp2}
 \Delta (V_\lambda+ \bar\phi_1+ \bar\phi_2) + \Delta \bar\psi -  [(\lambda_1-\epsilon)  (V_\lambda+\bar\phi_1+\bar\phi_2)+ f(V_\lambda+\bar\phi_1+\bar\phi_2)]=0.
\end{equation}
Choosing $\bar\phi_2$ as test function, since $\bar\phi_2 \in \mathcal{K}^\perp$,  $\psi \in \mathcal K$  we deduce that
\begin{equation}\label{eq3exp2}
\int_\Omega \nabla (V_\lambda+\bar\phi_1+\bar\phi_2) \cdot \nabla \bar\phi_2 \ dx - (\lambda_1-\epsilon) \int_\Omega (V_\lambda+\bar\phi_1+\bar\phi_2) \bar\phi_2 \ dx = \displaystyle \int_\Omega  f(V_\lambda+\bar\phi_1+\bar\phi_2)  \bar\phi_2 \ dx
\end{equation}
Thanks to (\ref{eq3exp2}) we rewrite (\ref{eq1exp2}) as
\begin{eqnarray}\label{eq4exp2}
\nonumber
J_\lambda(V_\lambda+\bar\phi_1+\bar\phi_2) - J_\lambda(V_\lambda+\bar\phi_1)
&=&- \frac{1}{2} \|\bar\phi_2\|^2   +\frac{\lambda_1-\epsilon}{2} |\bar\phi_2|_2^2  + \int_\Omega [f(V_\lambda+\bar\phi_1+\bar\phi_2)-f(V_\lambda+\bar\phi_1)]\bar\phi_2 \ dx \\
\nonumber
&&-  \int_\Omega [F(V_\lambda+\bar\phi_1+\bar\phi_2) -F(V_\lambda+\phi_1)-f(V_\lambda+\bar\phi_1)\bar\phi_2] \ dx\\
&=&A+B+C+D.
\end{eqnarray}

\textbf{A, B:}  Thanks to Proposition \ref{auxsolving}, for all sufficiently small $\epsilon$, we have $\|\bar\phi_2\|\leq c \epsilon^{\frac{\theta_2}{2}+ \sigma}$, for some $c>0$ and for some $\sigma>0$ depending only on $N$.  Hence we deduce that $A=O(\epsilon^{\theta_2+ 2 \sigma})$, $B=O(\epsilon^{\theta_2+ 2 \sigma})$.

 \textbf{C:} By Lemma \ref{lem1n2} and Lemma \ref{lem0n2} we get

\begin{eqnarray*}
\left|\int_\Omega [f(V_\lambda+\bar\phi_1+\bar\phi_2)-f(V_\lambda+\bar\phi_1)]\bar\phi_2 \ dx\right| &\leq & \int_{\O}| \bar\phi_2|^{p+1}\, dx +\int_{\O} |V_\lambda+\bar\phi_1|^{p+1}\bar\phi_2^2\, dx\\
&\leq & c \|\bar\phi_2\|^{p+1}+c |V_\lambda|^{p-1}_{p+1}|\bar\phi_2|_{p+1}^2+c|\bar\phi_1|^{p-1}_{p+1}|\bar\phi_2|_{p+1}^2\\&\leq& c\e^{\theta_2+\sigma}
\end{eqnarray*}
for all sufficiently small $\epsilon$.

\textbf{D:} Applying Lemma \ref{lem0n2} and H\"older inequality we get that
\begin{equation*}
\left|\int_\Omega [F(V_\lambda+\bar\phi_1+\bar\phi_2) -F(V_\lambda+\bar\phi_1)-f(V_\lambda+\bar\phi_1)\bar\phi_2] \ dx \right|
\leq c |V_\lambda|_{p+1}^{p-1}|\bar\phi_2|_{p+1}^2 + c |\bar\phi_1|_{p+1}^{p-1} |\bar\phi_2|_{p+1}^2  + c |\bar\phi_2|_{p+1}^{p+1}.
\end{equation*}

Since all the terms from $A$ to $D$ are high order terms with respect to $\epsilon^{{\theta_2}}$ the proof is complete.

\end{proof}

\section{Energy expansion of the approximate solution}
In order to prove our main results, it is important to understand critical points of the functional $$(\delta, \tau)\to J_\lambda(\p\U_\delta-\tau e_1)$$ Next we will find explicit asymptotic expressions for this functional.
\begin{proposition}\label{ridotto}
\begin{description}
The following facts hold:
\item[(i)] Let $N=4$.  For any $\eta>0$, as $\lambda\to \lambda_1^+$, the following expansion holds:
\begin{eqnarray}\label{espfunzridotto4}
\begin{array}{lll}
\displaystyle \tilde J_\lambda(s_1,s_2)&=&\displaystyle \frac{1}{4} S^{2} + \epsilon e^{-\frac 1\e}\left[ -b_1 g(s_2)^2+b_2 g(s_2) s_1-b_3 s_1^2\right] + o(\epsilon e^{-\frac 2\e}),
\end{array}
\end{eqnarray}
where $\e=\lambda-\lambda_1$, $b_1, b_2, b_3$ are positive known constants.
\item[(ii)] Let $N=5$. For any $\eta>0$,  as $\lambda\to \lambda_1^-$ it holds:
\begin{eqnarray}\label{espfunzridotto5}
\begin{array}{lll}
\displaystyle \tilde J_\lambda(d_1,d_2)&=&\displaystyle \frac{1}{5} S^{5/2} + \epsilon^{\frac 52}\left[ a_1 d_1^{2} -a_2 d_1^{\frac{10}{3}}\right] + O(\epsilon^{\frac 52 + \sigma}),
\end{array}
\end{eqnarray}
with
\begin{equation}\label{espfunzridottopart2}
O(\epsilon^{\frac 52 + \sigma}) = \epsilon^{\frac 52+\sigma} M_1(d_1)+ \epsilon^{3}\left[ a_3 d_1 d_2^{\frac 32} -a_4 d_2^2\right]
+  O\left( \epsilon^{3 + \sigma}\right),
\end{equation}
for some function $M_1$ depending only on $d_1$ (and uniformly bounded with respect to $\epsilon=\lambda_1-\lambda$), where $\sigma$ is some positive real number (depending only on $N$) and $a_j$, $j=1, 2, 3,4$ are some positive and known constants. \\
\end{description}
The expansions \eqref{espfunzridotto4}, \eqref{espfunzridotto5} and \eqref{espfunzridottopart2} are $C^0$-uniform with respect to $(s_1, s_2)$ or $(d_1,d_2)$ satisfying condition \eqref{limdjfj}.
\end{proposition}
\begin{remark}\label{remarkespfunzrid}
We point out that the term $M_1$ appearing in \eqref{espfunzridottopart2} does not depend on $d_2$ and this will be used in the sequel, in particular in \eqref{eq21propcrit}.
\end{remark}
\begin{proof}
We make some computations finding that
\begin{eqnarray*}
J_\lambda(\p\U_\delta-\tau e_1)&=& \left(\frac 12 -\frac{1}{p+1}\right)\int_\O \U_\delta^{p+1}\, dx +\frac 12 \int_\O \U_\d ^p \varphi_\d\, dx +\frac{\tau^2}{2}(\lambda_1-\lambda)\int_\O e_1^2\, dx\\
&&\hskip-2cm +\tau(\lambda-\lambda_1)\int_\O \p\U_\d e_1\, dx -\frac{\lambda}{2}\int_\O \p\U_\d^2\, dx -\frac{1}{p+1}\underbrace{\int_\O\left[|\U_\d-\varphi_\d|^{p+1}-\U_\d^{p+1}+(p+1)\U_\d^p\varphi_\d\right]\, dx}_{(I)}\\\\
&&\hskip-2cm-\frac{\tau^{p+1}}{p+1}\int_\O e_1^{p+1}\, dx +\tau\int_\O \p\U_\d^p e_1\, dx -\tau^p\int_\O \p\U_\d e_1^p\, dx \\\\
&&\hskip-2cm-\frac{1}{p+1}\underbrace{\int_\O \left[|\p\U_\d-\tau e_1|^{p+1}-\p\U_\d^{p+1}-\tau^{p+1}e_1^{p+1}+(p+1)\p\U_\d^p \tau e_1 -(p+1)\p\U_\d \tau^p e_1^p\right]\, dx}_{(II)}
\end{eqnarray*}
Then for $N=4, 5$
$$\left(\frac 12 -\frac{1}{p+1}\right)\int_\O \U_\d^{p+1}\, dx = \frac{1}{N}S_N^{N/2}+O(\d^N)$$ and $$\frac 12 \int_\O \U_\d^p\varphi_\d\, dx =O(\d^{N-2}).$$ Now if $N=4$, fixing a small $R>0$ such that $B_R\subset\subset \O$, we get
\begin{eqnarray*}
\int_\O\U_\d^2\, dx &=&\d^2\int_{|x|<R} \frac{\alpha_4^2}{(\d^2+|x|^2)^2}\, dx+\d^2\int_{\O\setminus \{|x|<R\}}\frac{\alpha_4^2}{(\d^2+|x|^2)^2}\, dx\\
 &=& \omega_4 \alpha_4^2\d^2\int_0^R\frac{r^3}{(\d^2+r^2)^2}\, dr+ O\left(\d^2\right)\\
  &=& \omega_4\alpha_4^2 \d^2 \log\frac 1\d +O(\d^2)
  \end{eqnarray*}
  where $\omega_4$ denotes the surface area of the unit sphere in $\mathbb R^4$. Instead, for $N=5$ we have
\begin{eqnarray*}
\int_\O \U_\d^2\, dx &=& \d^{-3}\int_\O \frac{\alpha_5^2}{\left(1+\left|\frac{x}{\d}\right|^2\right)^3}\, dx
=\delta^2 \int_{\mathbb R^5} \U^2\, dx  +O\left(\d^2\int_{\frac{1}{\d}}^{+\infty}\frac{r^4}{(1+r^2)^3}\, dr\right).
\end{eqnarray*}
Hence
\begin{equation*}
\int_\O \p\U_\d^2\, dx =\int_\O \U_\d^2\, dx +\int_\O \varphi_\d^2\, dx -2\int_\O \U_\d\varphi_\d\, dx \\
=\left\{\begin{array}{lr} \omega_4\alpha_4^2\d^2\log\frac{1}{\d}+O(\d^2)+O(|\varphi_\d|_\infty\int_\O \U_\d\, dx)\qquad \mbox{for}\,\, N=4\\\\
\d^2 \displaystyle \int_{\R} \U^2\, dx +O(\d^3)+O(|\varphi_\d|_2|\U_\d|_2)\qquad\qquad \mbox{for}\,\, N=5
\end{array}
\right.
\end{equation*}
and so
\begin{equation*}
\int_\O \p\U_\d^2\, dx =\left\{\begin{array}{lr} \omega_4\alpha_4^2\d^2\log\frac{1}{\d}+O(\d^2)\qquad \mbox{for}\,\, N=4\\\\
\d^2 \int_{\R} \U^2\, dx +O(\d^{\frac 52})\qquad \mbox{for}\,\, N=5.
\end{array}
\right.
\end{equation*}
Moreover

\begin{eqnarray*}
\int_\O \p\U_\d e_1\, dx &=& \int_\O e_1\left[\U_\d-\varphi_\d\right]\, dx \\
&=&\int_\O e_1\left[\alpha_N\frac{\d^{\frac{N-2}{2}}}{(\d^2+|x|^2)^{\frac{N-2}{2}}}-\alpha_N\d^{\frac{N-2}{2}}H(x, 0)+O(\d^{\frac{N+2}{2}})\right]\, dx\\
&=&\int_\O\alpha_N \d^{\frac{N-2}{2}} e_1\left[\frac{1}{|x|^{N-2}}-H(x, 0)\right]\, dx + O(\d^{\frac{N+2}{2}})\\
&=& \frac{\alpha_N}{\gamma_N} \d^{\frac{N-2}{2}}\int_\O e_1 G(x, 0)\, dx + O(\d^{\frac{N+2}{2}})\\
&=& \frac{\alpha_N}{\gamma_N\lambda_1}\d^{\frac{N-2}{2}}e_1(0)+O(\d^{\frac{N+2}{2}}),
\end{eqnarray*}
since $-\Delta e_1=\lambda_1 e_1$ and hence $e_1(0)=\lambda_1 \int_\O e_1(x)G(x, 0)\, dx$.\\\\
Moreover
\begin{eqnarray*}
\tau\int_\O e_1 \p\U_\d^p\, dx &=& \tau \int_\O e_1 \U_\d^p\, dx + \tau \int_\O e_1 (\p\U_\d^p-\U_\d^p)\, dx  \\
&=& \tau \d^{\frac{N-2}{2}}e_1(0)\int_{\R}\U^p\, dx + \begin{cases} O(\tau \d^{\frac{N+2}{2}}\log\frac{1}{\d}) & \mbox{if} \ N=4\\
 O(\tau \d^{\frac{N+2}{2}}) & \mbox{if} \ N=5\end{cases} \end{eqnarray*}
and
$$\tau^p\int_\O \p\U_\d e_1^p\, dx =\tau^p \d^{\frac{N-2}{2}}\frac{\alpha_N}{\gamma_N}\int_\O e_1^p G(x, 0)\, dx + O(\tau^p \d^{\frac{N+2}{2}})$$
Now
\begin{eqnarray*}
|(I)|&\leq & c\left( |\varphi_\d|_{p+1,\O}^{p+1}+\int_\O\U_\d^{p-1}\varphi_\d^2\, dx\right)\\
&\leq & c_1\left( \d^N +|\varphi_\d|^2_\infty\int_\O \frac{\d^{2}}{(\d^2+|x|^2)^2}\, dx\right)\\\\
&\leq & c_1\d^N + c_2\d^{N-2}\left\{\begin{array}{lr} C_0 \d^2 \log\frac{1}{\d}+O(1)\qquad\ \mbox{for}\,\, N=4\\\\\displaystyle \d^2 \int_\O \frac{1}{|x|^{N-2}} \ dx \ \ \ \qquad \mbox{for}\,\, N=5 \end{array}\right.\\\\
&\leq & c_3 \left\{\begin{array}{lr} \d^4 \log\frac{1}{\d}\qquad \mbox{for}\,\, N=4\\\\
\d^5\qquad\qquad\  \mbox{for}\,\, N=5\end{array}\right.
\end{eqnarray*}
and
\begin{eqnarray*}
|(II)|&\leq &|\int_{B_{\sqrt \d}(0)}\ldots\, dx|+|\int_{\O\setminus B_{\sqrt \d}(0)}\ldots\, dx|\\
&\leq & \int_{B_{\sqrt \d}(0)}\left||\p\U_\d-\tau e_1|^{p+1}-\p\U_\d^{p+1}+(p+1)\p\U_\d^p \tau e_1\right|\, dx\\
&&+\tau^{p+1}\int_{B_{\sqrt \d}(0)}e_1^{p+1}\, dx + \tau^p(p+1)\int_{B_{\sqrt \d}(0)} \p\U_\d e_1^p\, dx +\int_{\O\setminus B_{\sqrt \d}(0)} \p\U_\d^{p+1}\, dx  \\
&&+\tau (p+1) \int_{\O\setminus B_{\sqrt \d}(0)} \p\U_\d^p e_1\, dx+\int_{\O\setminus B_{\sqrt \d}(0)} \left||\p\U_\d-\tau e_1|^{p+1}-\tau^{p+1}e_1^{p+1}-(p+1)\tau^p e_1^p\p\U_{\d}\right|\, dx\\
&\leq & c_1\left( \tau^2\int_{B_{\sqrt \d}(0)} \p\U_\d^{p-1}e_1^2\, dx +\tau^{p+1}\int_{B_{\sqrt \d}(0)} e_1^{p+1}\, dx +\tau^p\int_{B_{\sqrt \d}(0)} \p\U_\d e_1^p\, dx \right.\\
&&\left. +\int_{\O\setminus B_{\sqrt \d}(0)} \p\U_\d^{p+1}\, dx +\tau^{p-1}\int_{\O\setminus B_{\sqrt \d}(0)} \p\U_\d^2 e_1^{p-1}\, dx +\tau\int_{\O\setminus B_{\sqrt \d}(0)} \p\U_\d^p e_1\, dx \right)\\
&\leq &c_2\left( \tau^2\d^2 |e_1|_\infty^2\int_{B_{\sqrt \d}(0)} \frac{1}{(\d^2+|x|^2)^2}\, dx +\tau^{p+1}|e_1|_\infty^{p+1}\int_0^{\sqrt \d}r^{N-1}\, dx\right. \\
&&+\tau^p |e_1|_\infty^p \d^{\frac{N-2}{2}}\int _{B_{\sqrt \d}(0)} \frac{1}{(\d^2+|x|^2)^{\frac{N-2}{2}}}\, dx+\d^N\int_{\O\setminus B_{\sqrt \d}(0)} \frac{1}{(\d^2+|x|^2)^N}\, dx \\
&&\left.+\tau^{p-1}|e_1|_\infty^{p-1}\d^{N-2}\int_{\O\setminus B_{\sqrt \d}(0)} \frac{1}{(\d^2+|x|^2)^{N-2}}\, dx +\tau|e_1|_\infty \d^{\frac{N+2}{2}}\int_{\O\setminus B_{\sqrt \d}(0)} \frac{1}{(\d^2+|x|^2)^{\frac{N+2}{2}}}\, dx\right) \\
&\leq & c_3\tau^2\d^2\left\{\begin{array}{lr}  \int_0^{\sqrt\d}\frac{r^3}{(\d^2+r^2)^2}\, dx \ \ \mbox{for}\,\, N=4\\\\ \int_0^{\sqrt \d}\, dr \qquad\qquad\mbox{for}\,\, N=5\end{array}\right. +c_4\tau^{p+1}\d^{\frac{N}{2}}+c_5\tau^p \d^{\frac{N-2}{2}}\int_{B_{\sqrt \d}(0)} \frac{1}{|x|^{N-2}}\, dx \\
&&+c_6\int_{\frac{1}{\sqrt \d}}^{+\infty}\frac{r^{N-1}}{(1+r^2)^N}\, dr +c_7\tau^{p-1}\left\{\begin{array}{lr}\int_{\O\setminus B_{\sqrt \d}(0)} \frac{\d^2}{(\d^2+|x|^2)^2}\, dx \qquad \mbox{for}\,\, N=4\\\\ \int_{\O\setminus B_{\sqrt \d}(0)} \frac{\d^3}{(\d^2+|x|^2)^3}\, dx\qquad \mbox{for}\,\, N=5\end{array}\right. \\
&&+ c_8\tau \d^{\frac{N-2}{2}}\int_{\frac{1}{\sqrt\d}}^{+\infty}\frac{r^{N-1}}{(1+r^2)^{\frac{N+2}{2}}}\, dr\\
&\leq & c_9\left( \tau^2\d^2 \left\{\begin{array}{lr} \log\frac{1}{\d}\ \ \ \mbox{for}\,\, N=4\\\\ \sqrt\d\qquad \mbox{for}\,\, N=5\end{array}\right.+\tau^{p+1}\d^{\frac N2}+\tau^p\d^{\frac N2}+\d^{\frac N2}+\tau^{p-1}\left\{\begin{array}{lr}\d^2\ \ \ \mbox{for}\,\, N=4\\\\ \d^{\frac 52}\ \ \ \mbox{for}\,\, N=5\end{array}\right. +\tau \d^{\frac N2}\right).
\end{eqnarray*}
Putting together all these estimates we get for $N=4$
\begin{equation}\label{energialimite4}
\begin{array}{lll}
\displaystyle J_\lambda(\p\U_\d-\tau e_1)&=& \displaystyle \frac 14 S^2 + O(\d^2) +\tau^2\frac{\lambda_1-\lambda}{2}\int_{\O}e_1^2\, dx+\tau(\lambda-\lambda_1)\frac{\alpha_N}{\gamma_N\lambda_1}e_1(0)\d\\
&&\displaystyle -\frac{\lambda_1}{2}\omega_4\alpha_4^2 \d^2\log\frac{1}{\d}+O\left(\d^2\log\frac{1}{\d}(\lambda-\lambda_1)\right)-\tau^{p+1}\frac{1}{p+1}\int_{\O}e_1^{p+1}\, dx\\
&&\displaystyle +\tau \d e_1(0)\int_{\mathbb R^4} \U^p\, dx + O(\tau\d^2)-\tau^p\d \frac{\alpha_N}{\gamma_N} \int_{\O}e_1^p G(x, 0)\, dx\\
&=&\displaystyle \frac 14 S^2+\e e^{-\frac 2\e}\left[-b_1 g(s_2)^2 +b_2 g(s_2)s_1-b_3 s_1^2\right]+ o(\e e^{-\frac 2\e})
\end{array}
\end{equation}
where $$b_1:=\frac 12 \int_{\O} e_1^2\, dx ;\qquad b_2:=e_1(0)\int_{\mathbb R^4}\U^p\, dx;\qquad b_3:=\frac{\lambda_1}{2}\omega_4\alpha_4^2$$
while for $N=5$ we get
\begin{equation}\label{energialimite5}
\begin{array}{lll}
\displaystyle J_\lambda(\p\U_\d-\tau e_1)&=&\displaystyle \frac 15 S^{\frac 5 2}_5 + O(\d^3) +\tau^2\frac{\lambda_1-\lambda}{2}\int_{\O}e_1^2\, dx+\tau(\lambda-\lambda_1)\frac{\alpha_N}{\gamma_N\lambda_1}e_1(0)\d^{\frac 32}\\
&&\displaystyle -\frac{\lambda_1}{2}\lambda_1 \int_{\mathbb R^5}\U^2\, dx \d^2+O\left(\d^2(\lambda-\lambda_1)\right)-\tau^{p+1}\frac{1}{p+1}\int_{\O}e_1^{p+1}\, dx\\
&&\displaystyle +\tau \d^{\frac 32} e_1(0)\int_{\mathbb R^5}\U^p\, dx + O(\tau\d^{\frac 52})-\tau^p\d^{\frac 32} \frac{\alpha_N}{\gamma_N}\int_{\O} e_1^p G(x, 0)\, dx\\
&=&\displaystyle \frac 15 S^{\frac 52}_5+\e^{\frac 52}\left[a_1 d_1^2-a_2 d_1^{\frac{10}{3}}\right]+\e^{\frac 52+\sigma}M_1(d_1)\\
&&\displaystyle +\e^3\left[a_3 d_1 d_2^{\frac 32}-a_4d_2^2\right]+ O(\e^{3+\sigma})
\end{array}
\end{equation}
and the result follows with $$a_1:=\frac 12 \int_{\O}e_1^2\, dx;\quad a_2:= \frac{1}{p+1}\int_{\O}e_1^{p+1}\, dx;\quad a_3:= e_1(0)\int_{\mathbb R^5}\U^p\, dx;\quad a_4:=\frac{\lambda_1}{2}\int_{\mathbb R^5}\U^2.$$
\end{proof}
\subsection{$C^1-$ estimate of the reduced functional in the case $N=4$}

Let $\Psi:\mathbb{R}^2_+\to \mathbb{R}$ the function defined by $$\Psi(s_1, s_2):=-b_1 g(s_2)^2+b_2 g(s_2) s_1-b_3 s_1^2,$$ where $b_j$, for $j=1,2,3$, are the positive constants appearing in \eqref{energialimite4} and $g$ is the function defined in \eqref{deltatau4}. The following result holds.
\begin{lemma}\label{espC1}
For any $\eta>0$ there exists $\e_0>0$ such that for any $\e\in (0, \e_0)$ it holds that $$\partial_{s_j}J_\lambda (V_\lambda+\bar\phi)= \e e^{-\frac 2\e}\partial_{s_j}\Psi(s_1, s_2)+ o(\e e^{-\frac 2\e})$$
$C^0$-uniformly with respect to $s_j$ satisfying \eqref{limdjfj}.
\end{lemma}
\begin{proof}
Let us fix a small $\eta>0$. By definition we have
\begin{eqnarray*}
\partial_{s_j}J_\lambda(V_\lambda+\bar\Phi_\lambda)-\e e^{-\frac 2\e}\partial_{s_j}\Psi(s_1, s_2)&=& \left(J'_\lambda (V_\lambda+\bar\Phi_\lambda)-J'_\lambda(V_\lambda)\right)[\partial_{s_j}V_\lambda]\\
&&+J_\lambda' (V_\lambda+\bar\Phi_\lambda)[\partial_{s_j}\bar\Phi_\lambda]+(\partial_{s_j}J_\lambda (V_\lambda)-\e e^{-\frac 2\e}\partial_{s_j}\Psi)
\end{eqnarray*}
Now $$(\partial_{s_j}J_\lambda (V_\lambda)-\e e^{-\frac 2\e}\partial_{s_j}\Psi)=o(\e e^{-\frac 2\e})$$
uniformly with respect to $s_j$ in compact sets of $\mathbb R_+$. Indeed,
\begin{eqnarray*}
\partial_{s_j}J_\lambda(V_\lambda)&=&J_\lambda'(V_\lambda)[\partial_{s_j}V_\lambda]=\int_\O\left(\nabla V_\lambda\nabla\partial_{s_j}V_\lambda-\lambda V_\lambda \partial_{s_j}V_\lambda-f(V_\lambda)\partial_{s_j}V_\lambda\right)\, dx\\
&=&\int_\O\nabla\p\U_\d\nabla\partial_{s_j}V_\lambda\, dx -\tau\int_\O \nabla e_1\nabla \partial_{s_j}V_\lambda\, dx -\lambda\int_\O \p\U_\d \partial_{s_j}V_\lambda\, dx \\
&&+\tau\lambda\int_\O e_1\partial_{s_j}V_\lambda\, dx -\int_\O f(V_\lambda)\partial_{s_j}V_\lambda\, dx.
\end{eqnarray*}
Now we recall that $$\partial_{s_j}V_\l =\left\{\begin{array}{lr}  \e e^{-\frac 1\e}\p Z_\d\qquad\,\, \ \  j=1\\\\-e^{-\frac 1\e}g'(s_2) e_1\quad j=2.\end{array}\right.$$
Hence for $j=1$
\begin{eqnarray*}
\partial_{s_1}J_\lambda (V_\lambda)&=& \e e^{-\frac 1\e}\left[\int_\O \nabla\p\U_\d\nabla\p Z_\d\, dx -\tau \int_\O \nabla e_1 \nabla\p Z_\d\, dx \right.\\
&&\left.-\lambda\int_\O \p\U_\d\p Z_\d\, dx +\tau\lambda \int_\O e_1 \p Z_\d\, dx -\int_\O f(V_\lambda)\p Z_\d\right]\\
&=& \e e^{-\frac 1\e}\left[\tau(\lambda-\lambda_1)\underbrace{\int_\O e_1\p Z_\d\, dx}_{(1)}-\lambda_1\underbrace{\int_\O \p \U_\d \p Z_\d\, dx}_{(2)}\right.\\
&&\left.-(\lambda-\lambda_1)\underbrace{\int_\O \p \U_\d \p Z_\d\, dx}_{(2)}-\underbrace{\int_\O\left[f(\p\U_\d-\tau e_1)-f(\p\U_\d)+f'(\p\U_\d)\tau e_1\right]\p Z_\d\, dx}_{(3)}\right.\\
&&\left.-\underbrace{\int_\O\left[f(\p\U_\d)-f(\U_\d)\right]\p Z_\d\, dx}_{(4)}+\underbrace{\int_\O f'(\U_\d)\tau e_1 \p Z_\d\, dx}_{(5)}+ \underbrace{\int_\O (f'(\p\U_\d)- f'(\U_\d))\tau e_1 \p Z_\d\, dx}_{(6)}\right]
\end{eqnarray*}
while for $j=2$ we get that
\begin{eqnarray*}
\partial_{s_2}J_\lambda (V_\lambda)&=& -g'(s_2) e^{-\frac 1\e}\left[\int_\O \nabla\p\U_\d\nabla e_1\, dx -\tau \int_\O |\nabla e_1|^2\, dx \right.\\
&&\left.-\lambda\int_\O \p\U_\d e_1\, dx +\tau\lambda \int_\O e_1^2\, dx -\int_\O f(V_\lambda) e_1\right]\\
&=& -g'(s_2) e^{-\frac 1\e}\left[\underbrace{\tau(\lambda-\lambda_1)\int_\O e_1^2\, dx}_{(I)} -\underbrace{(\lambda-\lambda_1)\int_\O \p\U_\d e_1\, dx }_{(II)} \right.\\
&&\left. -\underbrace{\int_\O \left(f(\p\U_\d-\tau e_1)-f(\p\U_\d)\right) e_1\, dx}_{(III)} -\underbrace{\int_\O \left(f(\p\U_\d)-f(\U_\d)\right)e_1\, dx}_{(IV)} -\underbrace{\int_\O f(\U_\d) e_1\, dx}_{(V)}\right]
\end{eqnarray*}
Now, by using the decompositions $\p Z_\d=Z_\d-\psi_\d$, $\p \U_\d=\U_\d - \varphi_\d$ we get that

 $$|(1)|=\left|\int_\O e_1 Z_\d\, dx -\int_\O e_1\psi_\d\, dx\right|\leq\alpha_N\int_\O \frac{e_1(x)}{|x|^{N-2}}\, dx +O(1)=O(1).$$
\begin{eqnarray*}
(2)&=& \int_\O \U_\d Z_\d +O(\d^{\frac{N-2}{2}})\\
&=& \alpha_N^2\d^{\frac{N-2}{2}}\int_\O \frac{\d^2-|x|^2}{(\d^2+|x|^2)^{N-1}}\, dx+O(\d^{\frac{N-2}{2}})\\
&=&\alpha_N^2 \d^{\frac{N-2}{2}}\int_{|x|<R}\ldots\, dx +\alpha_N\d^{\frac{N-2}{2}}\int_{\O\setminus\{|x|<R\}}\ldots+ O(\d^{\frac{N-2}{2}})\\
&=&\alpha_4^2\omega_4 \d\int_0^R \frac{r^3(\d^2-r^2)}{(\d^2+r^2)^3}\, dr+O(\d)\\
&=&\alpha_4^2\omega_4\d\ln\frac 1\d+O(\d)
\end{eqnarray*}
and by using elementary estimates we get
\begin{eqnarray*}
|(3)|&\leq & c_1\left(\int_\O |\U_\d^{p-2}\tau^2 e_1^2\p Z_\d|\, dx + \int_\O |\tau^p e_1^p\p Z_\d|\, dx \right)\\
&\leq & c_2\left(\tau^2 \int_\O |\U_\d^{p-2}Z_\d|\, dx + c\tau^2 \int_\O |\U_\d^{p-2}\psi_\d|\, dx+c\tau^p \int_\O|\p Z_\d|\, dx\right)\\
&\leq & c_3\tau^2\d\ln\frac 1\d.
\end{eqnarray*}
Moreover
\begin{eqnarray*}
|(4)|&\leq & c_1\left(\int_\O |\U_\d^{p-1}\varphi_\d \p Z_\d|\, dx + \int_\O |\varphi_\d^p\p Z_\d|\, dx+\int_\O |\U_\d^{p-2}\varphi_\d^2\p Z_\d|\, dx\right)\\
&\leq & c_2\d^2\ln\frac 1\d.
\end{eqnarray*}
Finally
\begin{eqnarray*}
(5)&=&\tau \alpha_4^p p \int_\O \frac{\d^2}{(\d^2+|x|^2)^2}e_1(x)(Z_\d-\psi_\d)\, dx\\
 &=& 3\alpha_4^p\tau e_1(0)\int_{\mathbb R^4}\frac{|y|^2-1}{(1+|y|^2)^4}\, dy +O\left(\tau \int_{\frac 1\d}^{+\infty}\frac{r^3(r^2-1)}{(1+r^2)^4}\, dr\right)+O\left(\d^2\tau\int_\O \frac{1}{(\d^2+|x|^2)^2}\, dx\right)\\
 &=&\frac 1 4\alpha_4^p \tau e_1(0) \omega_4+O(\tau\d^2\log\frac 1\d)
 \end{eqnarray*}
 since by making some easy computations one finds that
 $$ \int_{\mathbb R^4}\frac{|y|^2-1}{(1+|y|^2)^4}\, dy=\frac{\omega_4}{4}.$$
 For $(6)$, by the usual elementary inequalities and arguing as in $(4)$, we have
$$|(6)| = o(\e e^{-\frac 2\e}).$$

 Let us observe that $$b_2=e_1(0)\int_{\mathbb R^4}\U^p\, dx =\alpha_4^p e_1(0)\omega_4\int_0^{+\infty}\frac{r^3}{(1+r^2)^3}\, dr=\alpha_4^p e_1(0)\frac{\omega_4}{4}.$$
Hence
\begin{eqnarray*}
\partial_{s_1}J_\lambda(V_\lambda)&=&\e e^{-\frac 1\e}\left[e^{-\frac 1\e}\left(\underbrace{-\lambda_1 \alpha_4^2 \omega_4}_{:=2 b_3}s_1+\underbrace{\alpha_4^p e_1(0)\frac{\omega_4}{4}}_{:=b_2}g(s_2)\right)+o(e^{-\frac 1\e})\right]\\
&=& \e e^{-\frac 2\e}\partial_{s_1}\Psi+o(\e e^{-\frac 2\e}).
\end{eqnarray*}
Now
$$(I)=e^{-\frac 1\e}g(s_2) \e \int_\O e_1^2$$ $$|(II)|\leq C\e \d |e_1|_\infty=o(\e e^{-\frac 1\e})$$ $$|(III)|\leq C \left( \tau \int_\O|\p\U_\d|^{p-1}\, dx + \tau^2  \int_\O |\p\U_\d|^{p-2}\, dx + \tau^p\right)= o(\e e^{-\frac 1\e})$$ Similarly $$|(IV)|=o(\e e^{-\frac 1\e})$$ and finally $$(V)= \d\int_{\frac\O \d}\frac{\alpha_4^p}{(1+|y|^2)^p}e_1(\d y)\, dy = e_1(0)\d \int_{\mathbb R^4}\U^p\, dy + o(\d).$$
Hence
\begin{eqnarray*}
\partial_{s_2}J_\lambda(V_\lambda)&=&\e e^{-\frac 2\e}\left[-\underbrace{\int_\O e_1^2\, dx}_{:=2 b_1}g'(s_2) g(s_2)+\underbrace{e_1(0)\int_{\mathbb R^4}\U^p\, dx}_{:=b_2}g'(s_2)s_1+o(1)\right]\\
&=& \e e^{-\frac 2\e}\partial_{s_2}\Psi+o(\e e^{-\frac 2\e}).
\end{eqnarray*}

We shall prove now that
\beq\label{1}
\left[J_\lambda'(V_\lambda+\bar\Phi_\lambda)-J'_\lambda(V_\lambda)\right][\partial_{s_j}V_\lambda]=o(\e e^{-\frac 2 \e })
\eeq
and
\beq\label{2}
J_\lambda'(V_\lambda+\bar\Phi_\lambda)[\partial_{s_j}\bar\Phi_\lambda]=o(\e e^{-\frac 2 \e }).
\eeq
Let us prove \eqref{1}. We have
\begin{eqnarray*}
\left[J_\lambda'(V_\lambda+\bar\Phi_\lambda)-J'_\lambda(V_\lambda)\right][\partial_{s_j}V_\lambda]&=& \int_{\O}\nabla\bar\Phi_\lambda\nabla \partial_{s_j}V_\lambda\, dx-\l \int_{\O}\bar\Phi_\lambda\partial_{s_j}V_\l\, dx-\int_\O f'(V_\l)\partial_{s_j}V_\l \bar\Phi_\lambda\, dx\\
&&-\int_{\O}\left[f(V_\l+\bar\Phi_\lambda)-f(V_\l)-f'(V_\l)\bar\Phi_\lambda\right]\partial_{s_j}V_\l\, dx
\end{eqnarray*}

For $j=1$
$$\left|\int_{\O}\bar\Phi_\lambda\partial_{s_j}V_\l\, dx\right| =\e e^{-\frac 1\e}\left|\int_{\O}\bar\Phi_\lambda\p Z_\d\, dx\right|\leq \e e^{-\frac 1\e}|\bar\Phi_\lambda|_{4}|\p Z_\d|_{\frac 43}=o(\e e^{-\frac 2\e})$$
while for $j=2$ $$\l \int_{\O}\bar\Phi_\lambda\partial_{s_2}V_\lambda\, dx =\l e^{-\frac 1\e}g'(s_2)\int_{\O}\bar\Phi_\lambda e_1\, dx =0$$ since $\bar\Phi_\lambda\in \mathcal K^\bot$ and $$0=(\bar\Phi_\lambda, e_1)_{\zo}=\l_1\int_{\O}\bar\Phi_\lambda e_1\, dx .$$
Moreover $$\int_{\O}\nabla\bar\Phi_\lambda\nabla\partial_{s_j}V_\l dx=0$$ since $\bar\Phi_\lambda\in\mathcal K^\bot$. At the end for $j=1$
\begin{eqnarray*}
\left|\int_\O f'(V_\l)\partial_{s_1}V_\l \bar\Phi_\lambda\, dx\right|&\leq & c\left( \d \int_\O \p \U_\d^2 \p Z_\d \bar\Phi_\l\, dx\right.\\
&&\left.+\d\tau^2 \int_\O e_1^2 \p Z_\d \bar\Phi_\l\, dx +\d\tau \int_\O \p\U_\d e_1 \p Z_\d \bar\Phi_\l\, dx\right).
\end{eqnarray*}
Now
$$\d\int_\O \p\U_\d^2 \p Z_\d \bar\Phi_\l\, dx \leq \d |\p Z_\d|_\infty \int_\O \p\U_\d^2\bar\Phi_\l\, dx \leq c\frac{1}{\d}|\p\U_\d|_3^2\|\bar\Phi_\l\|\leq \d^2=o(\e e^{-\frac 2\e})$$
$$\d\tau^2 \int_\O e_1^2 \p Z_\d \bar\Phi_\l\, dx \leq C\d\tau^2 \|\p Z_\d\|\|\bar\Phi_\l\|\leq C \tau^2 \d =o(\e e^{-\frac 2\e})$$
$$\d\tau \int_\O \p\U_\d e_1 \p Z_\d \bar\Phi_\l\, dx \leq \d\tau |\p\U_\d|_3 |\p Z_\d|_3 |\bar\Phi_\l|_3\leq c \d^{\frac 43}\tau =o(\e e^{-\frac 2\e})$$
while for $j=2$
\begin{eqnarray*}
\left|\int_\O f'(V_\l)\partial_{s_2}V_\l \bar\Phi_\lambda\, dx\right|&\leq & c_1 \tau|e_1|_\infty |V_\l|_3^2\|\bar\Phi_\l\|\leq c_2\left( \tau\d^3+c\tau^2\d \right)=o(\e e^{-\frac 2\e})
\end{eqnarray*}
\begin{eqnarray*}
\left|\int_\O\left[f(V_\l+\bar\Phi_\lambda)-f(V_\l)-f'(V_\l)\bar\Phi_\l\right]\partial_{s_j}V_\l\right|&\leq &c\left( \int_\O |V_\l \bar\Phi_\lambda^2 \partial_{s_j}V_\l|\, dx + \int_\O |\bar\Phi_\lambda^3 \partial_{s_j}V_\l|\, dx\right)\\
&& \hskip-4.0cm \leq c \left\{\begin{array}{lr} \d |\p\U_\d|_4 \|\bar\Phi_\lambda\|^2 \|\p Z_\d\|+\d\tau |e_1|_4\|\bar\Phi_\lambda\|^2\|\p Z_\d\|+\d\|\bar\Phi_\lambda\|^3 \|\p Z_\d\|\qquad j=1\\\\
\tau|\p \U_\d|_4\|\bar\Phi_\lambda\|^2 |e_1|_4+\tau^2|e_1|_4^2\|\bar\Phi_\lambda\|^2+\|\bar\Phi_\lambda\|^3\tau |e_1|_4\qquad \qquad \qquad\ \  j=2\end{array}\right.\\\\
&& \hskip-4.0cm =o(\e e^{-\frac 2\e}).
\end{eqnarray*}

Now we prove \eqref{2}. Since $\Phi_\lambda$ is a solution of \eqref{aux} we have
 $$J_\l'(V_\l+\bar\Phi_\lambda)[\partial_{s_j}\bar\Phi_\lambda]=c_0(\p Z_\d, \partial_{s_j}\bar\Phi_\lambda)_{\zo}+c_1\underbrace{(e_1, \partial_{s_j}\bar\Phi_\lambda)_{\zo}}_{:=0}$$

and, from \eqref{numpag20} we deduce that
 $$(\p Z_\d, \partial_{s_j}\bar\Phi_\lambda)_{\zo}=\left\{\begin{array}{lr} O\left(\frac{\|\bar\Phi_\lambda\|}{\d}\right)\quad \mbox{for}\  j=1\\\\ 0 \qquad\qquad \ \ \  \mbox{for}\ j=2\end{array}\right.$$
So we have to analyze only the case $j=1$.\\
In order to estimate $c_0$ we write
 $$J'_\l (V_\l+\bar\Phi_\lambda)[\p Z_\d]=c_0 (\p Z_\d, \p Z_\d)_{\zo}+c_1(\p Z_\d, e_1)_{\zo}$$ 
 and
 $$J'_\l (V_\l+\bar\Phi_\lambda)[e_1]=c_0 (\p Z_\d, e_1)_{\zo}+c_1\underbrace{(e_1, e_1)_{\zo}}_{:=D_0>0}.$$ Now $$(\p Z_\d, \p Z_\d)_{\zo}=\|\p Z_\d\|^2 =A_0\d^{-2}+o(1); \qquad (\p Z_\d, e_1)_{\zo}=A_1+o(1),$$ for some positive constants $A_0, A_1$.\\
  Hence
  $$c_1:=\frac{J_\lambda'(V_\lambda+\bar\Phi_\lambda)[e_1]}{D_0}-c_0\left(\frac{A_1}{D_0}+o(1)\right)$$
  and by making some standard computations one finds that $$J_\lambda'(V_\lambda+\bar\Phi_\lambda)[e_1]=O(\d).$$
Now
\begin{eqnarray*}
J'_\l(V_\l+\bar\Phi_\lambda)[\p Z_\d] &=& -\int_\O \left[f(V_\l)-f(\U_\d)\right]\p Z_\d\, dx +\tau (\l-\l_1)\int_\O e_1 \p Z_\d\, dx\\
&&-\l\int_\O \p\U_\d \p Z_\d\, dx -\l\int_\O \bar\Phi_\lambda \p Z_\d\, dx -\int_\O \left[f(V_\l+\bar\Phi_\lambda)-f(V_\l)\right]\p Z_\d\, dx\\
&=& O\left(\sqrt{\log\frac 1 \d}\right).
\end{eqnarray*}
Hence $$c_0\left(A_0-\frac{A_1^2}{D_0}\d^2+o(\d^2)\right)=O\left(\d^2\sqrt{\log\frac 1\d}\right)$$
from which it follows that $$c_0 =O\left(\d^2 \sqrt{\log\frac 1 \d}\right)$$ and hence $$J_\l'(V_\l+\bar\Phi_\lambda)[\partial_{s_j}\bar\Phi_\lambda]=o(\e e^{-\frac 2\e}).$$

\end{proof}
\section{Proof of the main Theorems}

\begin{proof}[Proof of Theorem \ref{principale4}]
Let us fix a small $\eta>0$. Recalling that  $\e=\lambda-\lambda_1$, by Proposition  \ref{ridotto} (i), for $(s_1,s_2)$ satisfying \eqref{limdjfj} the reduced functional has the uniform expansion $$\tilde J_\lambda (s_1, s_2)= \frac 14 S^2 +\e e^{-\frac 2\e}\left[\Psi(s_1, s_2)\right]+o(\e e^{-\frac 2\e}),$$
where $$\Psi(s_1, s_2)=-b_1 g(s_2)^2+b_2 g(s_2)s_1-b_3 s_1^2.$$

It is easy to see that $\Psi$ has a non-trivial critical point in $\left(\frac{b_2}{2b_3}, 1\right)$. Moreover it is a non-degenerate maximum point if $b_2^2-4b_1b_3<0$. Hence, since the maximum points are stable under small perturbation, we get that the functional $\tilde J_\lambda(s_1, s_2)$ has a critical point in some $(\bar s_{1\lambda}, \bar s_{2\lambda})$ such that $$(\bar s_{1\lambda}, \bar s_{2\lambda})\to \left(\frac{b_2}{2b_3}, 1\right)$$ as $\lambda\to \lambda_1^+$. If instead $b_2^2-4b_1 b_3=0$, the point is a degenerate critical point but it is stable according to Definition \ref{def} since it is a maximum for $\Psi$. Indeed $$\Psi(s_1, s_2)-\Psi\left(\frac{b_2}{2b_3}, 1\right)<0\qquad \forall \,\, (s_1, s_2) \in \mathcal U$$ where $\mathcal U$ is a neighborhood of the point $\left(\frac{b_2}{2b_3}, 1\right)$, and we get the same conclusion by using also Lemma \ref{espC1}.\\\\ At the end if $b_2^2-4b_1 b_3>0$ then $\left(\frac{b_2}{2b_3}, 1\right)$ is a non degenerate critical point but we have a direction in which it is a maximum and a direction in which it is a minimum. However by Lemma \ref{espC1} we get the same conclusion.
\end{proof}

\begin{proof}[Proof of Theorem \ref{principale5}]
Let us set $G_1(d_1):= a_1 d_1^{2} -a_2 d_1^{10/3}$, where $a_1$, $a_2$  are the positive constants appearing in Proposition \ref{ridotto} (ii). 
It is elementary to see that the function $G_1:\mathbb R^+ \rightarrow \mathbb R$ has a strictly local maximum point at $\bar d_1=\left(\frac{3}{5}\frac{a_1}{a_2}\right)^{\frac{3}{4}}$.

Since $\bar d_1$ is a strictly local maximum for $G_1$, then, for any sufficiently small $\gamma>0$ there exists an open interval $I_{1,\sigma_1}$ such that $\overline I_{1,\sigma_1} \subset \mathbb R^+$, $I_{1,\sigma_1}$  has diameter $\sigma_1$, $\bar d_1 \in I_{1,\sigma_1}$ and  for all $d_1 \in \partial I_{1,\sigma_1}$
\begin{equation}\label{eq00propcrit}
G_1(d_1) \leq G_1(\bar d_1) - \gamma.
\end{equation}
Clearly as $\gamma \rightarrow 0$   we can choose $\sigma_1$ so that $\sigma_1 \rightarrow 0$.

We set $G_2(d_1,d_2):=a_3 d_1{d_2}^{\frac{3}{2}} -a_4 d_2^2$, $G_2: \mathbb
R^2_+ \rightarrow \mathbb R$, where $a_3, a_4$ are the positive constant appearing in Proposition \ref{ridotto} (ii). If we fix $d_1=\bar d_1$ then $\hat G_2(d_2):=G( \bar d_1, d_2)$ has a strictly local maximum point at $\bar d_2:=\left(\frac{3}{4}\frac{a_3}{a_4} \bar d_1\right)^{2}$. As in the previous case
there exists an open interval $I_{2,\sigma_2}$ such that $\overline I_{2,\sigma_2} \subset \mathbb R^+$, $I_{2,\sigma_2}$  has diameter $\sigma_2$, $\bar d_2 \in I_{1,\sigma_1}$ and  for all $d_2 \in \partial I_{2,\sigma_2}$
\begin{equation}\label{eq01propcrit}
\hat G_2(d_2) \leq \hat G_2(\bar d_2) - \gamma.
\end{equation}
As $\gamma \rightarrow 0$   we can choose $\sigma_2$ so that $\sigma_2 \rightarrow 0$.

Let us set $K:= \overline{ I_{1,\sigma_1} \times I_{2,\sigma_2}}$ and let $\eta>0$ be small enough so that $K\subset ]\eta,\frac{1}{\eta}[\times]\eta,\frac{1}{\eta}[$. Thanks to Proposition \ref{auxsolving}, for all sufficiently small $\epsilon$, $\tilde J_\lambda:\mathbb R_+^2 \rightarrow \mathbb R$ is defined  and it is of class $C^1$, where we recall that $\e=\lambda_1-\lambda$. By Weierstrass theorem we know there exists a global maximum point for $\tilde J_\lambda$ in $K$. Let  $(d_{1,\lambda},d_{2,\lambda})$ be that point, we want to show that there exists $\epsilon_1$ such that, for all $\epsilon < \epsilon_1$, $(d_{1,\lambda},d_{2,\lambda})$ lies in the interior of $K$.

Assume by contradiction there exists a sequence $\epsilon_n\rightarrow 0$ such that for all $n \in \mathbb N$ $$ (d_{1,\lambda_n},d_{2,\lambda_n}) \in \partial K.$$ There are only two possibilities:
\begin{description}
\item[(a)] $d_{1,\lambda_n} \in \partial I_{1,\sigma_1}$, $d_{2,\lambda_n} \in\overline I_{2,\sigma_2}$,
\item[(b)] $d_{1,\lambda_n} \in  \overline I_{1,\sigma_1}$, $d_{2,\lambda_n} \in \partial I_{2,\sigma_2}$.
\end{description}

Thanks to (ii) of Proposition \ref{ridotto} we have the uniform expansion
\begin{equation}\label{eq11propcrit}
\tilde J_\lambda (d_1, d_2) - \tilde J_\lambda ( \bar d_1,d_2)=  \epsilon^{\theta_1}\left[ G_1(d_1)-G_1(\bar d_1)\right]
+  o\left( \epsilon^{\theta_1}\right).
\end{equation}
for all $\epsilon<\epsilon_0$, $(d_1,d_2) \in K$. We point out that we have incorporated the other high order terms in $o\left( \epsilon^{\theta_1}\right)$.
Thanks to (\ref{eq00propcrit}) and (\ref{eq11propcrit}), for all sufficiently small $\epsilon$ we have
\begin{equation} \label{eq20propcrit}
 \tilde J_\lambda (d_1, d_2) - \tilde J_\lambda ( \bar d_1,d_2) <0,
\end{equation}
for all $d_1 \in \partial I_{1,\sigma_1}$, for all $d_2 \in \overline I_{2,\sigma_2}$. So for $n$ sufficiently large if (a) holds, since by definition $\tilde J_{\lambda_n} (d_{1,\lambda_n},d_{2,\lambda_n})=\max_{K} \tilde J_{\lambda_n}$, then
$$  \tilde J_{\lambda_n} (d_{1,\lambda_n},d_{2,\lambda_n}) \geq \tilde J_{\lambda_n} ( \bar d_1,d_{2,\lambda_n}),$$
which contradicts (\ref{eq20propcrit}).
Assume (b). Thanks to (ii) of Proposition \ref{ridotto} (see also Remark \ref{remarkespfunzrid}) we have the uniform expansion
\begin{equation}\label{eq21propcrit}
\tilde J_\lambda (d_1, d_2) - \tilde J_\lambda ( d_1,\bar d_2)=  \epsilon^{\theta_2}\left[ G_2(d_1,d_2)-G_2( d_1,\bar d_2)\right] + o\left(\epsilon^{\theta_2}\right),
\end{equation}
for all $\epsilon \in (0,\epsilon_0)$, for all $(d_1,d_2) \in K$.

For $n$ sufficiently large so that $\epsilon_n < \epsilon_0$ we have
\begin{equation}\label{eq22propcrit}
\begin{array}{lll}
 \displaystyle \tilde J_{\lambda_n}(d_{1,\lambda_n},d_{2,\lambda_n})-\tilde J_{\lambda_n} ( d_{1,\lambda_n},\bar d_2)&=&\displaystyle  \epsilon^{\theta_2}\left[ G_2(d_{1,\lambda_n},d_{2,\lambda_n})-G_2( d_{1,\lambda_n},\bar d_2)\right] + o\left(\epsilon^{\theta_2}\right)\\[12pt]
&=&\displaystyle  \epsilon^{\theta_2}\left[ G_2(d_{1,\lambda_n},d_{2,\lambda_n})-G_2(\bar d_1,d_{2,\lambda_n}) + G_2(\bar d_1,d_{2,\lambda_n})-G_2(\bar d_1,\bar d_{2})\right. \\[10pt]
&&\displaystyle \ \  \ \ \ \ \left.  G_2(\bar d_1,\bar d_{2}) - G_2( d_{1,\lambda_n},\bar d_2)\right] + o\left(\epsilon^{\theta_2}\right)\\[12pt]
&=&\displaystyle  \epsilon^{\theta_2}\left[  a_3  d_{2,\lambda_n}^{\frac{3}{2}} \left(d_{1,\lambda_n}-\bar d_{1}\right) + G_2(\bar d_1,d_{2,\lambda_n})-G_2(\bar d_1,\bar d_{2})\right. \\[18pt]
&&\displaystyle  \ \ \ \ \left. + \ a_3  \bar d_{2}^{\frac{3}{2}} \left(\bar d_{1}-d_{1,\lambda_n}\right) \right] + o\left(\e_n^{\theta_2}\right)
\end{array}
\end{equation}

We observe now that, up to a subsequence, $d_{1,\lambda_n} \rightarrow \bar d_1$ as $n\rightarrow + \infty$. This is a consequence of the uniform expansion given by  (ii) of Proposition \ref{ridotto}, in fact
\begin{equation}\label{eq31propcrit}
\tilde J_{\lambda_n} (d_{1,\lambda_n}, d_{2,\lambda_n}) - \tilde J_{\lambda_n} ( \bar d_1,\bar d_2)=  \epsilon_n^{\theta_1}\left[ G_1(d_{1,\lambda_n})-G_1(\bar d_1)\right]
+  o\left( \epsilon_n^{\theta_1}\right).
\end{equation}
 Since $(d_{1,\lambda_n}, d_{2,\lambda_n}) $ is the maximum point we have $\tilde J_\lambda (d_{1,\lambda_n}, d_{2,\lambda_n}) - \tilde J_\lambda ( \bar d_1,\bar d_2) \geq 0$, hence, dividing (\ref{eq31propcrit}) by $\epsilon_n^{\theta_1}$, for all sufficiently large $n$ we get that
 $ G_1(d_{1,\lambda_n})-G_1(\bar d_1) \geq - \frac{o\left( \lambda_n^{\theta_1}\right)}{\epsilon_n^{\theta_1}}$. On the other side, since $\bar d_1$ is the maximum of $G_1$, we get that
$ G_1(d_{1,\lambda_n})-G_1(\bar d_1) \leq 0$. So we have proved that
$$-\frac{o\left( \epsilon_n^{\theta_1}\right)}{\epsilon_n^{\theta_1}} \leq G_1(d_{1,\lambda_n})-G_1(\bar d_1) \leq 0, $$
and passing to the limit we deduce that $\lim_{n \rightarrow + \infty}  G_1(d_{1,\lambda_n})=G_1(\bar d_1)$. Hence, up to a subsequence, since $\bar d_1$ is a strict local minimum, the only possibility is  $d_{1,\lambda_n} \rightarrow \bar d_1$.

Since we are assuming (b), from (\ref{eq01propcrit}) we get that
$$ G_2(\bar d_1,d_{2,\lambda_n})-G_2(\bar d_1,\bar d_{2}) \leq - \gamma.$$
From this last inequality, (\ref{eq22propcrit}) and since $(d_{2,\lambda_n})_n$ is bounded, then, choosing $\bar n$ sufficiently large so that $a_3  d_{2,\epsilon}^{\frac{3}{2}} \left|d_{1,\epsilon}-\bar d_{1}\right|$ and $a_3  \bar d_{2}^{\frac{3}{2}} \left|\bar d_{1}-d_{1,\epsilon}\right|$ are small enough, we deduce that
$$ \tilde J_{\lambda_n}(d_{1,\lambda_n},d_{2,\lambda_n})-\tilde J_{\lambda_n} ( d_{1,\lambda_n},\bar d_2) < 0,$$
for all  $n> \bar n$. Since $(d_{1,\lambda_n},d_{2,\lambda_n})$ is the maximum point it also holds $$\tilde J_{\lambda_n}(d_{1,\lambda_n},d_{2,\lambda_n})-\tilde J_{\lambda_n} ( d_{1,\lambda_n},\bar d_2) \geq 0,$$
and we get a contradiction.

To complete the proof we point out that, as observed before, up to a subsequence $d_{1,\lambda} \rightarrow \bar d_1$ as $\epsilon \rightarrow 0$. With a similar argument we prove that $d_{2,\lambda} \rightarrow \bar d_2$. In fact, from the same argument of (\ref{eq22propcrit}), since $d_{1,\lambda} \rightarrow \bar d_1$ and $(d_{2,\lambda})_\epsilon$ is bounded,  we have
\begin{equation} \label{eq41propcrit}
\begin{array}{lll}
 \displaystyle 0\leq \frac{\tilde J_{\lambda}(d_{1,\lambda},d_{2,\lambda})-\tilde J_{\lambda} ( d_{1,\lambda},\bar d_2)}{ \epsilon^{\theta_2}}&=&\displaystyle  G_2(d_{1,\lambda},d_{2,\lambda})-G_2( d_{1,\lambda},\bar d_2) + \frac{o\left(\epsilon^{\theta_2}\right)}{ \epsilon^{\theta_2}}\\[12pt]
&=&\displaystyle   a_3  d_{2,\lambda}^{\frac{3}{2}} \left(d_{1,\lambda}-\bar d_{1}\right) + G_2(\bar d_1,d_{2,\lambda})-G_2(\bar d_1,\bar d_{2}) \\[18pt]
&&\displaystyle + \ a_3  \bar d_{2}^{\frac{3}{2}} \left(\bar d_{1}-d_{1,\lambda}\right) + \frac{o\left(\epsilon^{\theta_2}\right)}{ \epsilon^{\theta_2}}\\[18pt]
&=& \displaystyle o(1) + G_2(\bar d_1,d_{2,\lambda})-G_2(\bar d_1,\bar d_{2}).
\end{array}
\end{equation}
Since $\bar d_2$ is a local maximum point for $d_2 \rightarrow \hat G_2(d_2)$ we have $G_2(\bar d_1,d_{2,\epsilon})-G_2(\bar d_1,\bar d_{2})\leq 0$ and so from (\ref{eq41propcrit}) we get that
$$- o(1) \leq G_2(\bar d_1,d_{2,\lambda})-G_2(\bar d_1,\bar d_{2}) \leq 0. $$
Passing to the limit as $\epsilon \rightarrow 0$ we deduce that $ \hat G_2(d_{2,\lambda}) \rightarrow \hat G_2(\bar d_{2}) $.  Hence, up to a subsequence, since $\bar d_2$ is a strict local maximum, the only possibility is  $d_{2,\lambda} \rightarrow \bar d_2$. \\ Hence by (i) of Proposition \ref{ridotto} we have that $V_{\lambda}+\bar\phi_1+\bar\phi_2$ is a solution  of \eqref{BN}.\\

It remains to prove that the solution obtained is sign-changing. Let us set $\Phi=\Phi_\lambda:=\bar\phi_1+\bar\phi_2$.
Since $u_\lambda=V_\lambda + \Phi$ is a solution of \eqref{BN} then, by elementary computations, taking into account that by definition $-\Delta V_\lambda=\U_\delta^p - \lambda_1 \tau e_1$ (see \eqref{V}), we see that $\Phi$ solves
\begin{equation}\label{Presto}
\begin{cases}
-\Delta \Phi = \lambda \Phi + \lambda \p\U_\delta + \e \tau e_1 - \U_\delta^p +  f(u_\lambda) & \hbox{in}\ \Omega\\
\ \ \ \ \ \Phi=0 & \hbox{on}\ \partial \Omega.
\end{cases}
\end{equation}
Since $\Phi$ solves  \eqref{Presto}, then, arguing as in the proof of Lemma 3.9 in \cite{IacVair} (see also the proofs of Lemma  \ref{techlembehvphi1}, Proposition \ref{estimatelinfrem} in the present paper), we have that $|\Phi|_{\infty,\Omega}=o(\delta^{-\frac{N-2}{2}})=o(\e^{-9/4})$, \footnote{Thanks to the definition of $\delta$ and $\tau$  (see \eqref{deltatau5}) and since $d_1=d_{1,\lambda}\to \bar d_1>0$, $d_2=d_{2,\lambda}\to \bar d_2>0$, as $\e \to 0$, we have $\delta=O(\e^{3/2})$, $\tau=O(\e^{3/4})$, as $\e\to0$} for all sufficiently small $\e>0$. Hence, evaluating $u_\lambda$ at the origin, we have $u_\lambda(0)= c(N) \delta^{-\frac{N-2}{2}} - \tau e_1(0) + o(\delta^{-\frac{N-2}{2}})=c(N) d_{2,\lambda}^{-3/2} \e^{-9/4} + o(\e^{-9/4}) >0$ for all sufficiently small $\e>0$. On the contrary, thanks to Proposition \ref{estimatelinfrem}, if we fix a small ball $B_\rho$ centered at the origin, then, in $\Omega\setminus B_\rho$, we have $u_\lambda=O(\delta^{\frac{N-2}{2}}) - \tau e_1 + o(\tau) = - d_{1,\lambda} \e^{3/4} e_1 + o(\e^{3/4}) <0$, for all sufficiently small $\e>0$.  Hence $u_\lambda$ is sign-changing and the proof is complete.
\end{proof}

\begin{proposition} \label{estimatelinfrem}
Let $\Phi_\lambda$ be the remainder term appearing in Theorem \ref{principale5}. Then, for any compact subset $K$ of $\overline{\Omega}\setminus \{0\}$ we have
$$|\Phi_\lambda|_{\infty,K}=o\left((\lambda_1-\lambda)^{3/4}\right),$$
as $\lambda\to \lambda_1^-$.
\end{proposition}

\begin{proof}
 Let us set $\e:=\lambda_1-\lambda$, and let $\Phi=\Phi_\e:=\bar\phi_1+\bar\phi_2$ be the remainder term obtained in the proof of Theorem \ref{principale5}. We want to show that $|\Phi|_{\infty,K}=o(\e^{3/4})$, as $\e\to 0$. To this end, let us fix a positive number $\rho$ such that $B_{\rho}=B_{\rho}(0) \subset\subset \Omega$.

As observed in the proof of Theorem \ref{principale5} since $u_\lambda=V_\lambda + \Phi$ is a solution of \eqref{BN}, then, 
 $\Phi$ solves \eqref{Presto}.
We also point out that $\Phi$ is a smooth function since it is the difference between the two smooth functions $u_\lambda$  and $V_\lambda$. Let us set $\displaystyle\Psi= \Psi_\e:= \frac{\Phi}{\tau^{1+\gamma}}$, where $\gamma$ is a small positive number and $\tau$ is defined in \eqref{deltatau5} (see also the footnote 1). We want to prove that $|\Psi|_{\infty,\Omega\setminus B_\rho}=O(1)$, for all sufficiently small $\e>0$. By elementary computations we get that $\Psi$ solves
\begin{equation}\label{PrestoPsi}
\begin{cases}
-\Delta \Psi =\displaystyle \lambda \Psi + \lambda \frac{\p\U_\delta}{\tau^{1+\gamma}} + \frac{\e}{\tau^{1+\gamma}}  e_1 - \frac{\U_\delta^p}{\tau^{1+\gamma}} +  \tau^{p-1-\gamma}f\left(\frac{u_\lambda}{\tau}\right) & \hbox{in}\ \Omega\setminus B_\rho\\
\ \ \ \ \ \Psi=0 & \hbox{on}\ \partial \Omega,
\end{cases}
\end{equation}
We observe that in $\Omega\setminus B_\rho$ it holds $|\p\U_\delta|_{\infty, \Omega\setminus B_\rho} \leq c(N,\rho) \delta^{\frac{N-2}{2}}$, and hence, taking into account the choice of $\tau$ and $\delta$ we get that $\frac{|\p\U_\delta|_{\infty, \Omega\setminus B_\rho}}{\tau^{1+\gamma}} = o(1)$, as $\e\to 0$. By analogous computations we get that $\frac{|\U_\delta^p|_{\infty, \Omega\setminus B_\rho}}{\tau^{1+\gamma}} = o(1)$ and clearly it also holds $\frac{\e}{\tau^{1+\gamma}}\|e_1\|_{\infty,\Omega\setminus B_\rho} \leq \frac{\e}{\tau^{1+\gamma}}\|e_1\|_{\infty,\Omega}=o(1)$, as $\e\to 0$.

Let us set $M_\e:=|\Psi|_{\infty,\Omega\setminus B_{\rho}}$ and let $a_\e \in \Omega\setminus B_{\rho}$ such that $|\Psi(a_\e)| =|\Psi|_{\infty,\Omega\setminus B_{\rho}}$. Assume by contradiction that there exists a subsequence $\e_k \rightarrow 0$ (and consequently a sequence of points $a_{\e_k} \in \Omega\setminus B_{\rho}$)
 such that $M_{\e_k}= |\Psi_{\e_k}|_{\infty,\Omega\setminus B_{\rho}}=|\Psi_{\e_k}(a_{\e_k})| \to +\infty$, as $k\to + \infty$. In order to simplify the notation we shall omit the index $k$ and use the notation $\e$ to denote that subsequence. We consider the rescaled function
$$\widetilde\Psi(y):=
\frac{1}{M_\e}\Psi \left(a_\e +
\frac{y}{M_\e^{\beta}}\right),\qquad \beta=\frac{2}{N-2},$$
defined for $y \in \widetilde
\A_\e:=M_\e^{\frac{2}{N-2}}[(\Omega\setminus B_{\rho})-a_\e]$. Let us also set $\widetilde\Omega_\e:=M_\e^{\frac{2}{N-2}}(\Omega-a_\e)$
By elementary computations we see that $\widetilde\Psi$ solves

\begin{equation}\label{PrestoPsiRes}
\begin{cases}
-\Delta \widetilde \Psi =\displaystyle \lambda \frac{\widetilde\Psi}{M_\e^{2\beta}} + \lambda \frac{\p\U_\delta \left(a_\e +
\frac{y}{M_\e^{\beta}}\right)}{\tau^{1+\gamma}M_\e^{2\beta+1}} + \frac{\e}{\tau^{1+\gamma}M_\e^{2\beta+1}}  e_1 \left(a_\e +
\frac{y}{M_\e^{\beta}}\right) - \frac{\U_\delta^p \left(a_\e +
\frac{y}{M_\e^{\beta}}\right)}{\tau^{1+\gamma}M_\e^{2\beta+1}} & \\[12pt]
 \  \ \ \ \ \  \  \ \ \ \ \  \  \ \ \ \ \ \ \  +  \tau^{p-1-\gamma}f\left(\frac{u_\lambda\left(a_\e +
\frac{y}{M_\e^{\beta}}\right)}{\tau M_\e}\right ) & \hbox{in}\ \widetilde\A_\e\\
\ \ \ \ \ \widetilde\Psi=0 & \hbox{on}\ \partial \widetilde\Omega_\e,
\end{cases}
\end{equation}
As observed before, since we are assuming that $M_\e \to + \infty$, we have $$  \frac{\left|\p\U_\delta \left(a_\e +
\frac{y}{M_\e^{\beta}}\right)\right|_{\infty,  \widetilde\A_\e}}{\tau^{1+\gamma}M_\e^{2\beta+1}}=o(1),\ \  \frac{|\U_\delta^p \left(a_\e +
\frac{y}{M_\e^{\beta}}\right)|_{\infty,  \widetilde\A_\e}}{\tau^{1+\gamma}M_\e^{2\beta+1}}=o(1),  \ \frac{\e}{\tau^{1+\gamma}M_\e^{2\beta+1}}  e_1 \left(a_\e +
\frac{y}{M_\e^{\beta}}\right) =o(1),$$
as $\e\to0$. In particular, since $\widetilde \Psi$ is uniformly bounded we get that $\left| \lambda \frac{\widetilde\Psi}{M_\e^{2\beta}}\right|_{\infty,  \widetilde\A_\e}=o(1)$, and
\begin{eqnarray*}
&& \tau^{p-1-\gamma}\left|f\left(\frac{u_\lambda\left(a_\e +
\frac{y}{M_\e^{\beta}}\right)}{\tau M_\e}\right )\right|_{\infty,  \widetilde\A_\e} \\
&=&\tau^{p-1-\gamma}\left|f\left(\frac{\p\U_\delta \left(a_\e +
\frac{y}{M_\e^{\beta}}\right)}{\tau M_\e}- \frac{\tau e_1\left(a_\e +
\frac{y}{M_\e^{\beta}}\right)}{\tau M_\e} +  \tau^\gamma \widetilde \Psi \right ) \right|_{\infty,  \widetilde\A_\e}= o(1),
\end{eqnarray*}
as $\e\to 0$. Now, up to a subsequence, by standard elliptic theory $\widetilde \Psi$ converges in $C_{loc}^2(\Pi)$ to some function $\hat \Psi$ which satisfies
$- \Delta \hat\Psi=0$ in $\Pi$, where $\Pi$ is the limit domain of $ \widetilde\A_\e$. There are only three possibilities:
\begin{enumerate}
\item[(i)] $\Pi=\R$,
\item[(ii)] $\Pi$ is an half-space and $0$ lies in the interior of $\Pi$,
\item[(iii)] $\Pi$ is an half-space and $0 \in \partial \Pi$.
\end{enumerate}
We will show that (i), (ii) and (iii) bring to a contradiction.

Assume (i) or (ii). From Remark \ref{stimafixpoint1} and in Remark \ref{stimafixpoint2} we deduce that $\|\Psi\|_\Omega \to 0$ as $\e\to 0$, and hence, since $|\widetilde \Psi|_{2^*, \widetilde\A_\e} =| \Psi|_{2^*,\Omega\setminus B_\rho} \leq | \Psi|_{2^*,\Omega} \leq c \|\Psi\|_{\Omega} \to 0$, as $\e\to 0$, by Fatou's Lemma we deduce that
$$|\hat \Psi|_{2^*,\Pi}\leq \liminf_{\e\to 0 } |\widetilde\Psi|_{2^*, \widetilde\A_\e}=0.$$
Since $\hat\Psi$ is smooth, we deduce that $\hat\Psi\equiv 0$, but, since we are assuming (i) or (ii) then $0$ lies in the interior of $\Pi$, and by definition $\widetilde\Psi(0)=1$ (or $\widetilde\Psi(0)=-1$), and hence $\hat \Psi(0)=1$ (or  $\hat \Psi(0)=-1$), and we get a contradiction.

Assume (iii). Then $\partial \Pi$ is an hyperplane and $0 \in \partial \Pi$. We consider a closed ball $\overline B$ such that $\overline B \subset \overline\Pi$ and $\partial B$ is tangent at $\Pi$ in $0$. Since the limit domain of $ \widetilde\A_\e$ is $\Pi$ and thanks to the choice of $\overline B$ we get that $ \widetilde\A_\e \cap \overline B = \overline B$ for all sufficiently small $\e>0$. Since $\widetilde \Psi$ is smooth and uniformly bounded and thanks to the estimates made before, we deduce that the right-hand side of the equation in \eqref{PrestoPsiRes} is smooth (it is sufficient it is of class $C^{0,\alpha}$) and uniformly bounded. Hence, by standard elliptic theory (see Theorem 6.6 and Lemma 6.36 in \cite{GT}), we get that, up to a subsequence, the restriction of $\widetilde\Psi$ to $\overline B$ converges in $C^2(\overline B)$ to a function $\hat\Psi$. As before we have that $\hat\Psi\equiv 0$ in $B$, but, since we have the convergence in $C^2(\overline B)$, we also have $\hat\Psi(0)=1$ (or $\hat\Psi(0)=-1$) which contradicts the smoothness of $\hat\Psi$. Hence, we have that $M_\e$ is uniformly bounded and hence $|\Phi|_{\infty,\Omega\setminus B_\rho} = o(\tau)=o(\e^{3/4})$, as $\e \to 0$. The proof is complete.

\end{proof}
\begin{remark}\label{congABwN4}
We point out that, even for $N=4$, we can prove that for any compact subset $K$ of $\overline{\Omega}\setminus \{0\}$, the remainder term $\Phi_\lambda$ (appearing in Theorem \ref{principale4}) verifies $|\Phi_\lambda|_{\infty,K}=o(e^{-\frac{1}{\lambda-\lambda_1}})$, as $\lambda \to \lambda_1^+$. The key ingredient of the proof is that the remainder term verifies $\|\Phi_\lambda\|=O(\e e^{-\frac{1}{\e}})$, as $\e\to 0$ (see Proposition \ref{aux4solve}), and hence, considering, $\Psi:=\frac{\Phi_\lambda}{\e^\alpha e^{-\frac{1}{\e}}}$, where $\alpha$ is any fixed number in $(0,1)$, then, it still holds $\|\Psi\|\to 0$. Hence, arguing as in the previous proof, we get the same conclusion.
\end{remark}

\begin{remark}\label{congettura}
We believe that in the case $N=6$ the limit profile of a sign-changing solution of the problem \eqref{pb0} is given by $$u_\lambda(x)=\p\U_\d-v_\lambda(x)+\phi_\lambda$$ as $\lambda\to\bar\lambda\in (0, \lambda_1)$, where $v_\lambda$ is a positive solution of \eqref{pb0} whose existence is guaranteed by \cite{Brezis} and $\phi_\lambda$ is a remainder term which goes to zero. Moreover we have that $$\bar\lambda= 2 v_{\bar\lambda}(0)$$ and $$\lambda\to \bar\lambda^+.$$
\end{remark}

\end{document}